\documentclass[a4paper,10pt,leqno]{article}
\usepackage{amssymb}
\usepackage[mathscr]{eucal}

\def\NM{{\mathbb{N}}}
\def\RM{{\mathbb{R}}}

\def\QM{{\mathbb{Q}}}
\def\FM{{\mathbb{F}}}
\def\ZM{{\mathbb{Z}}}
\def\CM{{\mathbb{C}}}

\def\PG{{\mathfrak P}}

\def\uG{{\mathfrak u}}

\def \PG{{\mathfrak P}}

\def \aG{{\mathfrak a}}

\def \hG{{\mathfrak h}}

\def \jG{{\mathfrak j}}

\def \nG{{\mathfrak n}}

\def \uG{{\mathfrak u}}

\def\SC{{\mathcal S}}
\def\CC{{\mathcal C}}
\def\HC{{\mathcal H}}
\def\NC{{\mathcal N}}
\def\ZC{{\mathcal Z}}
\def\RC{{\mathcal R}}
\def\IC{{\mathcal I}}
\def\OC{{\mathcal O}}
\def\MC{{\mathcal M}}
\def\KC{{\mathcal K}}
\def\PC{{\mathcal P}}
\def\QC{{\mathcal Q}}
\def\LC{{\mathcal L}}
\def\EC{{\mathcal E}}
\def\FC{{\mathcal F}}
\def\VC{{\mathcal V}}
\def\GC{{\mathcal G}}
\def\UC{{\mathcal U}}
\def\JC{{\mathcal J}}

\def\ssi{si et seulement si}
\def\para{sous-groupe parabolique }
\def\paras{sous-groupes paraboliques }
\def\levi{sous-groupe de Levi }
\def\levis{sous-groupes de Levi }

\def\subsetneq{\varsubsetneq}

\def\simto{\buildrel\hbox{$\sim$}\over\longrightarrow}

\def\leq{\leqslant}
\def\geq{\geqslant}
\def\injo{\hookrightarrow}

\def\equ{\Leftrightarrow}

\def\ba{\backslash}
\def\wt{\widetilde}

\def\o#1{\overline{#1}}

\def\application#1#2#3#4#5{\begin{array}{rcl}
                            #1 \;\;\; #2 & \to &  #3 \\
                              #4 & \mapsto & #5 
                            \end{array}}

\def\To#1{\buildrel\hbox{\tiny{$#1$}}\over\longrightarrow}
\def\to{\rightarrow}

\def\ker{\mathop{\hbox{\sl ker}\,}}
\def\coker{\mathop{\hbox{\sl coker}\,}}
\def\im{\mathop{\hbox{\sl im}\,}}

\def\hom#1#2#3{\hbox{\sl Hom}_{#3}\>\!(#1,#2)}
\def\endo#1#2{\hbox{\sl End}_{#1}\>\!(#2)}
\def\aut#1#2{\hbox{\sl Aut}_{#1}(#2)}

\def\Mo#1#2{\mathop{\hbox {\sl Mod}_{#1}(#2)}}%cat des RG-modules lisses
\def\Irr#1#2{\mathop{\hbox {\sl Irr}_{#1}(#2)}}%RG-modules irreductibles
\def\cInd#1#2{\hbox {\sl ind}_{#1}^{#2}}
\def\cind#1#2#3{\hbox {\sl ind}_{#1}^{#2}\>\!(#3)} %ind a supports compacts
\def\ip#1#2#3{\hbox {\sl i}_{#1}^{#2}\>\!(#3)}  %induction parabolique
\def\ipn#1#2#3{\hbox {\sl I}_{#1}^{#2}\>\!(#3)}  %induction para normalisee
\def\Ip#1#2{\hbox {\sl i}_{#1}^{#2}}
\def\Ipn#1#2{\hbox {\sl I}_{#1}^{#2}}

\def\Res#1#2{\hbox {\sl Res}_{#1}^{#2}}              %restriction   
\def\rpn#1#2#3{\hbox {\sl R}_{#1}^{#2}\>\!(#3)}  %restriction para normalisee
\def\Rp#1#2{\hbox {\sl r}_{#1}^{#2}}  
\def\Rpn#1#2{\hbox {\sl R}_{#1}^{#2}}  
 
\def\dim{\mathop{\mbox{\sl dim}}\nolimits}

\def\gal{\mathop{\mbox{\sl Gal}}\nolimits}

\def \limi#1{\lim\limits_{\displaystyle\longrightarrow\atop {#1}}}
\def \limp#1{\lim\limits_{\displaystyle\longleftarrow\atop {#1}}}
\def \limproj{{\lim\limits_{\longleftarrow}}}
\def \limind{{\lim\limits_{\longrightarrow}}}

\makeatletter
\renewcommand{\subsection}{\@startsection{subsection}{3}{0mm}{-\baselineskip}{-0.01\baselineskip}{\bf}}%\normalfont\normalsize\itshape}}
\makeatother

\def\ali{\subsection{}}
\def\alin#1{\ali{\sl #1}\ : }

\def\ini{\setcounter{equation}{\value{subsection}}\addtocounter{subsection}{1}}

\newtheorem{theo}[subsection]{Th{\'e}or{\`e}me}
\newtheorem{lemme}[subsection]{Lemme}
\newtheorem{prop}[subsection]{Proposition}
\newtheorem{coro}[subsection]{Corollaire}

\newtheorem{fact}[subsection]{Fait}

\newtheorem{DEf}[subsection]{D{\'e}finition}

\newtheorem{rema}[subsection]{Remarque}

\newtheorem{ques}[subsection]{Question} 

\newcommand{\findem}{\hfill$\Box$\par\medskip}
\newcommand{\dem}{\noindent {\sl Preuve :} \rm }

\newenvironment{proof}{\dem}{\findem}

\usepackage[french]{babel}
\usepackage[arrow,matrix]{xy} 
\usepackage{a4wide}

\title{Finitude pour les repr\'esentations lisses de groupes $p$-adiques}
\author{J.-F. Dat}
\date{Juillet 2006}
\def\zp{{\ZM[\frac{1}{p}]}}

\def\eux{e_{U_{x}}}
\def\eubx{e_{\o{U}_{x}}}
\def\euxp{e_{U_{x}^+}}
\def\eubxp{e_{\o{U}_{x}^+}}
\def\ve{\varepsilon}
\def\la{\langle}
\def\ra{\rangle}

\begin{document}
\maketitle
\bibliographystyle{plain}

\begin{abstract}
We study basic properties of the category of smooth representations of a $p$-adic group $G$ with coefficients in any commutative ring $R$ in which $p$ is invertible. 
Our main purpose is to prove that Hecke algebras are noetherian whenever $R$ is ; a question left open since Bernstein's fundamental work \cite{bernstein} for $R=\CM$. 
In a first step, we prove that this noetherian property would follow from a generalization  of the so-called Bernstein's second adjointness property between parabolic functors for complex representations.  Then, to attack this second adjointness, we introduce and study   ``parahoric functors" between representations of groups of integral points of smooth integral models of $G$ and of their ``Levi" subgroups. Applying our general study  to Bruhat-Tits parahoric models, we get second adjointness for minimal parabolic groups.
For non-minimal parabolic subgroups, we have to restrict to classical and linear groups, and use smooth models associated with Bushnell-Kutzko and Stevens semi-simple characters.  According to recent announcements by Kim and Yu, the same strategy should also work for ``tame groups", using Yu's generic characters.

\end{abstract}

\footnotetext{Classification AMS : 20E50}

\tableofcontents

\section{Introduction}

\alin{Le probl\`eme}
Soit $K$  un corps local non archim{\'e}dien d'anneau des entiers $\OC$ et
  de corps r{\'e}siduel $k$ de caract{\'e}ristique $p$, et $\GC$ un groupe alg\'ebrique r\'eductif connexe d\'efini sur $K$. On pose\footnote{De mani\`ere g\'en\'erale, nous noterons les $K$-sch{\'e}mas par des lettres  calligraphi{\'e}es et leurs $K$-points par les lettres ordinaires correspondantes.
} $G:=\GC(K)$. 
Si $H$ est un sous-groupe ouvert compact de $G$, on peut former l'anneau de Hecke $\HC(G,H):=\ZM[H\ba G/ H]$ des doubles classes de $H$ dans $G$.
Un c\'el\`ebre th\'eor\`eme de Bernstein \cite{bernstein} affirme qu'apr\`es extension des scalaires de $\ZM$ \`a $\CM$, ces anneaux sont noeth\'eriens. La preuve de Bernstein est tr\`es indirecte car il \'etudie principalement la cat\'egorie $\Mo{\CM}{G}$ des repr\'esentations complexes lisses de $G$. Le principe fondamental qui soutient toute sa th\'eorie est qu'une
repr\'esentation irr\'eductible supercuspidale complexe est un objet projectif (``modulo le centre") de $\Mo{\CM}{G}$. Ainsi, la m\^eme preuve fournirait le m\^eme r\'esultat apr\`es extension des scalaires \`a n'importe quel corps de caract\'eristique suffisamment grande (banale dans la terminologie de Vign\'eras \cite{Vigsheaves}), mais ne fonctionne pas, ne serait-ce que pour un corps de caract\'eristique non-banale, sans parler du cas d'un anneau.
Il est pourtant naturel de s'attendre \`a ce que ces anneaux de Hecke v\'erifient une propri\'et\'e de type ``th\'eor\`eme de Hilbert" : $R$ noeth\'erien $\Rightarrow R[H\ba G/ H]$ noeth\'erien. C'est ce que nous \'etudions -- entre autres -- dans cet article, avec l'hypoth\`ese suppl\'ementaire que $p$ est inversible dans $R$, car notre m\'ethode aussi passe par l'\'etude de la cat\'egorie $\Mo{R}{G}$ des $RG$-modules lisses et que pour faire le lien avec les alg\`ebres de Hecke, on a besoin d'une mesure de Haar.

Bien-s\^ur, l'int\'er\^et de passer par la cat\'egorie $\Mo{R}{G}$ vient des foncteurs d'induction et restriction paraboliques qui permettent de faire des raisonnements par r\'ecurrence sur le rang semi-simple de $G$. Si $\PC$ est un $K$-sous-groupe parabolique de $\GC$ et $\MC$ un  \levi de $\PC$, la r\'eciprocit\'e de Frobenius nous dit que la ``restriction" parabolique $\Rp{G,P}{M}$ est adjointe \`a gauche de l'induction $\Ip{M,P}{G}$.
Dans un article non publi\'e \cite{bernunp} mais bien connu des sp\'ecialistes, Bernstein a d\'ecouvert (avec surprise) que pour les repr\'esentations {\em complexes}, il existe une deuxi\`eme propri\'et\'e d'adjonction, entre le foncteur $\Ip{M,P}{G}$ \`a gauche et le foncteur $\delta_P \Rp{G,\o{P}}{M}$ \`a droite, o\`u $\o{P}$ est le parabolique oppos\'e \`a $P$ par rapport \`a $M$ et $\delta_P$ le module de $P$. Bushnell a publi\'e \cite{Bu} une preuve diff\'erente de cette deuxi\`eme adjonction, mais chacune de ces preuves repose de mani\`ere cruciale sur la propri\'et\'e de noeth\'eriannit\'e des alg\`ebres de Hecke complexes.

\alin{Principaux th\'eor\`emes}
Dans le pr\'esent article nous proc\'edons dans l'autre sens. Dans la section \ref{noether} nous prouverons en effet, pour un anneau de coefficients $R$ noeth\'erien et o\`u $p$ est inversible, 
que {\em la deuxi\`eme adjonction implique la noeth\'eriannit\'e } :
\begin{theo} \label{theo1}
Si pour tout \para de tout \levi de $\GC$, les foncteurs paraboliques v\'erifient la seconde adjonction, alors la cat\'egorie $\Mo{R}{G}$ est noeth\'erienne, ainsi que les alg\`ebres de Hecke $\HC_R(G,H)$ pour tout pro-$p$-sous-groupe ouvert $H$ de $G$.
\end{theo}
La preuve de ce r\'esultat utilise les propri\'et\'es des repr\'esentations \`a coefficients dans des corps valu\'es \'etudi\'ees dans \cite{nutempere}. 

D\`es lors, la majeure partie de cet article vise \`a \'etablir la propri\'et\'e de seconde adjonction. Pour cela on utilise un nouvel outil baptis\'e {\em induction parahorique}.
Identifions l'immeuble \'etendu $B(\MC,K)$ de $\MC$ \`a un sous-$M$-espace de celui de $\GC$.
Si $x\in B(\MC,K)$, on d\'efinit un certain idempotent $\ve_{x,P}$ normalis\'e par $M_x$, dans l'alg\`ebre $\zp G_x$ des $\zp$-distributions sur $G_x$. Par 
produit tensoriel avec le $(M_x,G_x)$-bimodule $\ve_{x,P}\CC^\infty_R(G_x)$, on obtient des
foncteurs d'induction $I_{x,P}$ et restriction
$R_{x,P}$ entre $M_x$-modules lisses et
$G_x$-modules lisses {\`a} coefficients dans $R$. 
La motivation initiale pour introduire ces foncteurs vient 
des relations de commutation remarquables suivantes, dans le cas o\`u $P$ est {\em minimal}  ({\em cf } \ref{pfcommutmin} ii) et \ref{commut}) :
\ini
\begin{equation} \label{commutmin} \cInd{G_x}{G} \circ I_{x,P} \simeq \Ip{M,P}{G} \circ
\cInd{M_x}{M} \;\; \hbox{ et } \;\; \Res{M}{M_x}\circ \Rp{G,P}{M} \simeq
R_{x,P}\circ \Res{G}{G_x}.
\end{equation}
%\end{enumerate}
Notons que ceci est nouveau aussi pour $R=\CM$. Nous expliquons alors en \ref{cororescent} comment la seconde adjonction (pour les paraboliques minimaux) d\'ecoule de ces formules.

Malheureusement, il semble que ces relations de commutation entre foncteurs paraboliques et parahoriques soient sp\'ecifiques aux paraboliques minimaux. En g\'en\'eral on introduit la notion d'idempotent $P$-bon de $RM$, {\em cf} \ref{pbon}. Il s'agit grosso-modo d'un idempotent $\ve$ pour lequel il existe un $x\in B(\MC,K)$ tel que $\ve\in RM_x$ et
$$ \ve.\Res{M}{M_x}\circ \Rp{G,P}{M} \simeq \ve. R_{x,P}\circ \Res{G}{G_x} .$$
S'il existe ``suffisamment" de tels idempotents ({\em i.e.} s'ils engendrent la cat\'egorie $\Mo{R}{M}$), alors d'apr\`es \ref{cororescent}, la seconde adjonction est v\'erifi\'ee pour $P$.

Reste donc \`a produire des familles g\'en\'eratrices d'idempotents $P$-bons. Ceci semble une t\^ache ardue en g\'en\'eral, aussi difficile que la construction de strates/types ``raffin\'e(e)s". Pour un groupe lin\'eaire ou classique, nous montrons dans les sections \ref{lineaire} et \ref{classique} que la famille des idempotents associ\'es aux caract\`eres semi-simples de Stevens \cite{Stevens} convient. Nous devons au passage prouver des r\'esultats apparemment nouveaux dans cette th\'eorie (propositions \ref{p2} et \ref{p2c}). On obtient donc sur tout anneau de coefficients $R$ o\`u $p$ est inversible :

\begin{theo} \label{theo2}
Soit $G$ un groupe lin\'eaire, classique (on suppose alors $p\neq 2$), ou de rang relatif $1$, alors pour tout parabolique de tout sous-groupe de Levi, les foncteurs paraboliques associ\'es v\'erifient la seconde adjonction.  
\end{theo}

Outre la propri\'et\'e de noeth\'eriannit\'e \ref{theo1}, on a aussi les cons\'equences suivantes :

\begin{coro} Avec la m\^eme hypoth\`ese sur $G$,
\begin{enumerate}
        \item Pour $R$ noeth\'erien, les  foncteurs de restriction parabolique de Jacquet pr\'eservent la $R$-admissibilit\'e, {\em cf} \ref{cororescent} i).
  \item {\em Support uniforme} : Pour tout pro-$p$-sous-groupe $H$ de
    $G$ il existe un sous-ensemble  
    $S_H$ de $G$, compact modulo le centre et {\em ind{\'e}pendant de
      l'anneau $R$} supportant toutes les fonctions
    cuspidales dans $\CC^c_R(H\ba G/H)$, {\em cf} \ref{support}.
\item {\em Irr{\'e}ductibilit{\'e} g{\'e}n{\'e}rique} : Si $R$ est un corps
  alg{\'e}briquement clos, alors pour tout \para $\PC=\MC\UC$ et toute $\pi \in
  \Irr{R}{M}$, la famille $\ip{M,P}{G}{\pi\psi}$ pour $\psi:M/M^c \To{}
  R^\times$ est g{\'e}n{\'e}riquement irr{\'e}ductible, {\em cf} \ref{nudiscus}. 
 \end{enumerate}
\end{coro}
 Signalons que le dernier point pour $R$ de caract\'eristique positive n'est nouveau que lorsque $K$ est aussi de caract\'eristique positive, {\em cf} \cite{nutempere}. Quant au point ii), il peut \^etre utile \`a ceux qui s'int\'eressent aux congruences entre formes automorphes.

L'ingr\'edient essentiel pour produire des idempotents $P$-bons est une sorte de g\'en\'eralisation d'un r\'esultat de Howlett-Lehrer \cite{HL} o\`u l'on remplace les foncteurs paraboliques des groupes de Lie finis par nos foncteurs parahoriques pour des mod\`eles entiers de $\GC$, {\em cf} partie \ref{mod}, qui dans les applications seront les mod\`eles de Bruhat-Tits (cas minimal) ou les mod\`eles entiers associ\'es aux types de Bushnell-Kutzko-Stevens (cas g\'en\'eral pour les groupes classiques). On peut aussi appliquer cet ingr\'edient aux mod\`eles entiers de Yu \cite{Yumodel} associ\'es \`a ses types ``mod\'er\'es" \cite{Yumodel}, voir la partie \ref{modere}. En utilisant un r\'esultat d'exhaustivit\'e (analogue de \ref{p3} et \ref{p3c}) annonc\'e r\'ecemment par Yu et Kim pour les groupes {\em mod\'er\'es} (ceux dont tout tore se d\'eploie sur une extension mod\'er\'ee), on doit pouvoir prouver le th\'eor\`eme \ref{theo2} pour de tels groupes, avec formellement la m\^eme preuve que pour les groupes classiques.

Enfin, on obtient aussi quelques r\'esultats partiels sans conditions sur $\GC$ ; outre la seconde adjonction pour les paraboliques minimaux d\'eja mentionn\'ee, on prouve la noeth\'eriannit\'e de la sous-cat\'egorie pleine de $\Mo{R}{G}$ des objets ``de niveau z\'ero", ainsi que la seconde adjonction des foncteurs paraboliques restreints \`a ces sous-cat\'egories, {\em cf} \ref{niveau0}.

\alin{Organisation de l'article} 
L'induction parahorique est d\'efinie dans la premi\`ere section. C'est un cas particulier d'"induction" pour des groupes munis d'une d\'ecomposition d'Iwahori ``abstraite", et c'est par une discussion de cette situation g\'en\'erale que la section commence ; le r\'esultat fondamental est la proposition \ref{propIwa}.

La partie 3 \'etudie les propri\'et\'es de commutation entre foncteurs parahoriques et foncteurs paraboliques (notamment \ref{commut}). On y d\'egage la notion d'idempotent $P$-bon, et on prouve  la seconde adjonction en \ref{cororescent}, sous r\'eserve qu'il existe ``suffisamment" de tels idempotents.
Puis la partie 4 prouve que ``la seconde adjonction implique la noetheriannit\'e".

Dans la partie 5 on consid\`ere des mod\`eles entiers lisses de $\GC$ et on \'etend implicitement la notion d'induction parahorique \`a ces groupes. On prouve alors un \'enonc\'e  d'ind\'ependence du sous-groupe parahorique \ref{theomod}, analogue \`a l'\'enonc\'e principal de \cite{HL} sur l'ind\'ependance du sous-groupe parabolique pour l'induction parabolique des groupes finis. 

Dans la partie 6, on sp\'ecialise la pr\'ec\'edente aux mod\`eles de Bruhat-Tits et on applique les r\'esultats obtenus aux foncteurs paraboliques minimaux et aux repr\'esentations de niveau z\'ero. Dans les parties \ref{lineaire} et \ref{classique}, on sp\'ecialise la partie \ref{mod} aux mod\`eles entiers associ\'es aux caract\`eres semi-simples de Stevens, et on prouve le th\'eor\`eme \ref{theo2}. Dans la partie \ref{modere} on sp\'ecialise la partie \ref{mod} aux mod\`eles entiers de Yu et \`a ses caract\`eres g\'en\'eriques ; il ne manque qu'un r\'esultat crucial d'exhaustivit\'e pour en d\'eduire la noetheriannit\'e.

\medskip

\noindent{\em Remerciements :} je remercie B. Lemaire pour quelques discussions sur les strates et autres mod\`eles entiers, S. Stevens pour quelques explications sur sa th\'eorie, et G. Henniart pour son int\'er\^et dans ce travail. Je remercie aussi M.-F. Vign\'eras pour m'avoir transmis ce probl\`eme de  noetheriannit\'e, qu'elle a abord\'e dans \cite{vigtypmod}. Signalons aussi, outre les travaux non-publi\'es de Bernstein sur ce sujet (\`a peine \'evoqu\'es dans \cite[5.4, Rk 1]{bernunp}) une autre approche imagin\'ee par Bezrukavnikov, tr\`es naturelle mais qui \`a ma connaissance n'a pas abouti, consistant \`a essayer de prouver la noetheriannit\'e du gradu\'e des alg\`ebres de Hecke pour une certaine filtration ``g\'eom\'etrique" \cite[II-2]{Bez}.

\alin{Notations}
 Pour tout groupe localement compact totalement discontinu 
   $H$ et tout anneau commutatif unitaire $R$ on
  note : 
 \begin{itemize} 
 \item $\CC^{\infty,c}_R(H)$, le $R$-module des fonctions localement
 constantes {\`a} valeurs dans $R$ et {\`a} support compact. 
 \item  $RH$ le $R$-module des distributions \`a support compact 
 et \`a valeurs dans $R$.  Le produit de convolution en fait une $R$-alg\`ebre unitaire qui se plonge dans le commutant $\endo{RH}{\CC^{\infty,c}_R(H)}$ des translations \`a droite par $H$.
Lorsque $H$ est compact, ce plongement est un isomorphisme, et si $(H_n)_{n\in\NM}$ est une base de voisinages de l'unit\'e form\'ee de sous-groupes normaux, on a un isomorphisme de $R$-alg\`ebres $RH\simto \limproj\, R[H/H_n]$.
  
 \item  $\Mo{R}{H}$ la cat{\'e}gorie des repr{\'e}sentations lisses de $H$ {\`a}
  coefficients dans $R$. Tout objet de $\Mo{R}{H}$ est canoniquement muni d'une action de $RH$, ce qui donne un plongement pleinement fid\`ele $\Mo{R}{H} \injo \Mo{}{RH}$.
 S'il existe une  mesure de Haar $dh$ sur $H$ \`a valeurs dans $R$, alors le sous-$R$-module
 $\HC_R(H):=\CC^{\infty,c}_R(H)dh$ de $RH$ form\'e des distributions  localement  constantes est  un id\'eal bilat\`ere de $RH$ engendr\'e par ses idempotents, et 
l'image essentielle de $\Mo{R}{H}$ dans $\Mo{}{RH}$ s'identifie \`a la cat\'egorie des modules ``unitaux" sur $\HC_R(H)$. 
\end{itemize}  

Si $K$ est un sous-groupe compact de $H$ de pro-ordre inversible dans $R$, la mesure de Haar $dk$ de volume total $1$ sur $K$ d\'efinit un idempotent de $RH$ que nous noterons $e_K$. Son action sur un objet $(\pi,V)$ de $\Mo{R}{G}$ est 
%$\CC^{\infty,c}_R(H)$ est 
donc donn\'ee par $e_K*v = \int_K \pi(k)v dk$.

\section{D{\'e}compositions {\`a} la Iwahori et induction parahorique} \label{pfiwa}

Au d\'ebut de cette section, la lettre $G$ ne d{\'e}signe 
pas un groupe r\'eductif $p$-adique.

\begin{DEf} \label{defIwa}
Soit $G$ un groupe profini, muni de deux sous-groupes ferm\'es $U$ et $\o{U}$ normalis\'es par un troisi\`eme sous-groupe ferm\'e $M$. Nous dirons que le triplet $(U,M,\o{U})$ induit une {\em d\'ecomposition d'Iwahori} de $G$ si :
\begin{enumerate}
        \item L'application produit $U\times M\times \o{U} \To{} G$ est bijective.
        \item Il existe une base $(G_i)_{i\in\NM}$ de voisinages de $G$ form\'ee de sous-groupes ouverts normaux de la forme $G_i=(U\cap G_i)(M\cap G_i)(\o{U}\cap G_i)$. 
\end{enumerate}
 
\end{DEf}

On obtient par exemple de telles d\'ecompositions lorsque $\o{U}$, $U$, $M$ et
$G$ sont les points entiers de groupes alg{\'e}briques affines lisses $\o{\UC}$, $\UC$, $\MC$ et $\GC$
d{\'e}finis sur un anneau complet de valuation discr\`ete \`a corps r\'esiduel fini tels que l'application produit $\o{\UC}\times \MC \times \UC \To{} \GC$ est une immersion ouverte induisant un isomorphisme des fibres sp\'eciales (voir aussi \ref{dilatIwa}).

\begin{prop} \label{propIwa}
Soit $G=UM\o{U}$ un groupe profini muni d'une d\'ecomposition d'Iwahori comme ci-dessus. On suppose que $M$ contient un sous-groupe ouvert normal $M^\dag$ tel que l'ensemble $G^\dag:=UM^\dag\o{U}$ soit un pro-$p$ sous-groupe (ouvert) de $G$.
Alors il existe une unique distribution centrale inversible $z_{U,\o{U}}\in \ZC({\zp}M)^\times$ telle que
  la distribution  $z_{U,\o{U}}^{-1} e_Ue_{\o{U}}$ soit un idempotent de l'anneau  %$=z_M.e_Ue_{\o{U}}e_Ue_{\o{U}}$ dans 
  $ \zp G$.
\end{prop}

\begin{proof}
Comme premi\`ere cons\'equence de la d\'ecomposition d'Iwahori, l'application
\ini\begin{equation} \label{eqIwa}
\application{}{ \zp M}{e_U \zp G e_{\o{U}}}{f}{e_Ufe_{\o{U}}}
\end{equation}
est un isomorphisme de $\zp M$-modules. Elle est en effet injective par l'axiome i) de \ref{defIwa}, et lorsque $G$ est fini elle est surjective, puisqu'alors la multiplication induit un isomorphisme de $\zp$ modules $\zp[\o{U}]\otimes \zp[M] \otimes \zp[U]\simto \zp[G]$. Pour $G$ profini, on se ram\`ene au cas fini gr\^ace \`a  l'axiome ii) de \ref{defIwa} qui donne une pr\'esentation de $G$ sous la forme 
$ G \simeq \limproj G/G_i$ o{\`u} les $G/G_i$ sont des groupes finis munis de d\'ecompositions d'Iwahori $G/G_i=U/(U\cap G_i).M/(M\cap G_i).\o{U}/(\o{U}\cap G_i)$. L'isomorphisme \ref{eqIwa} est alors la limite projective des isomorphismes correspondants pour les $G/G_i$.

%o\`u $U_i:=G_i\cap U$, $M_i:=G_i\cap M$ et $\o{U}_i:=G_i\cap\o{U}$.

Par \ref{eqIwa}, il existe un unique \'el\'ement $z_{U,\o{U}}\in \zp M$ tel que $$e_Ue_{\o{U}}e_Ue_{\o{U}}= e_U z_{U,\o{U}} e_{\o{U}}=z_{U,\o{U}} e_U e_{\o{U}}.$$
Par unicit\'e et puisque $M$ normalise $U$ et $\o{U}$, cet \'el\'ement est {\em central} dans $\zp M$. Reste \`a prouver qu'il est inversible. Par application de ce qui pr\'ec\`ede \`a $G^\dag$, on constate que $z_{U,\o{U}}\in RM^\dag$, et ceci nous ram\`ene au cas o\`u $G$ est pro-$p$. Il nous suffit alors de prouver que pour tout $i$ la distribution lisse $z_{U,\o{U}}*e_{M\cap G_i}$ est un \'el\'ement inversible dans l'anneau $\zp[M/M\cap G_i]$, et ceci nous ram\`ene au cas o\`u $G$ est {\em fini}.

Supposons dor\'enavant que $G$ est un $p$-groupe fini. Pour montrer que $z_{U,\o{U}}$ est inversible, il suffit de prouver que la   multiplication par $z_{U,\o{U}}$ dans le
$\zp$-module libre de type fini $\zp[M]$ a un d{\'e}terminant inversible, car alors elle
 sera surjective et son image contiendra  l'unit{\'e}. Il suffit donc de montrer que pour tout corps $R$ de caract\'eristique diff\'erente de $p$, l'image $z_{U,\o{U}}^R$ de $z_{U,\o{U}}$ dans $R[M]$ est inversible.

Pour un tel corps $R$,  les cat\'egories $\Mo{R}{M}$ et $\Mo{R}{G}$ sont semi-simples. 
Notons $\Ipn{M}{G}$ le foncteur 
$$\application{\Ipn{M}{G}:\;}{\Mo{R}{M}}{\Mo{R}{G}}
{W}{R[G]e_Ue_{\o{U}}\otimes_{R[M]} W}$$
o{\`u} l'on consid\`ere $R[G]e_Ue_{\o{U}}$ comme  $R[M]$-module {\`a} droite
et $R[G]$-module {\`a} gauche par la formule  $(g,m).f:=gfm$.
L'isomorphisme \ref{eqIwa} induit l'isomorphisme
$$\application{}{ R[M]}{e_U R[G]e_Ue_{\o{U}}}{f}{e_Ufe_Ue_{\o{U}}}$$
puisque $e_U f=e_U f e_U$. %Il en va {\'e}videmment de m{\^e}me si on remplace $\zp$ par $R$. 
On a donc $e_U\ipn{M}{G}{W}\simeq_{M} W$   pour tout $W\in \Mo{R}{M}$. Par ailleurs, 
puisque $$R[G]e_UR[G]e_Ue_{\o{U}}=R[G]e_Ue_{\o{U}},$$ le sous $R$-module
$e_U\ipn{M}{G}{W}$ engendre  $\ipn{M}{G}{W}$ en tant que $RG$-module. 
Par semi-simplicit{\'e} de $\Mo{R}{G}$, on en d\'eduit que $\Ipn{M}{G}$ envoie irr\'eductibles sur irr\'eductibles.

Maintenant, nous pr\'etendons que l'application
$$\application{}{ R[M]}{e_{\o{U}} R[G]e_Ue_{\o{U}}}{f}{e_{\o{U}}fe_Ue_{\o{U}}}$$
est aussi un isomorphisme. En effet, on en construit un inverse comme ceci : la restriction des fonctions $R[G] \To{} R[M\o{U}]$ est une application $M\o{U}$-\'equivariante \`a droite et \`a gauche qui induit une application $M$-\'equivariante \`a droite et \`a gauche $e_{\o{U}}R[G] e_{\o{U}} \To{} e_{\o{U}}R[M\o{U}]e_{\o{U}}\simeq R[M]$. Celle-ci induit l'inverse cherch\'ee.

Il s'ensuit que pour tout $W\in\Mo{R}{G}$, on a $e_{\o{U}}\ipn{M}{G}{W}\simeq_M W$.
 Supposons de plus $W$  irr{\'e}ductible, alors le $R$-module 
 $e_{\o{U}}\ipn{M}{G}{W}$ \'etant non-nul, il engendre le $RG$-module irr\'eductible $\ipn{M}{G}{W}$. 
  Ainsi l'inclusion $R[G] e_{\o{U}}e_Ue_{\o{U}} \subset R[G]e_Ue_{\o{U}}$ induit 
 pour tout $W\in \Irr{R}{M}$ (qui rappelons-le est $R[M]$-projectif) un isomorphisme :
$$ R[G] e_{\o{U}}e_Ue_{\o{U}} \otimes_{ R[M]} W \simto   R[G]
e_Ue_{\o{U}}\otimes_{ R[M]} W.$$
Par semi-simplicit{\'e} de $\Mo{R}{M}$, cet isomorphisme est valable pour tout $W\in
\Mo{R}{M}$ et en particulier pour $ R[M]$, ce qui nous fournit
l'{\'e}galit{\'e}  $$ R[G] e_{\o{U}}e_Ue_{\o{U}} = R[G]e_Ue_{\o{U}}.$$ 
En
particulier, il existe $f\in  R[G]$ telle que
$fe_{\o{U}}e_Ue_{\o{U}}=e_Ue_{\o{U}}$. On peut choisir $f$ telle que $f=e_U f
e_{\o{U}}$ et on peut alors {\'e}crire $f = e_U f_M^R e_{\o{U}}= f_M^R e_U
e_{\o{U}}$ avec $f_M^R \in  R[M]$. On a alors $f_M^Rz_{U,\o{U}}^R e_U e_{\o{U}} = f_M^R
e_Ue_{\o{U}}e_Ue_{\o{U}}=e_Ue_{\o{U}}$. Par l'isomorphisme \ref{eqIwa}, ceci montre que $f_M^R$ est un inverse de $z_{U,\o{U}}^R$ dans $R[M]$. 

\end{proof}

\alin{Exemple de calcul de l'\'el\'ement $z_{U,\o{U}}$ pour $SL(2)$} Nous
incluons ce calcul explicite en r\'eponse \`a une question  de
G. Henniart : on suppose que $G$ est
le pro-$p$-radical du sous-groupe d'Iwahori ``standard'' de $SL(2,\QM_p)$,
c'est \`a dire le groupe form\'e des matrices \`a
coefficients entiers dont la r\'eduction modulo $p$ est unipotente sup\'erieure.
On prend alors pour $U$ les matrices de $G$ qui sont unipotentes sup\'erieures et pour 
$\o{U}$ les unipotentes inf\'erieures. Pour $M$ on prend les
diagonales. On a donc des isomorphismes
$$\application{m:\;}{1+p\ZM_p}{M}{z}{\left(\begin{array}{cc} z^{-1} & 0 \\
      0 & z \end{array}\right)},
\application{u:\;}{\ZM_p}{U}{x}{\left(\begin{array}{cc} 1 & x \\
      0 & 1 \end{array}\right) }\hbox{ et }
\application{\o{u}:\;}{\ZM_p}{\o{U}}{y}{\left(\begin{array}{cc} 1 & 0 \\
      py & 1 \end{array}\right).}$$ 
Un calcul \'el\'ementaire montre que l'unique \'el\'ement de $M$ tel
que $\o{u}(y)u(x)\in Um\o{U}$ est $m(1+pxy)$.
On peut alors exprimer la distribution $z_{U,\o{U}}$ sous la forme 
$$ z_{U,\o{U}} =\int_{\ZM_p \times \ZM_p} (1+pxy)dxdy, $$
ce qui signifie que si $\phi$ est une fonction lisse sur $M$, alors
$<z_{U,\o{U}},\phi>=\int_{\ZM_p \times \ZM_p} \phi(1+pxy)dxdy$ o\`u $dx$ et
$dy$ sont des mesures de Haar normalis\'ees.
En particulier si $M_n$ est l'image de $1+p^{n+1}\ZM_p$ dans $M$ par
l'isomorphisme pr\'ec\'edemment d\'ecrit, on obtient en notant $1_?$
la fonction carat\'eristique de $?$ :
$$z_{U,\o{U}}*1_{M_n} = \sum_{z\in \ZM/p^n\ZM} a(z) 1_{m(1+pz)M_n}$$
o\`u $a(z)=\frac{1}{p^{2n}} \#\{(x,y)\in \ZM/p^n\ZM^2,\;\; xy=z\}$.

\alin{Foncteurs associ\'es}
On consid\`ere par la suite un groupe profini $G$ muni d'une d\'ecomposition d'Iwahori $G=UM\o{U}$ {\em satisfaisant l'hypoth\`ese de la proposition \ref{propIwa}}. %telle que $U$ et $\o{U}$ soient pro-$p$, ce qui permet de d\'efinir les idempotents $e_U$ et $e_{\o{U}}$.
 Au $(RG,RM)$-bimodule lisse $E_{U\o{U}}:= \CC^\infty_R(G)e_U e_{\o{U}} $ est associ\'e le couple de foncteurs adjoints $(I_{U\o{U}},R_{U\o{U}})$  d\'efini par :
$$ \application{I_{U\o{U}}:\;}{\Mo{R}{M}}{\Mo{R}{G}}{A}{E_{U\o{U}}\otimes_{RM}A} \;\hbox{ et }\; 
\application{R_{U\o{U}}:\;}{\Mo{R}{G}}{\Mo{R}{M}}{B}{\hom{E_{U\o{U}}}{B}{RG}^\infty} $$
et au $(RM,RG)$-bimodule $E'_{U\o{U}}:= e_Ue_{\o{U}}\CC^\infty_R(G)$ est associ\'e le couple  adjoint $(R'_{U\o{U}},I'_{U\o{U}})$ o\`u
$$ \application{I'_{U\o{U}}:\;}{\Mo{R}{M}}{\Mo{R}{G}}{A}{\hom{E'_{U\o{U}}}{A}{RM}^\infty} \;\hbox{ et }\; 
\application{R'_{U\o{U}}:\;}{\Mo{R}{G}}{\Mo{R}{M}}{B}{E'_{U\o{U}} \otimes_{RG} B}. $$
%o\`u $E^*_{U\o{U}}$ est le $(RM,RG)$-bimodule obtenu en composant les actions avec $m\mapsto m^{-1}$ et $g\mapsto g^{-1}$.
Le signe $\infty$ d\'esigne la partie lisse du $RG$ ou $RM$-module concern\'e. %(si l'objet $A$ ou $B$  est de type fini, ce signe est superflu, mais en g\'en\'eral il est n\'ecessaire).

\begin{coro} \label{corIwa} 
Les foncteurs $R_{U\o{U}}$ et $R'_{U\o{U}}$ sont isomorphes au foncteur $B\mapsto e_Ue_{\o{U}} B$, et les foncteurs $I_{U\o{U}}$ et $I'_{\o{U}U}$ sont isomorphes. 
\end{coro}

\begin{proof}
L'assertion sur $R$ et $R'$ est une cons\'equence imm\'ediate de la propri\'et\'e d'idempotence de $z_{U,\o{U}}^{-1}e_Ue_{\o{U}}$ et de l'inversibilit\'e de $z_{U,\o{U}}^{-1}$.

Pour l'assertion concernant $I$ et $I'$, on remarque d'abord que pour tout $i$ le diagramme
 suivant est commutatif (``infl'' d{\'e}signe l'inflation et on pose $U_i:=U\cap G_i$, etc...) :
 $$\xymatrix{ \Mo{R}{M/M_i} \ar[r]^{infl} \ar[d]^{I_{U/U_i\o{U/U_i}}} &
               \Mo{R}{M} \ar[d]^{I_{U\o{U}}}  \\
            \Mo{R}{G/G_i} \ar[r]^{infl} & \Mo{R}{G} }$$
ainsi que son analogue pour $I'$. Ceci nous ram\`ene au cas o\`u $G$ est {\em fini}. 
Dans ce cas l'application
$$ \varphi :\; \application{}{R[G]\otimes_R A }{\hom{R[G]}{A}{R}}{\sum_{g\in G} a_g g}{\left(g\mapsto a_{g^{-1}}\right)} $$
est un isomorphisme de $(RG,RM)$-bimodules qui, puisque  $z_{U,\o{U}}^{-1}e_Ue_{\o{U}}$ est idempotent, induit un isomorphisme de $(RG,RM)$-bimodules
$$ \varphi :\;\; R[G]e_Ue_{\o{U}} \otimes_R A  \simto \hom{e_{\o{U}}e_U R[G]}{A}{R} .$$
%En appliquant l'idempotent $e_M$ de $RM$, on obtient l'isomorphisme $I_{U\o{U}}(A) \simto I_{\o{U}U}'(A)$ cherch\'e.
Maintenant $I_{U\o{U}}(A) \simeq (R[G]e_Ue_{\o{U}} \otimes_R A )_M$ et $I_{\o{U}U}'(A) \simeq \hom{e_{\o{U}}e_U R[G]}{A}{R}^M$, donc la compos\'ee de $\varphi$ avec la trace sur $M$ induit  un morphisme des co-invariants vers les invariants
$$ \varphi \circ \hbox{Tr}_M :\;\; I_{U\o{U}}(A) \To{} I_{\o{U}U}'(A).$$
Puisque $R[G]e_Ue_{\o{U}}$ est projectif sur $R[M]$, le lemme suivant montre que c'est un isomorphisme.
\end{proof}

\begin{lemme} \label{lemmebete}
Soit $M$ un groupe fini et $R$ un anneau commutatif. Pour toute paire de $R[M]$-modules $(B,A)$, si $B$ est projectif alors l'application trace  
$$ \application{\hbox{Tr}_{M} :\;}{(B\otimes_R A)_{M}}{(B\otimes_R
A)^{M}}{b\otimes a}{\sum_{m} mb\otimes ma} $$
est un isomorphisme de $R$-modules.
\end{lemme}
\begin{proof}
On se ram{\`e}ne au cas
o{\`u} $B$ est libre de rang $1$ sur $R[M]$ que l'on v{\'e}rifie ``\`a la main".
\end{proof}

On \'etudie maintenant la d\'ependance des foncteurs ci-dessus en la d\'ecomposition d'Iwahori de $G$ :

\begin{lemme} \label{lemIwa} %M\^eme hypoth\`ese sur $G$ que dans le corollaire \ref{corIwa}.
Soit $G=VM\o{V}$  une seconde d{\'e}composition d'Iwahori de $G$
 telle que $G^\dag=VM^\dag \o{V}$, et ``compatible" \`a la d\'ecom\-position $G=UM\o{U}$,  au
sens o{\`u} $V=(V\cap U)(V\cap \o{U})$, $U=(U\cap V)(U\cap \o{V})$ etc... Alors   l'application
$$\application{}{\CC^\infty_R(G)e_{U}e_{\o{U}}}{\CC^\infty_R(G)e_{V}e_{\o{V}}}{f}{f*e_{\o{V}}}$$
est un isomorphisme de $(RG,RM)$-bimodules. 
%d'inverse l'application $f\mapsto e_U*f$. 
\end{lemme}
\begin{proof}
En effet, on a
\begin{eqnarray*}
  e_Ue_{\o{U}}e_{\o{V}}=  e_U e_{V} e_{\o{U}} e_{\o{V}} = e_U e_V e_{\o{V}}
  \end{eqnarray*}
ce qui montre que l'application est bien d{\'e}finie 
 et qu'elle est surjective puisque les inclusions
$$ \CC^\infty_R(G) e_{V}e_{\o{V}}e_V e_{\o{V}}   \subseteq \CC^\infty_R(G) e_U e_V
e_{\o{V}} \subseteq \CC^\infty_R(G) e_{V}e_{\o{V}} $$
sont des \'egalit\'es par la proposition pr\'ec\'edente.
Elle est de plus injective, car l'application 
$f\mapsto f*z_{U,\o{U}}^{-1} e_U e_{\o{U}}$ en est un inverse {\`a} droite.
\end{proof}

Comme cas particulier de ce lemme, les bimodules $E_{U\o{U}}$ et $E_{\o{U}U}$ sont isomorphes.
Posons alors $\Ipn{M}{G}:=I_{U\o{U}}$ et $\Rpn{G}{M}:=R_{U\o{U}}$. Ces foncteurs ne d\'ependent pas, \`a isomorphisme pr\`es, de l'ordre entre $U$ et $\o{U}$.
%Le corollaire \ref{corIwa} et la preuve de la proposition \ref{propIwa} montrent :

\begin{coro} \label{cor2Iwa} Sous les hypoth\`eses de la proposition \ref{propIwa} et avec les notations ci-dessus,
\begin{enumerate}
\item
Les foncteurs $\Ipn{M}{G}$ et $\Rpn{G}{M}$ sont adjoints {des deux c\^ot\'es}. 
\item Le compos\'e  $\Rpn{G}{M}\circ \Ipn{M}{G}$ est isomorphe au foncteur identit\'e de $\Mo{R}{M}$.
%Si l'ordre de $M/M^\dag$ est inversible dans $R$, alors l
\item Le foncteur $\Ipn{M}{G}$ envoie objets irr\'eductibles sur objets irr\'eductibles.
\end{enumerate}
\end{coro}

\begin{proof}
La premi\`ere assertion d\'ecoule du  corollaire \ref{corIwa} et du lemme pr\'ec\'edent. Pour la deuxi\`eme assertion, il suffit de v\'erifier que le morphisme de $(RM,RM)$-bimodules
$$ \application{}{\CC^\infty_R(M)}{e_Ue_{\o{U}}\CC^\infty_R(G)e_Ue_{\o{U}}}{f}{e_Ue_{\o{U}}fe_Ue_{\o{U}}}$$
est un isomorphisme. Or, la d\'ecomposition d'Iwahori montre qu'il est surjectif, et la propri\'et\'e \ref{propIwa} de presqu'idempotence de $e_Ue_{\o{U}}$ montre qu'il est injectif.
La troisi\`eme assertion
a \'et\'e d\'emontr\'ee au cours de la preuve de la proposition \ref{propIwa} dans le cas o\`u $G$ est fini d'ordre une puissance de $p$. Le cas o\`u $G$ est pro-$p$ s'ensuit, par l'axiome ii) de \ref{defIwa}. Dans le cas g\'en\'eral, il faut commencer par v\'erifier  la propri\'et\'e de commutation aux foncteurs d'oublis  $\Res{G}{G^\dag}\circ \Ipn{M}{G} \simeq \Ipn{M^\dag}{G^\dag}\circ \Res{M}{M^\dag}$, ce que nous laisserons au lecteur.
Fixons ensuite $W\in\Mo{R}{M}$ irr\'eductible. Sa restriction \`a $M^\dag$ est semi-simple.
Par  le point iii) pour $G^\dag$ et le point ii), on a la propri\'et\'e suivante :
pour tout sous-$RG^\dag$-module $V$ de $\Ipn{M}{G}(W)$, on a 
$\hbox{long}_{RG^\dag}(V) = \hbox{long}_{RM^\dag}(\rpn{G}{M}{V})$.
Ainsi, si $V$ est un sous $RG$-module non-nul, alors le $RM$-module  $\rpn{G}{M}{V}$ est non-nul et, par exactitude de $\Rpn{G}{M}$, est un sous-module  de $\rpn{G}{M}{\ipn{M}{G}{W}}\simeq W$ donc lui est  \'egal. On en d\'eduit que  $\hbox{long}_{RG^\dag}(\ipn{M}{G}{W}) = \hbox{long}_{RM^\dag}(V)$ et donc $V=\Ipn{M}{G}{W}$.
\end{proof}

\alin{Induction parahorique}\def\GM{\mathbb G} \label{notaimmeuble}
Reprenons les notations de l'introduction. Au $K$-groupe r\'eductif $\GC$, Bruhat et Tits ont  associ\'e un immeuble affine ``\'etendu" $B(\GC,K)$, muni d'une action ``affine" de $G$. Cet immeuble n'est d\'efini qu'\`a isomorphisme pr\`es ; en fait le groupe des automorphismes affines et $G$-\'equivariants de $B(\GC,K)$ s'identifie canoniquement au $\RM$-espace vectoriel $a_G:=\hom{\GM_m}{\ZC(\GC)}{K-gr}\otimes \RM$. %\hom{X^*(\GC)_K}{\RM}{\ZM}.

Si $\MC$ est un $K$-\levi de $\GC$, alors le sous-ensemble $M.A(\GC,\SC,K)$ des points de $B(\GC,K)$ dont la $M$-orbite rencontre l'appartement $A(\GC,\SC,K)$ associ\'e \`a un  $K$-tore d\'eploy\'e maximal $\SC$ de $\MC$ ne d\'epend pas du choix de ce $K$-tore et est un immeuble \'etendu pour $\MC$ relativement \`a $K$. Nous le noterons donc $B(\MC,K)$ ; il est muni d'une action de $a_M$ commutant \`a $M$ que nous noterons $(x,\lambda)\mapsto x+\lambda$ o\`u $\lambda\in a_M$ et $x\in B(\MC,K)$.

Si $x\in B(\GC,K)$, on note $G_x$ son stabilisateur dans
  $G$ ; c'est un sous-groupe ouvert compact de $G$.
   La th{\'e}orie de Bruhat et Tits (notre r{\'e}f{\'e}rence sera
  \cite[3.4]{Tits}) associe {\`a} $x$ un mod\`ele entier lisse  $\GC_x$ de $\GC$ tel que $\GC_x(\OC)=G_x$. 
Par lissit\'e, la r{\'e}duction  $ \varpi : \;\; G_x=\GC_x(\OC) \To{} \GC_x(k)$ est
surjective. %On notera $G_x^0:= \varpi^{-1} ( \o\GC^0_x(k) )$ et
On pose $G_x^+:=\varpi^{-1} ({^u\GC_{x,k}}(k))$ o\`u ${^u\GC_{x,k}}$ d\'esigne le radical unipotent de la fibre sp\'eciale $\GC_{x,k}$ de $\GC_x$ ; c'est un pro-$p$-sous-groupe ouvert et normal
de $G_x$, parfois appel\'e {\em pro-$p$-radical}. 
% Pour $\dag=\emptyset,0,+$, les groupes
% $G_x^{\dag}$ sont ouverts compacts dans $G$ et $G_x^+$ est un
% pro-$p$-groupe. 
Supposons de plus que $x\in B(\MC,K)$ ; le m\^eme proc\'ed\'e que ci-dessus fournit un pro-$p$-radical $M_x^+$ qui fort heureusement co\"incide avec l'intersection $M\cap G_x^+$ (voir aussi \ref{defmod} pour une situation plus g\'en\'erale).
% En fait, l'immersion $\MC\injo \GC$ se prolonge de mani\`ere unique en une immersion $\MC_x \injo \GC_x$.
  De mani{\`e}re g{\'e}n{\'e}rale, pour  tout sous-groupe $H$ de $G$, on notera $H_x:=H\cap
    G_x$  et $H_x^+:=H\cap G_x^+$.

Soit $\PC=\MC\UC$ un $K$-\para et $\o{\PC}=\MC\o\UC$ son oppos\'e par rapport \`a $\MC$ ({\em i.e } l'unique \para de $\GC$ tel que $\PC\cap \o\PC=\MC$).
Si $x\in B(\MC,K)$, on sait  \cite[4.6.8]{BT2} (voir aussi \ref{dilatIwa}) que $G_x^+$ admet une d\'ecomposition d'Iwahori $G_x^+=U_x^+M_x^+ \o{U}_x^+$ au sens de \ref{defIwa}. Il s'ensuit que le groupe $G_{x,P}:= G_x^+ P_x$ admet la d\'ecomposition d'Iwahori $U_x M_x\o{U}_x^+$ et satisfait l'hypoth\`ese de la proposition \ref{propIwa} avec $M_x^+$ pour $M^\dag$. Celle-ci nous fournit donc un \'el\'ement 
\ini\begin{equation}
z_{x,P} \in \ZC(RM_x)^\times \hbox{ tel que } z_{x,P}\eux\eubxp \hbox{ est un idempotent de } RG_{x,P}.
\end{equation}

On obtient aussi deux foncteurs reliant $\Mo{R}{M_x}$ et $\Mo{R}{G_x}$
qui sont adjoints des deux c\^ot\'es :
\ini\begin{equation}
I_{x,P}:=\cInd{G_{x,P}}{G_x}\circ \Ipn{M_x}{G_{x,P}}  \hbox{ et } R_{x,P}:= \Rpn{G_{x,P}}{M_x}\circ \Res{G_x}{G_{x,P}}
\end{equation}
avec les notations du corollaire \ref{cor2Iwa}.

Par un l{\'e}ger abus de langage que nous expliquons ci-dessous, nous dirons que $G_{x,P}$ est un
sous-groupe parahorique de $G_x$. 
Le foncteur $I_{x,P}$  m\'erite alors le nom  d'{\em induction parahorique}
en ce qu'il relie les repr\'esentations de $M_x$ et $G_x$ en passant par le groupe parahorique $G_{x,P}$, comme l'induction parabolique relie les repr\'esentations de $M$ \`a $G$ en passant par $P$. L'innocent foncteur d'inflation des repr\'esentations de $M$ \`a $P$ est ici remplac\'e par le foncteur $\Ipn{M_x}{G_{x,P}}$. Ce dernier n'est pas si facile \`a d\'efinir  mais  partage certaines des propri\'et\'es  de l'inflation : il envoie irr\'eductibles sur irr\'eductibles et la compos\'ee avec son adjoint est l'identit\'e de $\Mo{R}{M_x}$, {\em cf} \ref{corIwa}.

Dans  la terminologie de Bruhat-Tits  \cite[5.2.6]{BT2},  un
sous-groupe parahorique de $G$ est un 
sous-groupe de la forme $\varpi^{-1}(\PC_k(k))$, o{\`u} $\PC_k$ est un
$k$-\para  de la fibre sp\'eciale {\em connexe} $\GC^\circ_{x,k}$.  En  posant $M_x^\circ = \varpi^{-1}(\MC_{x,k}^\circ(k))$, on a une application surjective
$$
\begin{array}{ccc}\{ K-\hbox{ss-gr paraboliques contenant } \MC\} & \To{} &
\{k-\hbox{ss-gr paraboliques contenant }
  \MC_{x,k}^\circ\}  \\ 
& & || \\ & & \{\hbox{ss-gr parahoriques de } G_x \hbox{ contenant }
M_x^\circ \}
\end{array} 
$$
qui n'est  injective  que  si $x$ est un
sommet hypersp{\'e}cial. 
Dans la terminologie usuelle, le sous-groupe parahorique  associ{\'e} {\`a} un parabolique $\QC=\MC\VC$ s'{\'e}crit
$G_x^+M_x^\circ V_x$, tandis que dans cet article, nous appelons 
sous-groupe parahorique les sous-groupes de la forme
$G_{x,Q}:=G_x^+M_xV_x$ ; les deux terminologies peuvent diff{\'e}rer mais sont
``en bijection''. Quoiqu'il en soit, le petit diagramme ci-dessus montre
qu'un m{\^e}me parahorique peut  {\^e}tre associ{\'e} {\`a}  deux
paraboliques contenant $M$ diff{\'e}rents. 
Malgr{\'e} les apparences, nos
foncteurs d'induction-restriction ne d{\'e}pendent, eux, que du
parahorique :
\begin{lemme} \label{dependance}
  Si $\QC=\MC.\VC$ et $\PC=\MC.\UC$ sont deux \paras tels que
  $G_{x,P}=G_{x,Q}$, les foncteurs $I_{x,P}$ et $I_{x,Q}$ sont canoniquement isomorphes.
  %les deux $(M_x,G_x)$-bimodules $E_{x,P}$ et
  %$E_{x,Q}$ sont isomorphes. 
\end{lemme}
\begin{proof}
D'apr\`es  \cite[3.1]{Tits}, les d\'ecompositions d'Iwahori $G_{x,P}=U_x M_x \o{U}_x^+$ et $G_{x,Q}= V_x M_x \o{V}_x^+$ sont ``compatibles" au sens du lemme \ref{lemIwa}. Il suffit donc d'appliquer ce lemme. 
\end{proof}

Une fois r\'egl\'e ce probl\`eme, une autre question appara\^it naturellement : {\em les foncteurs parahoriques d\'ependent-ils du groupe parahorique ?}
Dans le cas des groupes r{\'e}ductifs {\em finis} (sur $\FM_p$),
Howlett et Lehrer ont montr{\'e} dans
\cite{HL} que l'induction parabolique pour les repr{\'e}sentations {\`a}
coefficients dans $\zp$ 
ne d{\'e}pend pas du choix du parabolique contenant un Levi donn{\'e}. 
En s'inspirant de leur preuve, on est amen\'e  {\`a} la question suivante :

\begin{ques} \label{question}
  Soit $\PC=\MC\UC$ un $K$-\para de $\GC$. A-t-on pour tout $x\in B(\MC,K)$ 
  $$\euxp\eubx \in (\zp
  G_x)\eux\eubx \;\hbox{ et }\; \eux\eubxp \in   \eux\eubx (\zp G_x) ?$$ 
\end{ques}
Si cette question avait une r{\'e}ponse affirmative, alors l'application 
$$\application{}{\eux\eubxp\CC^\infty(G_x)}
                {\eubx\euxp\CC^\infty(G_x)}
                {f}{\eubx*f}$$
serait un isomorphisme de
 $(M_x,G_x)$-bimodules et induirait des isomorphismes $I_{x,P}\simto I_{x,\o{P}}$.
Mais l'int\'er\^et principal de la question \ref{question} apparaitra dans la prochaine section.
Le seul cas o\`u nous pouvons lui donner une r\'eponse positive ``sans restrictions"  est celui o\`u $\PC$ est minimal, {\em cf \ref{pfcommutmin}}.

\section{Commutation aux foncteurs paraboliques} \label{rescentpf}

Soit $\PC=\MC\UC$ un sous-groupe parabolique.
Rappelons que le foncteur de
  Jacquet (ou restriction parabolique) associe {\`a} $V\in \Mo{R}{G}$ 
   le $M$-module lisse $V_U:=V/V(U)$ o{\`u} $V(U):=\bigcup_{K\subset U} \ker
   e_K$ o{\`u} $K$ d{\'e}crit l'ensemble des sous-groupe ouverts compacts de
   $U$. On notera aussi parfois  $\Rp{G,P}{M}$ ce foncteur et 
  $\Ip{M,P}{G}$ son adjoint {\`a} droite (induction parabolique), ainsi que
  $j_U : V\To{} V_U$ la surjection canonique. Nos conventions sur les immeubles sont celles du paragraphe \ref{notaimmeuble}, et on note toujours $\o\PC=\MC\o\UC$ le parabolique de $\GC$ oppos\'e \`a $\PC$ par rapport \`a $\MC$.

\begin{prop} \label{rescent}
  Soit $\PC=\MC\UC$ un \para et  
 $y \in B(\MC,K)$. Soit $\ve$ un idempotent de $RM_y$ et $\wt\ve:= z_{y,P}^{-1}e_{U_y^{}}e_{\o{U}_y^+}\ve$ l'idempotent de $RG_y$ associ\'e.
Supposons que pour tout $x\in y+a_M$, on a 
$$ \euxp\eubx\ve \in (RG_x)\eux\eubx\ve \;\hbox{  et }\; \ve\eux\eubxp \in
  \ve\eux\eubx (RG_x).$$ 
   Alors pour tout objet $V\in \Mo{R}{G}$, la projection canonique
  $ V \To{j_U} V_U $ induit un isomorphisme $\wt\ve V \simto \ve V_U $ de $R$-modules.
\end{prop}

Rappelons que pour tout $x\in y+a_M$, on a $M_x=M_y$ de sorte que l'idempotent $\ve$, vu dans $RG_x$, commute aux idempotents $e_{U_x}$, $e_{\o{U}_x}$ etc...

\medskip

\begin{proof}
Le fait que $j_U$ envoie $\wt\ve V$ dans $\ve V_U$ est une simple cons\'equence de sa $M_y$-\'equivariance.
Fixons un {\'e}l{\'e}ment $z_M$ du centre de $M$ dont l'action par
conjugaison sur $\o{U}$ (resp. $U$) soit strictement contractante
(resp. dilatante). L'action de $z_M$
sur $B(\MC,K)$ est la translation par un  vecteur  not{\'e} $\o{z}_M\in a_M$. On
s'int{\'e}resse au comportement de $U_x^{(+)}$ et $\o{U}_x^{(+)}$ lorsque
$x$ d{\'e}crit la demi-droite d'origine $y$ et de vecteur directeur $\o{z}_M$.
On note pour cela $x(t):=y+t\o{z}_M$.

\begin{lemme} \label{l1} Il existe un ensemble discret
  $I_y=\{0\leq t_1 < t_2 <\cdots < t_n <\cdots\}\subset \RM_+$ tel que
  si $t\in ]t_i,t_{i+1}[$ alors :
  \begin{itemize}
  \item $U_{x(t)}=U^+_{x(t)}=U_{x(t_i)}=U^+_{x(t_{i+1})}$
  \item
    $\o{U}_{x(t)}=\o{U}^+_{x(t)}=\o{U}^+_{x(t_i)}=\o{U}_{x(t_{i+1})}$ 
  \end{itemize}
et tel que pour tout $i$, on a $U^+_{x(t_i)} \subsetneq U_{x(t_i)}$ et
$\o{U}^+_{x(t_i)} \subsetneq \o{U}_{x(t_i)}$. 
\end{lemme}

\begin{proof}
  Voir \cite[6.7]{MP2}.
\end{proof}

En particulier, les applications
$t\mapsto U_{x(t)}^{(+)}$ et $t \mapsto \o{U}_{x(t)}^{(+)}$ sont
respectivement croissante et d{\'e}croissante. Rappelons aussi la propri\'et\'e suivante :
\begin{fact} \label{f1}
  On a $\bigcup_{t\geq 0} U_{x(t)}^{(+)} = U$ et $\bigcap_{t\geq 0}
  \o{U}_{x(t)}^{(+)} ={1}$.
\end{fact}

Soit maintenant $V\in \Mo{}{G}$, on veut montrer d'abord que $\wt\ve V\cap V(U) = 0$. Soit $w$ un {\'e}l{\'e}ment de cette intersection, alors par
d{\'e}finition de $V(U)$ et par le fait pr{\'e}c{\'e}dent, l'ensemble
$$ I_w:=\{ t\in \RM_+ , e_{U_{x(t)}}w =0 \} $$
est non vide. Il admet donc une borne inf{\'e}rieure $t_w$ qui, par \ref{l1}
(semi-continuit{\'e} sup{\'e}rieure) appartient encore {\`a} $I_w$ et donc
n{\'e}cessairement {\`a} $I_y$. Toujours par
\ref{l1} (m{\^e}me notation) on a $t_w=0$ ou $t_w \in I_y\setminus\{0\}$. Le premier
cas est trivial : on obtient $e_{U_y} w = w =0$. Dans le deuxi{\`e}me cas, on
consid{\`e}re l'{\'e}l{\'e}ment $w'=e_{U_{x(t_w)}^+} w$, alors
\begin{itemize}
\item $w'\neq 0$ car il existe $\varepsilon >0$ tel que
  $U_{x(t_w)}^+=U_{x(t_w-\varepsilon)}$ d'apr{\`e}s \ref{l1}.
\item $w' \in e_{U_{x(t_w)}^+}e_{\o{U}_y^+} \ve V \subset
  e_{U_{x(t_w)}^+}e_{\o{U}_{x(t_w)}} \ve V$ car $\o{U}_{x(t_w)} \subset
  \o{U}_{y}^+$, toujours par \ref{l1}.
\end{itemize}

On a donc $w'=e_{U_{x(t_w)}^+}e_{\o{U}_{x(t_w)}} \ve v$ pour un certain
$v\in V$ et $e_{U_{x(t_w)}} w' =0$. Or par hypoth\`ese il
existe $\phi \in R G_{x(t_w)}$ tel que $\phi e_{U_{x(t_w)}} w'= w'$
ce qui contredit la non-nullit{\'e} de $w'$.
On a donc obtenu l'injectivit{\'e} de la restriction de $j_U$ \`a $\wt\ve V$.

Pour la surjectivit{\'e}, fixons $v \in V$ et
remarquons que pour $t$ assez grand, pr{\'e}cis{\'e}ment lorsque $e_{\o{U}^+_{x(t)}} v 
= v$, l'{\'e}l{\'e}ment $e_{\o{U}^+_{x(t)}} v$ de $V$ a la m{\^e}me image 
que $v$ dans $V/V(U)$. On a donc $V_U = \cup_{t \in \RM_+}
j_U(e_{\o{U}^+_{x(t)}}V)$.
En cons{\'e}quence, si on fixe maintenant $w \in \ve V_U$, l'ensemble 
$$ J_w:=\{ t\in \RM_+ ,w \in j_U(e_{\o{U}^+_{x(t)}}\ve\,V) \} $$
est non vide et admet une borne inf{\'e}rieure $t_w$ qui comme
pr{\'e}c{\'e}demment appartient {\`a} $J_w$ et {\`a} $I_y$. Supposons $t_w\neq 0$ et fixons
 $v\in V$ tel que $j_U(e_{\o{U}^+_{x(t_w)}} \ve v)=w$.
D'apr{\`e}s notre hypoth\`ese, il existe $f \in RG_{x(t_w)}$ telle que
$e_{U_{x(t_w)}}e_{\o{U}^+_{x(t_w)}}\ve =
  e_{U_{x(t_w)}}e_{\o{U}_{x(t_w)}} \ve f$. Il s'ensuit que 
    \begin{eqnarray*}
      w = j_U(e_{U_{x(t_w)}}e_{\o{U}^+_{x(t_w)}} \ve v)=
      j_U(e_{U_{x(t_w)}}e_{\o{U}_{x(t_w)}}\ve f v) = j_U(e_{\o{U}_{x(t_w)}} \ve f v)
    \end{eqnarray*}
ce qui contredit la minimalit{\'e} de $t_w$, car $\o{U}_{x(t_w)}^+\subsetneq
  \o{U}_{x(t_w)}=\o{U}_{x(t_w-\varepsilon)}^+$, pour $\varepsilon$
  assez petit. On en d{\'e}duit que $t_w=0$, donc $\ve V_U = j_U(e_{\o{U}_y^+}\ve V)=j_U(\wt\ve 
  V)$.

\end{proof}

\begin{rema}  \label{rem} Comme cas particulier de la proposition pr\'ec\'edente, 
l'application $$\application{}
{ \CC^{\infty,c}_{R}(G)}{\CC^{\infty,c}_{R}(U\ba
  G)}{f}{(g \mapsto      \int_U f(ug)du)}$$ 
induit un  isomorphisme ${\wt\ve\, \CC^{\infty,c}_{R}(G)}\simto{\ve\,\CC^{\infty,c}_{R}(U\ba
  G)}$ de $RG$-modules. 
\end{rema}
En effet, l'application $f \mapsto (g\mapsto  \int_U f(ug)du)$ 
se factorise par $j_U$ et  induit un isomorphisme
$\CC_R^{\infty,c}(G)_U\simto\CC_R^{\infty,c}(U\ba G)$.

Gardons les notations de la proposition pr\'ec\'edente et notons $\Mo{}{\ve RM_x \ve}$ la cat\'egorie des modules sur l'anneau $\ve RM_x \ve = \ve RM_y \ve$, o\`u $x\in y+ a_M$. On a la paire de foncteurs  
$$ \application{T_\ve :\;}{\Mo{}{\ve RM_x \ve}}{\Mo{R}{M_x}}{B}{\CC^\infty_R(M_x) \ve \otimes_{\ve RM_x\ve} B} \,\hbox{ et }\, \application{\ve :\;}{\Mo{R}{M_x}}{\Mo{}{\ve RM_x\ve}}{V}{\ve. V}. $$

\begin{coro} \label{commut}
Sous les m\^emes hypoth\`eses que la proposition \ref{rescent}, pour tout $x\in y+a_M$ on a
 \begin{enumerate}
 \item Les foncteurs $\ve.R_{x,P} \circ \Res{G}{G_x}$ et $\ve. \Res{M}{M_x}\circ \Rp{G,P}{M}$ sont isomorphes.
    \item   %Soit $T_\ve :\; \Mo{}{\ve RM_x \ve} \To{} \Mo{R}{M_x}$ le foncteur $B\mapsto (\CC^\infty_R(M_x) \ve \otimes_{\ve RM_x\ve} B)$   
 Les foncteurs $ \Ip{M,P}{G}\circ \cInd{M_{x}}{M} \circ T_\ve$  et $\cInd{G_x}{G}\circ I_{x,P}\circ T_\ve$
%de $\Mo{}{\ve RM_x\ve}$ dans $\Mo{R}{G}$ 
sont isomorphes. 
%     Si $\tau:= \CC^\infty_R(M_x) \ve \in \Mo{R}{M_x}$, alors
%  $\ip{M,P}{G}{\cind{M_x}{M}{\tau}} \simeq
%  \cind{G_x}{G}{I_{x,P}(\tau)}$.
   \end{enumerate}
 \end{coro}

\begin{proof}
Le premier point est une cons\'equence imm\'ediate de la proposition \ref{rescent} appliqu\'ee \`a $x$ au lieu de $y$ (l'hypoth\`ese de \ref{rescent} est en effet ``invariante" par translation sous $a_M$) et du fait que l'application $j_U$ est $M_x$-\'equivariante.

Le deuxi\`eme point {\em ne} se d\'eduit {\em pas} formellement du premier par adjonction.
Soit  $\tau_x$ un objet de $\Mo{R}{M_x}$. Rappelons que par d{\'e}finition, on a
$$ \ip{M,P}{G}{\cind{M_x}{M}{\tau_x}} \simeq (\tau_x\otimes_R
\CC_R^{\infty,c}(U\ba G))^{M_x}$$
o{\`u} $M_x$ agit diagonalement sur le produit tensoriel. Pour prouver le point ii), il nous faut d'abord modifier cette expression en rempla{\c c}ant ``invariants" par ``coinvariants". Pour cela 
rappelons que l'action de la fonction caract\'eristique $1_{M_x}$ de $M_x$ induit un morphisme du foncteur des coinvariants vers celui des covariants : pour tout objet $W$ de $\Mo{R}{M_x}$, on a donc une application ``trace" $R$-lin\'eaire (d\'ependant du choix d'une mesure de Haar sur $M_x$) $\application{\hbox{Tr}_{M_x}:\;}{W_{M_x}}{W^{M_x}}{w}{1_{M_x}.w}$.
Nous allons prouver que pour $W= (\tau_x\otimes_R \CC^{\infty,c}_R(U\ba G))$, cette application est bijective, bien que le pro-ordre de $M_x$ ne soit pas inversible dans $R$.
Pour cela, on peut supposer $\tau_x$ de type fini et choisir un pro-$p$-sous-groupe ouvert normal  $M_{x,r}$ de $M_x$ dans le noyau de $\tau_x$ et dont on note $\o{M_x}:=M_x/M_{x,r}$ le quotient. On a alors une factorisation $\hbox{Tr}_{M_x} = \hbox{Tr}_{\o{M_x}}\circ \hbox{Tr}_{M_{x,r}}$. Puisque $M_{x,r}$ est pro-$p$ et $p$ est inversible dans $R$, $\hbox{Tr}_{M_{x,r}}$ est un isomorphisme, et on a donc 
$$ \hbox{Tr}_{M_{x,r}} :\; (\tau_x\otimes_R \CC^{\infty,c}_R(U\ba G))_{M_{x,r}} \simto (\tau_x\otimes_R \CC^{\infty,c}_R(U\ba G))^{M_{x,r}} \simeq  \tau_x\otimes_R \CC^{\infty,c}_R(M_{x,r}U\ba G).$$
On est donc ramen\'e \`a \'etudier la trace sous le groupe fini $\o{M_x}$. 
D'apr\`es le lemme \ref{lemmebete}, il nous suffira  de prouver que le $R[\o{M_{x}}]$-module $\CC^{\infty,c}_R(M_{x,r}U\ba G)$ est projectif. Celui-ci est somme directe des $R[\o{M_x}]$-sous-modules
$$\CC^{\infty,c}_R(M_{x,r}U\ba M_xUgH)$$ pour $g\in M_xU\ba G/H$, et o\`u $H$ d\'esigne ici un sous-groupe ouvert compact fix\'e.
Soit $(H_n)_{n\in\NM}$  une suite d\'ecroissante et d'intersection triviale de
 pro-$p$-sous-groupes ouverts et normaux dans $H$ assez petits pour que $gH_ig^{-1}\cap (M_x\setminus M_{x,r})U =\emptyset$ pour tout $i$.
  Alors, pour tout $i$, l'action de $\o{M_x}$ sur l'ensemble $M_{x,r}U\ba M_xUgH /H_i$ est {\em libre}, et donc le $R[\o{M_x}]$-module $B_i:=R[M_{x,r}U\ba M_xUgH /H_i]$ est libre. Comme $H_i$ est pro-$p$, l'inclusion $B_i \injo B_{i+1}$ est scind\'ee sur $R[\o{M_x}]$ et par r\'ecurrence on construit une base de $\limi{i} B_i = \CC^{\infty,c}_R(M_{x,r}U\ba M_xUgH)$.
  
Finalement, on a donc prouv\'e que
$$ \ip{M,P}{G}{\cind{M_x}{M}{\tau_x}} \simeq (\tau_x\otimes_R
\CC_R^{\infty,c}(U\ba G))_{M_x} = \tau_x\otimes_{RM_x}
\CC_R^{\infty,c}(U\ba G)$$
pour tout objet $\tau_x\in \Mo{R}{M_x}$.
Prenons alors $\tau_x$ de la forme $T_\ve(B)$ pour un $\ve RM_x  \ve$-module $B$ ; on obtient
$$ \ip{M,P}{G}{\cind{M_x}{M}{T_\ve(B)}} %\simeq ((RM_x\ve\otimes_{\ve RM_x\ve} B)\otimes_R
%\CC_R^{\infty,c}(U\ba G))_{M_x} 
\simeq B\otimes_{\ve RM_x\ve} \ve\,\CC_R^{\infty,c}(U\ba G)$$

D'autre part, toujours pour un objet $\tau_x\in \Mo{R}{M_x}$ g\'en\'eral on a par d\'efinition
$$\cind{G_x}{G}{I_{x,P}(\tau_x)} \simeq
\left(\tau_x\otimes_{RM_x}\eux\eubxp\CC^{\infty}_R(G_x)\right)\otimes_{RG_x}
\CC_R^{\infty,c}(G)  \simeq \tau_x\otimes_{RM_x} \eux\eubxp \CC_R^{\infty,c}(G).    $$
En appliquant ceci \`a $\tau_x$ de la forme $T_\ve(B)$, on obtient
$$ \cind{G_x}{G}{I_{x,P}(T_\ve(B))} \simeq B \otimes_{\ve RM_x\ve} \wt\ve\,\CC^{\infty,c}_R(G).$$
On conclut donc gr\^ace \`a la $\ve RM_x \ve$-\'equivariance de l'isomorphisme de la remarque \ref{rem}.
\end{proof}

Pour \'enoncer le corollaire suivant, 
donnons quelques d\'efinitions : %un idempotent de $RM$ sera dit {\em admissible} si son centralisateur dans $M$ est ouvert. U
\begin{DEf} \label{pbon} Soit $\PC=\MC \UC$ un sous-groupe parabolique.
\begin{itemize}
        \item Un idempotent $\ve$ de $RM$ sera dit $P$-bon s'il   existe $y_\ve\in B(M,K)$ tel que $\ve\in RM_{y_\ve}$ et les hypoth\`eses de la proposition \ref{rescent} sont v\'erifi\'ees.
        \item Une famille $\EC$ d'idempotents de $RM$ sera dite {\em g\'en\'eratrice} si $\CC^{\infty,c}_R(M)= \sum_{\ve\in \EC} \CC^{\infty,c}_R(M)\ve$. 
\end{itemize}
\end{DEf}
Remarquons simplement qu'un idempotent $P$-bon est aussi $\o{P}$-bon et que la famille \`a un \'el\'ement $\{1\}$ est g\'en\'eratrice.

\begin{coro} \label{cororescent}
Soit $\PC=\MC \UC$ un sous-groupe parabolique. Supposons qu'il existe  une famille  g\'en\'eratrice $\EC$ d'idempotents $P$-bons de $RM$.
Alors
  \begin{enumerate}
  \item {\em Propri{\'e}t{\'e}s de finitude} :  l'induction parabolique $\Ip{M,P}{G}$ respecte la propri{\'e}t{\'e} d'{\^e}tre
    {\em de  type fini}, et la restriction parabolique $\Rp{G,P}{M}$ celle d'{\^e}tre
    {\em admissible}. De plus, pour tout pro-$p$-sous-groupe ouvert 
    $H$ de $G$ admettant une $(P,\o{P})$-d{\'e}composition, l'application
    canonique $V^H\To{j_U}\Rp{G,P}{M}(V)^{H\cap M}$ est surjective.
  \item {\em Seconde adjonction} : le foncteur $\delta_P^{-1}\Rp{G,P}{M}$ est adjoint
    {\`a} droite du foncteur ${\Ip{M,\o{P}}{G}}$ d'induction par
    rapport au parabolique oppos{\'e} $\o{P}$ tordu par le module.
   De plus, pour tout objet $V$ de $\Mo{R}{G}$, on a $\Rp{G,P}{M}{V}=0$ \ssi\
   ${\Rp{G,\o{P}}{M}}V=0$. 
  \end{enumerate}
 \end{coro}

\begin{proof}
Remarquons pour commencer, que si une telle famille $\EC$ existe, alors on peut en d\'eduire  une autre,  g\'en\'eratrice aussi, et dont les \'el\'ements sont des idempotents {\em lisses}, {\em i.e.} dans $\CC^{\infty,c}_R(G)$ : en effet, si $\ve\in \EC$ et $M_{y_\ve,r}$ est un pro-$p$-sous-groupe ouvert normal de $M_{y_\ve}$, alors $\ve$ commute \`a $e_{M_{y_\ve,r}}$ (qui est central dans $RM_{y_\ve}$) et le produit $\ve e_{M_{y_\ve,r}}$ est donc un idempotent. La famille obtenue en prenant pour chaque $\ve$ une base de voisinages de l'unit\'e form\'ee de tels sous-groupes ouverts normaux dans $M_{y_\ve}$ est g\'en\'eratrice et satisfait l'hypoth\`ese de l'\'enonc\'e. Nous supposerons d\'esormais que les idempotents de $\EC$ sont lisses.

Soit $\ve\in \EC$ et $\wt\ve$ l'idempotent de $RG_{y_\ve}$ associ\'e. Par construction cet idempotent est lisse dans $RG$ puisque $\ve$ l'est dans $RM$. Il s'ensuit que la repr\'esentation $\CC^{\infty,c}_R(G)\wt\ve$ est de type fini. Or, d'apr\`es le corollaire \ref{commut} ii) appliqu\'e au $\ve RM_{y_\ve} \ve$-module libre de rang $1$, on a
$$\ip{M,P}{G}{\CC^{\infty,c}_R(M)\ve} \simeq \CC^{\infty,c}_R(G)\wt\ve,$$
de sorte que l'induite parabolique de la repr\'esentation $\CC^{\infty,c}_R(M)\ve$ est de type fini. D'apr\`es notre hypoth\`ese, la famille des $\CC^{\infty,c}_R(M)\ve$ est g\'en\'eratrice dans la cat\'egorie $\Mo{R}{M}$. Puisque l'induction parabolique est un foncteur exact, on en d\'eduit qu'elle envoie objets de type fini sur objets de type fini.

Pour montrer l'"admissibilit\'e" de la restriction parabolique, il suffit bien-s\^ur de prouver la derni\`ere assertion du point i). Fixons donc $H=(H\cap \o{U})H_M(H\cap U)$ un pro-$p$-sous-groupe ouvert de $G$. Comme $H_M$ est pro-$p$, il nous suffira de prouver que la restriction de $j_U$ \`a $V^{H\cap\o{U}}\To{} V_U$ est {\em surjective}. Pour cela, il suffit de voir que pour tout $\ve\in\EC$, l'image $j_U(V^{H\cap\o{U}})$ contient $\ve V_U$.
Choisissons alors $y_\ve$ tel que $\o{U}_{y_\ve}^+\supset H\cap \o{U}$. On a 
$$ j_U(V^{H\cap\o{U}}) \supset j_U(e_{\o{U}_{y_\ve}^+}V) \supset j_U(\wt\ve V)=\ve V_U, $$
d'o\`u la surjectivit\'e recherch\'ee.

Venons-en maintenant {\`a} la propri{\'e}t{\'e} de seconde adjonction. En d\'eroulant la d\'efinition du foncteur d'induction parabolique, on peut en trouver un adjoint \`a droite. Ceci est d\^u \`a Bernstein \cite{bernunp} et n\'ecessite sa notion de ``compl{\'e}tion'' :
pour $V\in
\Mo{R}{G}$ on d\'efinit sa compl{\'e}tion $\hat{V}$ par 
$$ \hat{V} := \hom{\HC_R(G)}{V}{G} $$ %\simeq \limproj_H{V^H} $$
Le $R$-module obtenu $\hat{V}$ est canoniquement un $RG$-module
puisque $\HC_R(G)$ est un id\'eal bilat\`ere de $RG$. %$RG= \endo{G}{\HC_R(G)}$. 
L'action de $G$ sous-jacente n'est en
g{\'e}n{\'e}ral pas lisse, mais en prenant la partie lisse, on retombe sur
$V$, {\em i.e.} $\hat{V}^{\infty}=\HC_R(G)\hat{V}=V$.
En fait, un peu d'{\em abstract nonsense} montre que le 
le foncteur ``compl\'etion" de $\Mo{R}{G}$ dans $\Mo{}{RG}$ est un adjoint \`a droite du foncteur ``lissification". En clair, on a
la relation d'adjonction
$\hom{V^\infty}{W}{G}\simeq 
\hom{V}{\hat{W}}{G}$ pour tout $RG$-module $V$ et tout $RG$-module lisse $W$. 

Soit maintenant $(\pi,V)\in \Mo{R}{M}$ et
notons encore $\pi$ l'inflation de $\pi$ {\`a} $\o{P}$. Par d\'efinition on a $\ip{M,\o{P}}{G}{V}=\CC^\infty(G,V)^{\o{P}}$ o\`u $\o{P}$ agit sur une fonction $f$ par $(\o{p}f)(g)=\pi(\o{p})(f(g\o{p}))$. L'application $R$-lin\'eaire
$$\application{\int:\;}{\CC^{\infty,c}(G,V)}{\CC^{\infty}(G,V)^{\o{P}}}{f}{\left(g\mapsto \int_{\o{P}} \pi(\o{p})f(g\o{p})d\o{p}\right)}$$
induit une application $G$-\'equivariante $\CC^{\infty,c}(G,\delta_{\o{P}}^{-1}V)_{\o{P}} \To{} \CC^\infty(G,V)^{\o{P}}$, o{\`u} $\delta_{\o{P}}=\delta_P^{-1}$ d{\'e}signe le caract{\`e}re-module de
$\o{P}$. En utilisant la d\'ecom\-position de Cartan ($G=G_0\o{P}$ pour $G_0$ compact sp\'ecial), on v\'erifie que c'est un isomorphisme.

 On a donc pour tout $W\in \Mo{R}{G}$ : 
\begin{eqnarray*}
\hom{{\Ip{M,\o{P}}{G}}(V)}{W}{G} & \simeq &
%\hom{(RG\otimes_{R\o{P}}\delta_P V)^\infty}{W}{G} \\
\hom{(\CC^{\infty,c}_R(G)\otimes_R \delta_P V)_{\o{P}}}{W}{G} \\
%& \simeq & \hom{RG\otimes_{R\o{P}}\delta_PV}{\hat{W}}{G} \\
& = & \hom{\CC^{\infty,c}_R(G)\otimes_R \delta_P V}{W}{R}^{G\times \o{P}} \\
%& \simeq & \hom{\delta_PV}{\hat{W}}{\o{P}} \\
& \simeq &\hom{\delta_P V}{\hom{\CC^{\infty,c}_R(G)}{W}{R}}{R}^{G\times \o{P}} \\
& = & \hom{\delta_P V}{\hom{\CC^{\infty,c}_R(G)}{W}{G}}{\o{P}} \\
& = & \hom{\delta_P V}{\hat{W}}{\o{P}} \\
& = & \hom{\delta_P V}{\hat{W}^{\o{U}}}{M} \\
& = & \hom{\delta_P V}{(\hat{W}^{\o{U}})^\infty}{M}
\end{eqnarray*}
(Prendre garde que dans la derni{\`e}re ligne, $\hat{W}$ d{\'e}signe la
compl{\'e}tion de $W$ en tant que $G$-module et $(\hat{W}^{\o{U}})^\infty$
d{\'e}signe la partie lisse  de $\hat{W}^{\o{U}}$ en tant que $M$-module.)
Ceci montre que le foncteur $\Ip{M,\o{P}}{G}$ est adjoint \`a gauche du foncteur $W\mapsto \delta_P^{-1} \, (\hat{W}^{\o{U}})^\infty$. 
Nous allons maintenant expliciter une fl\`eche $M$-\'equivariante  $\hat{W}^{\o{U}} \To{} W_U$ fonctorielle en $W$. Pour ce faire, prenons un sous-groupe ouvert compact $U_c$ de $U$, et consid\'erons la compos\'ee
$$ \o{j}_U:\; (\hat{W}^{\o{U}})^\infty  \To{e_{U_c}} W \To{j_U} W_U. $$
La premi\`ere fl\`eche est donn\'ee par action de $e_{U_c}$ \`a gauche ; son image est dans la partie {\em lisse} $W$ de $\hat{W}$, car pour tout $w\in (\hat{W}^{\o{U}})^\infty$, on peut trouver $H=(H\cap U)H_M(H\cap \o{U})$ tel que $(H\cap U)\subset U_c$ et $H_M$ fixe $w$, qui est aussi fix\'e par $\o{U}$ donc par $\o{U}\cap H$. Par d\'efinition de $j_U$, la compos\'ee $\o{j}_U$ ci-dessus est {\em ind\'ependante du choix de $U_c$}, et par suite est $M$-\'equivariante et fonctorielle en $W$. 

Puisque la famille $\EC$ est suppos\'ee g\'en\'eratrice, pour montrer que $\o{j}_U$ est un isomorphisme, il suffit de v\'erifier que pour tout $\ve\in \EC$, la restriction de $\o{j}_U$ induit un $R$-isomorphisme $\ve (\hat{W}^{\o{U}})^\infty \simto \ve W_U$. Choisissons un point $y$ de $B(M,K)$ tel que $\ve\in RM_y$ et que l'hypoth\`ese de la proposition \ref{rescent} soit v\'erifi\'ee.
On a la factorisation
$$ \ve(\hat{W}^{\o{U}})^\infty \To{e_{U_y}} \wt\ve W \To{j_U} \ve W_U $$
et on sait d\'eja par la proposition \ref{rescent} que la fl\`eche de droite est un isomorphisme.
Pour \'etudier la fl\`eche de gauche, remarquons tout d'abord que $\hat{W}^{\o{U}} \simeq \hom{\CC^{\infty,c}_R(G)_{\o{U}}}{W}{G}$. %, toujours gr\^ace \`a l'isomorphisme $\CC^{\infty,c}_R(G/\o{U})\simeq \CC^{\infty,c}_R(G)_{\o{U}}$. 
On a donc $\ve\hat{W}^{\o{U}} \simeq \hom{\CC^{\infty,c}_R(G)_{\o{U}}\ve}{W}{G}$ et par ailleurs, on a aussi $\wt\ve W \simeq \hom{\CC^{\infty,c}_R(G)\wt\ve}{W}{G}$. Via ces identifications, on v\'erifie que la fl\`eche not\'ee $e_{U_y}$ ci-dessus devient la fl\`eche 
$$ \hom{\CC^{\infty,c}_R(G)_{\o{U}}\ve}{W}{G} \To{j_{\o{U}}^*}  \hom{\CC^{\infty,c}_R(G)\wt\ve}{W}{G}$$
duale de la fl\`eche $j_{\o{U}} : \CC^{\infty,c}_R(G)\wt\ve \To{} \CC^{\infty,c}_R(G)_{\o{U}}\ve$.
Or celle-ci, par la proposition \ref{rescent} est un isomorphisme ; en effet les groupes $U$ et $\o{U}$ jouent un r\^ole parfaitement sym\'etrique dans les hypoth\`eses, et donc dans la conclusion aussi.

\end{proof}

Avant d'\'enoncer le prochain r\'esultat, rappelons qu'une fonction $\varphi\in \CC^\infty_R(G)$ est dite cuspidale si pour tout \para $\PC=\MC\UC$ de $\GC$, on a $\int_U f(ug)du=0 $ quel que soit $g\in G$.

\begin{prop} \label{support}
Supposons que pour tout \para $\PC=\MC\UC$ de $\GC$, il existe une famille g\'en\'eratrice  d'idempotents $P$-bons de $RM$. Alors pour tout pro-$p$-sous-groupe ouvert $H$ de $G$ il existe un sous-ensemble  
    $S_H$ de $G$, compact modulo le centre et {\em ind{\'e}pendant de
      la $R$-alg\`ebre $\RC$} supportant toutes les fonctions
    cuspidales dans $\CC^c_\RC(H\ba G/H)$.
\end{prop}

\begin{proof}
Fixons pour commencer un \para $\PC=\MC\UC$. 
D'apr\`es notre hypoth\`ese, il existe un ensemble fini $\EC_H$ d'idempotents $P$-bons de $RM$ et des distributions localement constantes \`a support compact $f_\ve\in\HC_R(M)$ pour $\ve\in \EC_H$ telles qu'on ait l'\'egalit\'e 
$e_{H\cap M}= \sum_{\ve\in \EC_H} \ve f_\ve$ dans $\HC_R(M)$. Choisissons un sous-groupe ouvert compact $\o{U}_H$ de $\o{U}$ suffisamment petit pour avoir les \'egalit\'es $e_{\o{U}_H} \ve f_\ve e_H= \ve f_\ve e_H$ dans $\HC_R(G)$ pour tout $\ve\in\EC_H$, puis choisissons pour chaque $\ve\in \EC_H$ un point $y_{\ve} \in B(\MC,K)$ tel que $\o{U}_{y_\ve} \subset \o{U}_H$. Choisissons enfin un sous-groupe ouvert compact $U_H$ suffisamment grand pour contenir chaque $U_{y_\ve}$. On a alors l'\'egalit\'e dans $\HC_R(G)$ :
$$ e_{U_H} e_H = e_{U_H} \sum_{\ve\in\EC_H} \wt\ve_{y_\ve} f_\ve e_H .$$
Appliquons cette \'egalit\'e \`a un \'el\'ement $H$-invariant $v$ d'une repr\'esentation $V\in\Mo{R}{G}$ dans le noyau de la projection $j_U$ sur $V_U$. Chaque $f_\ve v $ est encore dans $\ker j_U$, donc par d\'efinition de ``$P$-bon" et par \ref{rescent}, on obtient $e_{U_H} v = 0 $.

Soit maintenant $\varphi\in\CC_R^c(H\ba G/ H)$ une fonction cuspidale. Pour tout parabolique $\PC$, on a donc $e_{U_H} *\varphi =0$. Il s'ensuit que pour tout \'el\'ement $g\in G$ tel que $U_H\subset gHg^{-1}$, on a
$$e_H * g^{-1}\varphi(1)=\int_H \varphi(gh)e_H=\varphi(g)=0.$$
Il nous reste \`a utiliser le r\'esultat suivant sur la g\'eom\'etrie de $G$ : {\em l'ensemble des $g\in G$ tels que pour tout $P$, on a $ gHg^{-1} \setminus U_H\neq \emptyset
$ est compact modulo le centre de $G$}.

\end{proof}

Pour terminer cette section, voici une \'etape facile pour construire des familles g\'en\'eratrices d'idempotents $P$-bons. N\'eanmoins, on ne s'en servira pas dans la suite.

\begin{lemme} \label{lemfamille}
Supposons que pour tout \para $\QC=\NC\VC$ contenant $\PC$, il existe une famille $\EC_Q$ d'idempotents $Q$-bons de $RN$ qui engendre la partie cuspidale de $\Mo{R}{N}$, dans le sens suivant : pour tout objet $V$ cuspidal de $\Mo{R}{N}$, il existe $\ve\in\EC$ tel que $\ve V\neq 0$. Alors il existe une famille g\'en\'eratrice $\EC$ d'idempotents $P$-bons de $RM$.
\end{lemme}

\begin{proof}
Pour tout $\QC$ et tout $\ve\in \EC_Q$, on choisit un point $y\in B(N,K)$ adapt\'e \`a $\ve$ et 
tel que $(V\cap M)_y = (V\cap M)_y^+$ et $(\o{V}\cap M)_y=(\o{V}\cap M)_y^+$.
On pose alors $\ve_M := z_{y,Q\cap M}^{-1} e_{(V\cap M)_y} e_{(\o{V}\cap M)_y^+} \ve$. C'est un idempotent de $RM_y$.
Si $x \in y+A_M$, on a
$$ \euxp\eubx \ve_M = e_{V_x^+} e_{\o{V}_x} z^{-1} \ve  \in RG_x e_{V_x} e_{\o{V}_x} \ve = RG_x \eux \eubx \ve_M $$
et de m\^eme on v\'erifie $\ve_M \eux\eubxp \in \ve_M\eux\eubx RG_x$. L'idempotent $\ve_M$ est donc $P$-bon. Appelons $\EC$ l'ensemble des $\ve_M$ obtenus en faisant varier $\QC$, $\ve\in \EC_Q$ et en saturant par $M$-conjugaison. Il nous reste \`a v\'erifier que cette famille $\EC$ est g\'en\'eratrice. Puisqu'on a satur\'e par conjugaison, il suffit de prouver que la famille des $\CC^\infty_c(M)\ve_M$ engendre la cat\'egorie $\Mo{R}{M}$, et pour cela, il suffit de montrer que pour tout objet $W$ de $\Mo{R}{M}$, il existe un $\ve_M$ tel que $\ve_M W\neq \{0\}$. Soit alors $\QC=\VC\NC$ un parabolique contenant $\PC$ et maximal pour la propri\'et\'e $W_{V\cap M}\neq 0$. Alors la repr\'esentation $W_{V\cap M}\in \Mo{R}{N}$ est cuspidale et par hypoth\`ese il existe $\ve\in \EC_Q$ tel que $\ve (W_{V\cap M})\neq 0$. Or par la proposition \ref{rescent}, on a $\ve_M W \simto \ve (W_{V\cap M})$.

\end{proof}

\section{Noeth{\'e}riannit{\'e}} \label{secnoe} \label{noether}

Dans toute cette section et sauf pr{\'e}cision suppl{\'e}mentaire, $R$ d{\'e}signe
toujours une $\zp$-alg{\`e}bre noeth{\'e}rienne. Nous allons prouver que {\em la seconde adjonction implique la noetheriannit\'e}. Commen{\c c}ons par une mise au point :

\begin{lemme}
Pour un objet $V$ de $\Mo{R}{G}$, les propri\'et\'es suivantes sont \'equivalentes :
\begin{enumerate}
        \item Pour tout pro-$p$-sous-groupe ouvert 
$H$, le $e_H\HC_R(G)e_H$-module $V^H$ est noeth{\'e}rien.
\item Tout sous-objet $W$ de $V$ engendr\'e par ses invariants sous un sous-groupe ouvert suffisamment petit est de type fini sur $G$.
\end{enumerate}
Une repr\'esentation $V$ satisfaisant ces propri\'et\'es sera dite {\em localement noeth\'erienne}.
\end{lemme}
\begin{proof}
$i) \Rightarrow ii)$ : soit $W\subset V$ et $H$ tel que $W^H$ engendre $W$. Alors $W^H\subset V^H$ est un sous-$e_H\HC_R(G)e_H$-module donc est de type fini. Or, tout sous-ensemble de $W^H$ engendrant $W^H$ sur $e_H\HC_R(G)e_H$ engendre $W$ sur $G$.

$ii)\Rightarrow i)$ : fixons $H\subset G$ pro-$p$ et ouvert, et $M$ un sous-$e_H\HC_R(G)e_H$-module de $V^H$. Le $G$-module $\la M\ra_G$ engendr\'e par $M$ satisfait $(\la M\ra_G)^H = M$. Soit $v_1,\cdots, v_n$ des g\'en\'erateurs de $\la M\ra_G$. On peut les \'ecrire $v_i=\sum g_{ij}m_{ij}$ pour des \'el\'ements $g_{ij}\in G$ et $m_{ij}\in M$ adequates. Donc les $m_{ij}$ obtenus (en nombre fini) engendrent $M$ sur $e_H\HC_R(G)e_H$.

\end{proof}

\begin{lemme} \label{cus}
  Si $V\in \Mo{R}{G}$ est {\em cuspidale} et de type fini,
   alors $V$ est localement noeth{\'e}rienne.
\end{lemme}
\begin{proof}
%(voir aussi l'appendice A.1.1 de Vign\'eras dans \cite{typmod})
(voir aussi \cite[A.1.1]{vigtypmod}). Soit $W\subset V$ un sous-objet et $H$ suffisamment petit pour que $W$ et $V$ soient respectivement engendr\'es par $W^H$ et $V^H$.
 Si  $v \in V^H$, alors par d{\'e}finition de la cuspidalit{\'e}, la fonction
$g\mapsto e_Hgv$ est {\`a} support compact modulo le centre $Z$ de $G$, donc le
$e_H\HC_R(G)e_H$-module engendr{\'e} par $v$ est de type fini sur l'anneau
$R[Z/Z\cap H]$.  Ainsi, puisque $V^H$ est de type fini sur $e_H\HC_R(G)e_H$, il est aussi de type fini sur $R[Z/Z\cap H]$. Alors par 
 noeth{\'e}riannit{\'e} de $R[Z/Z\cap H]$, le $R$-module $W^H$ est de type fini sur $R[Z/Z\cap H]$, donc {\em a fortiori} de type fini sur $e_H\HC_R(G)e_H$. Ainsi $W$ est bien de type fini sur $G$.
\end{proof}

Dans la suite de cette section, on soumet le
groupe r{\'e}ductif $\GC$ {\`a} l'hypoth{\`e}se suivante que nous appellerons (Adj) : 
\begin{center} {\em Pour tout \para $Q=NV$ de tout \levi $M$ de $G$, le foncteur $\Ip{N,\o{Q}}{M}$ est adjoint \`a gauche du foncteur $\delta_Q^{-1}\Rp{M,{Q}}{N}$.}

\end{center}
D'apr\`es le corollaire \ref{cororescent},
% et le lemme \ref{lemfamille}, 
cette hypoth\`ese est  v\'erifi\'ee si pour tout parabolique $\PC=\MC\UC$ de $\GC$, il existe une famille g\'en\'eratrice d'idempotents $P$-bons de $RM$, au sens de \ref{pbon}.
% qui engendre la partie cuspidale de $\Mo{R}{M}$.
Nous allons prouver :
\begin{prop}  \label{locnoe} Sous l'hypoth{\`e}se (Adj),
 tout $V\in \Mo{R}{G}$  de type fini est localement
noeth{\'e}rien.
\end{prop}
%On en d{\'e}duit alors l'assertion \ref{Cons} iv) en prenant
En prenant $V=\cind{H}{G}{1}$, on en d\'eduit :
\begin{coro}
(Toujours sous la m\^eme hypoth\`ese) Pour tout pro-$p$-sous-groupe ouvert, l'alg\`ebre de Hecke $e_H\HC_R(G)e_H$ est noeth\'erienne.
\end{coro}

Par ailleurs, on \'etendra dans l'appendice
\ref{decomposition}, des r\'esultats de Vign\'eras, Moy et Prasad  de d\'ecomposition de la
cat\'egorie $\Mo{R}{G}$ par le ``niveau". Ces r\'esultats montrent qu'une repr\'esenta\-tion lisse de type fini $V$ de $G$ est localement noeth\'erienne \ssi\ elle est noeth\'erienne, de sorte que :  
\begin{coro} \label{corodecomposition} Sous l'hypoth\`ese (Adj) et sous l'hypoth\`ese de validit\'e des constructions de Moy et Prasad,  la cat\'egorie $\Mo{R}{G}$ est noeth\'erienne, {\em i.e.} tout sous-objet  d'un objet de type fini est de type fini.
\end{coro}

 Comme l'{\'e}nonc{\'e} de la proposition \ref{locnoe} est clairement
 vrai pour un groupe de rang semi-simple nul, on fera l'hypoth{\`e}se de r{\'e}currence
(HR) qu'il est aussi connu pour tout sous-groupe de Levi strict de $G$. 
Pour manier cette hypoth{\`e}se de r{\'e}currence, il sera
commode d'utiliser le langage des \paras standards. On fixe donc un $K$-\para minimal $\PC_0=\MC_0\UC_0$ de $\GC$ ; les \paras standards sont ceux qui contiennent $\PC_0$ et les \levis  standards sont leurs composantes de Levi qui contiennent $\MC_0$.
Comme les \paras standards sont uniquement
d{\'e}termin{\'e}s par leur composante de Levi standard, on ommettra de
pr{\'e}ciser le parabolique dans les notations : on notera $\Rp{G}{M}$
pour $\Rp{G,P}{M}$ et $\o{\Rp{G}{M}}$ pour $\Rp{G,\o{P}}{M}$, etc...
On notera de plus $M<G$ pour ``$M$ Levi standard de $G$''.
%Ceci
%s'applique donc au moins aux groupes de rang relatif $1$.

\begin{lemme} \label{tf}
Les foncteurs paraboliques pr\'eservent la type-finitude et  le fait d'\^etre engendr\'e par ses invariants sous un sous-groupe ouvert.
\end{lemme}
\begin{proof}
Le fait que la restriction parabolique respecte la type-finitude est une cons\'equence imm\'ediate et classique de la d\'ecomposition de Cartan. Pour voir que l'induction $\Ip{M}{G}$ respecte la type finitude, rappelons qu'un objet $V\in\Mo{R}{G}$ est de type fini \ssi\ 
pour tout syst\`eme inductif filtrant (d\'enombrable) $(V_i)_{i\in I}$, l'application canonique
$ \limind \hom{V}{V_i}{G} \To{\gamma_V} \hom{V}{\limind V_i}{G} $
est surjective. Si $V = \ip{M}{G}{W}$ pour $W$  dans $\Mo{R}{M}$, alors par la seconde adjonction et le fait que $\o{\Rp{G}{M}}$ commute aux limites inductives (puisqu'il est adjoint \`a gauche de $\o{\Ip{M}{G}}$), on a un diagramme commutatif :
%pour tout syst\`eme inductif filtrant $(V_i)_{i\in I}$ dans $\Mo{R}{G}$ :
$$\xymatrix{ \hom{V}{\limind V_i}{G} \ar[r]^{\gamma_V} \ar@{=}[d] & \limind \hom{V}{V_i}{G} \ar@{=}[d] \\
\hom{W}{\delta_M (\limind \o{\Rp{G}{M}}(V_i))}{M} \ar[r]_{\gamma_W} &
  \limind \hom{W}{\delta_M \o{\Rp{G}{M}}(V_i)}{M} }
.  $$
Si $W$ est de type fini dans $\Mo{R}{M}$, l'application $\gamma_W$ est surjective donc $\gamma_V$ aussi et $V$ est de type fini.

La deuxi\`eme propri\'et\'e annonc\'ee dans l'\'enonc\'e est une cons\'equence de la premi\`ere :
v\'erifions-le pour $\Ip{M}{G}$ par exemple : un objet $W \in \Mo{R}{M}$ est engendr\'e par ses $H_M$-invariants s'il existe un $R$-module $U$ et un \'epimorphisme $\cind{H_M}{M}{1}\otimes U \To{} W$. Soit alors $H$ ouvert dans $G$ assez petit pour qu'il existe un \'epimorphisme $\cind{H}{G}{1}^n \To{} \ip{M}{G}{\cind{H_M}{M}{1}}$. Puisque $\Ip{M}{G}$ admet un adjoint \`a droite, il commute aux limites inductives et on obtient un \'epimorphisme $\cind{H}{G}{1}\otimes U^n \To{} \ip{M}{G}{W}$. 
\end{proof}

 Notons $\LC(G)$ l'ensemble fini des  \levis
standards. %(qui co{\"\i}ncide avec l'ensemble des classes de conjugaisons de
%\levis). 
On le munit d'un ordre total $\leq$ raffinant l'ordre
partiel  d{\'e}fini par l'inclusion. % $\MG \leq \MG'$ \ssi\ il existe $M,M'\in
%\MG\times \MG'$ tels que $M\subset M'$. 
On obtient donc une num\'erotation $\LC(G)=\{M_0,\cdots, M_g=G\}$ telle que $M_i\subset M_j \Rightarrow i<j$.

\begin{lemme} \label{filt}
   Soit $V\in \Mo{R}{G}$, il existe une filtration
   $$\{0\}=\FC_{-1}(V)\subseteq \FC_0(V) \subseteq \FC_1(V) \subseteq \cdots \subseteq \FC_g(V) = V $$   
 %  $(V_M)_{M\in \LC(G)}$
  fonctorielle en $V$ telle que :
   \begin{enumerate}
   \item  pour tout $i=0,\cdots, g$, le gradu\'e
   $\GC_i(V):= \FC_i(V)/\FC_{i-1}(V)$ est un quotient d'un objet de la forme $\Ip{M_i}{G}(W)$ pour un quotient {\em cuspidal} $W$ de $\delta_{M_i}\o{\Rp{G}{M_i}}(V)$. 
   \item les $\GC_i(V)$ sont  de type fini, resp. engendr\'es par leurs invariants sous un sous-groupe ouvert suffisamment petit,    si $V$ l'est.
\end{enumerate}
\end{lemme}
\begin{proof}
 Pour $M$ \levi standard de $G$, %posons $\delta_M$ le module du parabolique standard $\o{P}$ contenant $M$ et notons  
 notons
$$ G_M : =
  \coker ({\Ip{M}{G}}\circ \delta_M\o{\Rp{G}{M}} \To{adj} 1_{\Mo{R}{G}}) $$
le conoyau dans la cat\'egorie (ab\'elienne) des endofoncteurs de $\Mo{R}{G}$ de la fl{\`e}che d'adjonction entre ${\Ip{M}{G}}$
et $\delta_M \o{\Rp{G}{M}}$ conform\'ement \`a notre hypoth\`ese (Adj).
Remarquons que  la fl{\`e}che 
$\o{\Rp{G}{M}}\circ \Ip{M}{G}\circ \delta_M\o{\Rp{G}{M}}
\To{\o{\Rp{G}{M}}(adj)} \o{\Rp{G}{M}}$ est un {\'e}pimorphisme de
foncteurs, puisque par d{\'e}finition de l'adjonction et torsion par $\delta_M^{-1}$, la compos\'ee 
$ \o{\Rp{G}{M}} \To{(adj)'\o{\Rp{G}{M}}} \o{\Rp{G}{M}}\circ \Ip{M}{G}\circ \delta_M\o{\Rp{G}{M}}
\To{\o{\Rp{G}{M}}(adj)} \o{\Rp{G}{M}}$ est l'identit\'e. Ainsi, $\o{\Rp{G}{M}}\circ
G_M =0$.
Posons alors pour tout $i\geq 0$
$$ \FC_i:= \ker \left(1_{\Mo{R}{G}} \To{} G_{M_i}\circ G_{M_{i-1}} \circ
\cdots \circ G_{M_0}\right)
$$
On obtient une filtration du foncteur identit{\'e} de $\Mo{R}{G}$ dont les
quotients s'identifient 
\begin{eqnarray*}
\FC_i/\FC_{i-1} & \simeq & \ker (G_{M_{i-1}}\circ\cdots \circ G_{M_0} \To{}
G_{M_i}\circ\cdots \circ G_{M_0}) \\
& \simeq & \im(\Ip{M_i}{G}\circ \delta_{M_i}\o{\Rp{G}{M_i}}\circ G_{M_{i-1}}\circ\cdots
\circ G_{M_0} \To{adj} G_{M_{i-1}}\circ\cdots \circ G_{M_0})
\end{eqnarray*}
Puisque $G_{M_i}\circ\cdots \circ G_{M_0}$ est un quotient de
$G_{M_{i-1}}\circ\cdots \circ G_{M_0}$, on voit par une r{\'e}currence
imm{\'e}diate que
 $\o{\Rp{G}{M_j}} \circ G_{M_i}\circ \cdots \circ G_{M_0}=0$ pour tout $j\leq i$.
En particulier, $\o{\Rp{G}{M_i}}\circ G_{M_{i-1}}\circ\cdots \circ
G_{M_0}(V)$ est un $M_i$-module cuspidal pour tout $V\in\Mo{R}{G}$. %, de sorte que
%$F_{M_i}/F_{M_{i-1}}$ a son image dans la cat{\'e}gorie $\Mo{R}{G}_{\MG_i}$.
Pour un tel $V$, la filtration $\FC_*(V)$ remplit le cahier des charges  du point i)  et le point ii) en d\'ecoule par le lemme \ref{tf}.
\end{proof}

D'apr\`es ce lemme, 
pour  montrer que tout $RG$-module $V$ de type fini est localement
noeth{\'e}rien, il suffit de le faire pour $V$ de la
forme $\ip{M}{G}{W}$ avec $W$  cuspidal de type fini. 
De plus, d'apr\`es le lemme \ref{cus} on peut supposer $M\neq G$.

Soit alors $U\subset \ip{M}{G}{W}$ un sous $RG$-module engendr\'e par ses invariants sous un sous-groupe ouvert suffisamment petit ; on veut
montrer que $U$ est de type fini. 
Pour $0\leq i < g$, le gradu\'e $\GC_i(U)$ est  d'apr\`es \ref{filt} i) un quotient de $\Ip{M_i}{G}\circ \delta_{M_i}\o{\Rp{G}{M_i}}(U)$. Or, $\o{\Rp{G}{M_i}}(U)$ est un sous-objet de $\o{\Rp{G}{M_i}}(\ip{M}{G}{W})$, engendr\'e par ses invariants sous un sous-groupe ouvert suffisamment petit. 
Par l'hypoth\`ese de r\'ecurrence (HR) appliqu\'ee \`a $\delta_{M_i}\o{\Rp{G}{M_i}}(\ip{M}{G}{W})$ (qui est de type fini sur $N$) et le lemme \ref{tf}, $\GC_i(U)$ est donc de type fini.
 Il reste \`a prouver que le dernier quotient  de la filtration $\GC_g(U)$ (le quotient cuspidal) est aussi de type fini.

Nous allons utiliser le fait que $\ip{M}{G}{W}$ est muni d'une
structure suppl{\'e}mentaire : si $z$ est un {\'e}l{\'e}ment du
 centre  de $M$, il
agit par fonctorialit{\'e} de fa{\c c}on $G$-{\'e}quivariante 
sur $\ip{M}{G}{W}$ ; on note $\ip{M}{G}{z}$ cette action. 
Fixons alors un r\'eseau cocompact $Z_M$ du centre de $M$ ; on obtient ainsi une structure
de $R[Z_M]G$-module  lisse sur $V=\ip{M}{G}{W}$.

Il n'y a aucune raison
pour que $U$ soit stable par cette action de $R[Z_M]$, mais on peut consid\'erer la sous-repr\'esentation $\wt{U}:=R[Z_M].U$ de $\ip{M}{G}{W}$ somme des translat\'es de $U$ sous $Z_M$. Pour $0\leq i \leq g$, la sous-repr\'esentation $\FC_i(\wt{U})$ est stable sous $R[Z_M]$ et contient $\FC_i(U)$, par fonctorialit\'e de la filtration \ref{filt}. 
Par la seconde adjontion et par d\'efinition, la repr\'esentation $\FC_{g-1}(\wt{U})$ n'a aucun quotient cuspidal et on a donc $\FC_{g-1}(\wt{U})\cap U = \FC_{g-1}(U)$ ; en d'autre termes, $\GC_g(U)$ est un sous-objet de $\GC_g(\wt{U})$.
Ainsi d'apr\`es \ref{cus}, si on montre que $\GC_g(\wt{U})$ est de type fini, on pourra en d\'eduire que $\GC_g({U})$ est de type fini et on aura termin\'e la preuve de la proposition \ref{locnoe}.

D'apr\`es la preuve du lemme \ref{cus}, $W$ est une $R[Z_M]M$-repr\'esentation admissible. Puisque l'induction $\Ip{M}{G}$ respecte l'admissibilit\'e,  $\ip{M}{G}{W}$ est elle-aussi $R[Z_M]$-admissible. Ainsi, 
  pour tout  pro-$p$-sous-groupe ouvert $H$, les invariants $(\GC_g(\wt{U}))^H$ forment un $R[Z_M]$-module de type
  fini. Si l'on montre que la restriction de ce module  {\`a} $R[Z_G]$, o\`u $Z_G:=Z_M\cap \ZC(G)$ est  de type fini, alors en prenant $H$ suffisamment petit pour que $(\GC_g(\wt{U}))^H$ engendre $\GC_g(\wt{U})$, on en d\'eduira la type-finitude de $\GC_g(\wt{U})$ sur $G$.

On peut effectuer une r\'eduction suppl\'ementaire en pr\'esentant $W$
comme un quotient $$W_0 \otimes R[M/M^c]=
\cind{M^c}{M}{R\otimes{W_0}_{|M^c}}  \To{} W$$
o\`u $W_0 \in \Mo{\zp}{M}$ est de type fini {\em et $\zp$-admissible}.
Pour cela il suffit de choisir un ensemble fini de
$RM$-g\'en\'erateurs de $W$, de consid\'erer la
$\zp M^c$-repr\'esentation $W_c$ qu'ils engendrent (et qui est $\zp$-admissible
puisque $W$ est cuspidale et $Z_M\cap M^c=\{1\}$), d'\'etendre
celle-ci trivialement \`a $M^cZ_M$ puis de poser
$W_0:=\cind{M^cZ_M}{M}{W_c}$.
L'avantage de cette r\'eduction apparaitra dans la preuve du lemme
suivant. Remarquons seulement que $\ip{M}{G}{W_0\otimes R[M/M^c]}$ est
munie d'une action de $R[M/M^c]$ qui prolonge l'action naturelle de
$R[Z_M]$. Maintenant,  $\GC_g(\wt{U})$ s'identifie \`a un
sous-quotient $R[Z_M]$-\'equivariant de $\ip{M}{G}{W_0\otimes
  R[M/M^c]}$, mais le m\^eme argument que ci-dessus nous dit qu'il
suffit de prouver la type-finitude de son $R[M/M^c]$-satur\'e.
 
En d'autres termes, il nous suffit maintenant de  prouver :
\begin{lemme} 
  Soit $M$ \levi standard et $W_0\in \Mo{\zp}{M}$ cuspidale de
  type fini et admissible.
 Alors tout $R[M/M^c]G$-sous-quotient cuspidal de
 $\ip{M}{G}{W_0\otimes R[M/M^c]}$ est $R[Z_G]$-admissible.
\end{lemme}
\begin{proof} 
Commen{\c c}ons par fixer un pro-$p$-sous-groupe ouvert $H$ de $G$.
Si $\PG$ est un id\'eal premier de $R[M/M^c]$, nous noterons $K_\PG$
le corps r\'esiduel du localis\'e $R[M/M^c]_\PG$.

\medskip

{\em Premi\`ere \'etape : si $\PG$ est dans le support de $X^H$ pour un
$R[M/M^c]G$-sous-quotient cuspidal $X$ de $\ip{M}{G}{W_0\otimes
  R[M/M^c]}$, alors l'induite $\ip{M}{G}{W_0\otimes K_\PG}$ poss\`ede 
 un $K_\PG G$-sous-quotient cuspidal non-nul.}

Comme le foncteur  ``localisation en $\PG$''  est exact, le
$R[M/M^c]_\PG G$-module $X_\PG$ est un sous-quotient cuspidal non-nul
de $Y_\PG:=\ip{M}{G}{W_0\otimes R[M/M^c]_\PG}$. Soient $U\subset
V\subset Y_\PG $ tels que $V/U= X_\PG$. Par le lemme de Nakayama, on a
$U^H\cap \PG V^H\subsetneq V^H$ et par le lemme d'Artin-Rees, il existe un entier $k$ tel
que $\PG^k Y_\PG^H \cap V^H \subset \PG V^H$. Ainsi le quotient
$U/(\PG V+U)$ est un sous-quotient cuspidal non nul de
$Y_\PG/\PG^k=\ip{M}{G}{W_0\otimes R[M/M^c]_\PG/\PG^k}$.

Posons $Y_k := Y_\PG/\PG^k$ et fixons $U_k\subset V_k\subset Y_k$ de
quotient non-nul et cuspidal. On a une suite exact 
$$ (V_k \cap \PG Y_k)/(U_k \cap \PG Y_k) \injo V_k/U_k
\twoheadrightarrow (V_k + \PG Y_k)/(U_k +\PG Y_k) $$
dans laquelle le terme de gauche est un sous-quotient de $\PG Y_k$ et
celui de droite est un sous-quotient de $Y_1=Y_k/\PG$.
 Choisissons un syst\`eme de $r$
g\'en\'erateurs ($r\in \NM$) de l'id\'eal $\PG$ ; il lui est associ\'e
un \'epimorphisme $Y_{k-1}^r\To{} \PG Y_k$ pour chaque $k>1$. On voit
donc par r\'ecurrence descendante que $Y_1$ admet un sous-quotient
cuspidal non-nul.

\medskip

{\em Deuxi\`eme \'etape : si l'induite $\ip{M}{G}{W_0\otimes K_\PG}$ poss\`ede 
 un $K_\PG G$-sous-quotient cuspidal non-nul, alors $K_\PG$ est de
 caract\'eristique positive $l\neq p$, et il existe une
 repr\'esentation irr\'eductible $\sigma\in \Mo{\FM_l}{M}$ telle que
 l'induite $\ip{M}{G}{\sigma\otimes_{\FM_l} K_\PG}$ poss\`ede 
 un $K_\PG G$-sous-quotient cuspidal non-nul. }

D'apr\`es la th\'eorie classique de Bernstein, on sait que sur un
corps de caract\'eristique nulle, une induite parabolique n'a pas de
sous-quotient cuspidal, d'o\`u l'assertion sur la
caract\'eristique. Le reste est \'evident puisque une repr\'esentation
cuspidale admissible de type fini sur un corps est de longueur finie.

\medskip

{\em Troisi\`eme \'etape : si  $\ip{M}{G}{\sigma\otimes_{\FM_l} K_\PG}$ poss\`ede 
 un $K_\PG G$-sous-quotient cuspidal non-nul, alors pour tout parabolique $Q$ contenant $M$, l'induite $\ip{M,Q}{G}{\sigma\otimes_{\FM_l} K_\PG}$ poss\`ede 
 un $K_\PG G$-sous-quotient cuspidal non-nul.
} 

Notons $\wt{R}:=R/(\PC\cap R)$ ; c'est une $\FM_l$-alg\`ebre  int\`egre et on note $\wt{K}$ son corps des fractions. Par hypoth\`ese, l'induite $\ip{M}{G}{\sigma\otimes \wt{R}[M/M^c]}$ a un sous-quotient cuspidal dont le support contient $\PG$. Par l'argument de \cite[7.3]{nutempere}, la repr\'esentation $\sigma\otimes \wt{K}(M/M^c)$ est $(P,Q)$-r\'eguli\`ere au sens de \cite[2.10]{nutempere}, et 
il existe donc un op\'erateur d'entrelacement {\em non-nul} :
$$ J_{Q|P}(\sigma\otimes \wt{K}(M/M^c)) : \ip{M}{G}{\sigma\otimes \wt{K}(M/M^c)} \To{} \ip{M,Q}{G}{\sigma\otimes \wt{K}(M/M^c)}. $$
Cet op\'erateur est {\em injectif}, par le th\'eor\`eme d'irr\'eductibilit\'e g\'en\'erique \cite[5.1]{nutempere} que l'on peut appliquer gr\^ace au lemme \ref{nudiscus} ci-dessous. En effet, soit $\o{K}$ une cl\^oture alg\'ebrique de $\wt{K}$. Choisissons une d\'ecomposition de la repr\'esentation $\sigma\otimes\o\FM_l$ en une somme directe de repr\'esentations $\o\sigma_i\in\Mo{\o\FM_l}{M}$ irr\'eductibles. Chaque $\o\sigma_i\otimes_{\o\FM_l}\o{K}(M/M^c)$ est $(P,Q)$-r\'eguli\`ere et lui est donc associ\'e un op\'erateur d'entrelacement. Par d\'efinition, l'op\'erateur $J_{Q|P}(\sigma\otimes \wt{K}(M/M^c))$ induit par extension des scalaires de $\wt{K}$ \`a $\o{K}$ la somme directe des op\'erateurs $J_{Q|P}(\o\sigma_i\otimes \o{K}(M/M^c))$. Maintenant par \cite[5.1]{nutempere} chaque $\ip{M}{G}{\o\sigma_i\otimes \o{K}(M/M^c)}$ est irr\'eductible, d'o\`u l'injectivit\'e cherch\'ee (l'op\'erateur $J_{Q|P}$ \'etant non-nul).

Puisque la repr\'esentation $\ip{M}{G}{\sigma\otimes\wt{R}[M/M^c]}$ est de type fini, il existe un \'el\'ement $f\in \wt{R}[M/M^c]$ tel que l'op\'erateur  $f.J_{Q|P}$ induise le $\wt{R}[M/M^c]G$-morphisme {\em injectif} :
$$ f.J_{Q|P}(\sigma\otimes \wt{K}(M/M^c)) : \ip{M}{G}{\sigma\otimes \wt{R}[M/M^c]} \To{} \ip{M,Q}{G}{\sigma\otimes \wt{R}[M/M^c]}. $$
On en d\'eduit l'existence d'un sous-quotient cuspidal de support contenant $\PG$ dans l'induite $\ip{M,Q}{G}{\sigma\otimes \wt{R}[M/M^c]}$ et par la premi\`ere \'etape, l'existence d'un sous-quotient cuspidal non-nul dans l'induite $\ip{M,Q}{G}{\sigma\otimes K_\PG}$.

\medskip

{\em Quatri\`eme \'etape : si $R[M/M^c]/\PG$ n'est pas  fini
  sur $R[Z_G]$ alors il existe une valuation $\nu :K_\PG \To{}
  \RM$ et un $m \in M\cap [G,G]$ tels que $\nu(\psi(m))\neq 0$, o\`u 
  %$\o{m}$ d\'esigne l'image de $m$ dans $K_\PG$.}
  $\psi: M \To{} K_\PG^\times$ d\'esigne le caract\`ere tautologique.}

Pour abr\'eger, notons $R_G$ l'image de $R[Z_G]$ dans $R[M/M^c]/\PG$
et $K_G$ son corps de fractions.
L'anneau $R[M/M^c]/\PG$ est fini
sur l'anneau engendr\'e par $R_G$ et  $\psi(M\cap [G,G])$.
Soit $d$ le degr\'e de transcendance de $K_\PG$ sur $K_G$. Deux cas se pr\'esentent :

Si $d>0$, alors il existe $m\in M\cap [G,G]$ tel que $\psi(m)$ soit transcendant sur $K_G$. On pose $\nu(\psi(m))=1$ et on \'etend de mani\`ere arbitraire la valuation obtenue \`a $K_\PG$.

 Si $d=0$, notons $\wt{R}_G$ la cl\^oture int\'egrale de $R_G$  dans $K_\PG$. 
 Puisque $R[M/M^c]/\PG$ est de type fini sur $R_G$, il n'est pas inclus dans $\wt{R}_G$ (sinon, il serait fini puisqu'entier et contredirait notre hypoth\`ese).
L'anneau $\wt{R}_G$ n'est pas n\'ecessairement noeth\'erien, mais il est ``de Krull" \cite[1.4
  Corollaire]{AC7}, 
donc est l'intersection des localis\'es en ses id\'eaux premiers de hauteur $1$. On peut donc trouver un \'el\'ement $m\in M\cap[G,G]$ et un tel localis\'e ne contenant pas $\psi(m)$.
Or ce localis\'e est normal lui aussi  donc est un anneau de valuation discr\`ete. La valuation associ\'ee est n\'ecessairement non-nulle sur $\psi(m)$ et s'\'etend au corps de fractions qui par hypoth\`ese n'est autre que $K_\PG$.

 \medskip 
  
 {\em Fin de la preuve :} Supposons que $R[M/M^c]/\PG$ n'est pas fini sur $R[Z_G]$ et choisissons une valuation de $\KC_\PG$ comme dans l'\'etape pr\'ec\'edente. 
Ainsi la compos\'ee $\nu \circ \psi$ est un \'el\'ement non nul de l'espace vectoriel ${a_M^G}^*:= \hom{M/M^c}{\RM}{}/\hom{G/G^c}{\RM}{}$.
Nous renvoyons \`a \cite[2.2]{nutempere} pour la d\'efinition des chambres de Weyl $(a_Q^*)^+$ dans ${a_M^G}^*$ associ\'ees aux paraboliques $Q$ dont $M$ est une composante de Levi. Il existe un $Q$ tel que $\nu\circ \psi$ soit dans l'adh\'erence $\o{(a_Q^*)^+}$ de $(a_Q^*)^+$. Cette adh\'erence est r\'eunion disjointe \cite[2.3]{nutempere} de chambres de Weyl de paraboliques $O$ contenant $Q$, {\em i.e.}
$$\o{(a_Q^*)^+} = \bigsqcup_{Q\subseteq O \subseteq G} (a_O^*)^+. $$
Soit $O$ l'unique parabolique contenant $Q$ tel que $\nu\circ \psi\in (a_O^*)^+$. Puisque $\nu\circ \psi\neq 0$, on a $O\neq G$. Soit $N$ sa composante de Levi contenant $M$. Alors 
d'apr\`es le lemme \ref{nudiscus} ci-dessous et  \cite[3.16]{nutempere} joint \`a l'argument de \cite[7.3]{nutempere}, les repr\'esentations $\ip{M,Q}{G}{\sigma\otimes \o{K_\PG}}$ et $\ip{M,N\cap Q}{N}{\sigma\otimes \o{K_\PG}}$ ont la m\^eme longueur ($\o{K_\PG}$ d\'esigne une cl\^oture alg\'ebrique de $K_\PG$); plus pr\'ecis\'ement, tout sous-quotient irr\'eductible de la seconde s'induit irr\'eductiblement par le foncteur $\Ip{N,O}{G}$. Il s'ensuit que l'induite 
$\ip{M,Q}{G}{\sigma\otimes K_\PG}$ {\em ne peut pas avoir de $K_\PG G$-sous-quotient cuspidal}, ce qui contredit l'\'enonc\'e de l'\'etape 3.  
  Cela conclut la preuve du lemme et de la proposition \ref{locnoe}. 
\end{proof}

Dans la preuve pr\'ec\'edente, nous avons \`a deux reprises fait appel \`a certains r\'esultats de \cite{nutempere}, comme par exemple la propri\'et\'e {\em d'irr\'eductibilit\'e g\'en\'erique}. Dans \cite{nutempere}, ceux-ci sont \'enonc\'es sous la condition que $G$ poss\`ede un sous-groupe discret cocompact, mais cette condition ne sert qu'\`a assurer qu'une certaine propri\'et\'e ($\nu$-discret implique cuspidal, {\em cf} ci-dessous) soit vraie.
Le lemme suivant montre comment l'hypoth\`ese (Adj) implique directement cette propri\'et\'e.

\begin{lemme} \label{nudiscus}
Soit $\KC$ un $R$-corps et $\nu$ une valuation discr\`ete de $\KC$. Soit $(\pi,V)$ une repr\'esentation admissible dans $\Mo{\KC}{G}$ dont les coefficients matriciels {\em tendent essentiellement vers $0$ \`a l'infini}  pour la norme associ\'ee \`a $\nu$ ({\em cf} \cite[3.18]{nutempere} o\`u une telle repr\'esentation est dite $\nu$-discr\`ete). Alors $\pi$ est cuspidale.
\end{lemme}
\begin{proof}
Soit $\OC$ l'anneau de la valuation $\nu$ et $\varpi$ une uniformisante.
Par \cite[Prop 6.3]{nutempere}, il existe un sous-$\OC$-module $G$-stable et $\OC$-admissible $\omega\subset V$ qui engendre $V$ sur $\KC$. Puisque les coefficients matriciels tendent essentiellement vers $0$, il en est de m\^eme des application $f_v :\; g\in G \mapsto e_Hgv \in \omega$ pour $v\in\omega$ et $H$ prop-$p$-sous-groupe ouvert, en le sens suivant : 
$$\forall n\in\NM,  f_v^{-1}(\omega \setminus \varpi^n \omega) \hbox{ est compact modulo le centre}. $$
En d'autres termes chaque, $\omega/\varpi^n\omega$ est cuspidal. Puisque les foncteurs $\Rp{G}{M}$ ont des adjoints \`a gauche, il commutent aux limites projectives et par cons\'equent la limite $\limp{n} \omega/\varpi^n\omega$ est cuspidale elle-aussi.
Or, par admissibilit\'e la fl\`eche canonique $\omega \To{} \limp{n} \omega/\varpi^n\omega$ est injective et $\omega$ et donc $\pi$ sont cuspidales.

\end{proof}

\def\dag{\dagger}
\def\la{\langle}
\def\ra{\rangle}

\section{Mod\`eles entiers lisses} \label{mod} 
\def\TC{{\mathcal T}}
\def\UC{{\mathcal U}}
\def\SC{{\mathcal S}}

\def\u{\underline}

On note toujours $\GC$ un groupe r\'eductif connexe sur $K$. Nous commen{\c c}ons par donner les \'enonc\'es principaux de cette section, puis nous passerons aux preuves.

\alin{D\'efinitions et principaux \'enonc\'es} \label{defmod}
Soit $\u\GC$ un mod\`ele lisse et connexe de $\GC$ sur $\OC_K$, {\em i.e.} un sch\'ema en groupes lisse sur $\OC_K$ \`a fibres connexes, et muni d'une identification de sa fibre g\'en\'erique avec $\GC$. Pour un sous-groupe ferm\'e $\HC$ de $\GC$ nous noterons $\u\HC$ son adh\'erence sch\'ematique dans $\u\GC$ ; ses points entiers sont donc donn\'es par $\u\HC(\OC_K)=\HC(K)\cap \u\GC(\OC_K)$.

Si $\PC$ est un \para de $\GC$, nous dirons qu'il est {\em $\u\GC$-admissible} s'il admet une composante de Levi $\MC$ de la forme $\ZC_\GC(\SC)$ pour un sous-tore d\'eploy\'e  $\SC$ de $\GC$ qui se prolonge en un sous-tore de $\u\GC$ (un tel prolongement est  n\'ecessairement ferm\'e par \cite[Exp. VIII. Cor 5.7]{SGA3} donc \'egal \`a $\u\SC$). Le sous-groupe de Levi $\MC$ est alors dit lui aussi {\em $\u\GC$-admissible}, et son adh\'erence $\u\MC$ est le centralisateur $\ZC_{\u\GC}(\u\SC)$ de $\u\SC$ dans $\GC$ ; celui-ci est lisse sur $\OC_K$ par \cite[Exp. XI, Cor. 5.3]{SGA3} et connexe puisque les centralisateurs de tores dans les groupes connexes sur un corps sont connexes.

Nous utiliserons souvent la technique \'el\'egante de {\em dilatation} d\^ue \`a Raynaud, \cite[3.2]{BLR}, et introduite dans le pr\'esent contexte par Yu \cite{Yumodel}.
Soit $\u\GC^\dag \To{\nu_{\u\GC}} \u\GC$ la dilatation dans $\u\GC$ du radical unipotent ${^u\u\GC_k}$ de la fibre sp\'eciale de $\u\GC$. Par d\'efinition d'une dilatation, on sait que $\nu_{\u\GC}$  induit un isomorphisme des fibres g\'en\'eriques $\u\GC^\dag_K\simto \u\GC_K\simeq \GC$ et que $\u\GC^\dag(\OC_K)=\{g\in\u\GC(\OC_K),\;g\hbox{ mod }\varpi \in {^u\GC_k}(k)\}$. On sait aussi que $\u\GC^\dag$ est un sch\'ema en groupes lisse \cite[3.2 Prop 3]{BLR}, donc un mod\`ele lisse de $\GC$. 
Convenons de noter $\u\HC^\dag$ l'adh\'erence sch\'ematique dans $\u\GC^\dag$ d'un sous-groupe ferm\'e $\HC$ de $\GC$. La  restriction de $\nu_{\u\GC}$ induit un morphisme $\u\HC^\dag\To{} \u\HC$ qui s'identifie \`a la dilatation dans $\u\HC$ de $\u\HC_k\cap{^u\u\GC_k}$.

Dans le cas o\`u $\HC=\MC$ est un Levi admissible, alors $\u\MC_k$ est de la forme $\ZC_{\u\GC_k}(\u\SC_k)$ pour un tore $\SC_k$, et on a l'\'egalit\'e ${^u\u\MC_k}=\u\MC_k\cap {^u\u\GC_k}$ d'apr\`es \cite[Exp XIX, 1.3]{SGA3}.
 Le morphisme  $\u\MC^\dag\To{\nu_{\u\GC}} \u\MC$ co\"incide donc avec le morphisme de dilatation $\nu_{\u\MC}$ du radical unipotent de la fibre sp\'eciale de $\u\MC$.

Notons maintenant avec des lettres droites les ensembles de points entiers, par exemple $\u{H}=\u\HC(\OC_K)$. Notons que $\u{H}$ est pro-$p$ \ssi\ ${^u\u\HC_k}$ est unipotent. En particulier,
$\u{H}^\dag:=\u{\HC}^\dag(\OC_K)$ est pro-$p$.

\begin{lemme} \label{defniv0}
Soit $R$ un anneau commutatif unitaire o\`u $p$ est inversible.
Pour un idempotent central $\ve$ de $R\u{G}$, les propri\'et\'es  suivantes sont \'equivalentes :
\begin{enumerate} 
\item $\ve\in R\u{G}^\dag e_{\u{U}^\dag}  R\u{G}^\dag$  pour tout \para $\u\GC$-admissible $\PC=\MC\UC$ de $\GC$. 
\item $\ve\in R\u{G}^\dag e_{\u{U}^\dag} e_{\o{\u{U}}^\dag} R\u{G}^\dag$  pour une paire de sous-groupes paraboliques oppos\'es $\u\GC$-admissibles minimaux   $(\PC,\o\PC)$ de $\GC$
\end{enumerate}
Un idempotent satisfaisant ces propri\'et\'es sera dit {\em essentiellement de niveau z\'ero}.
\end{lemme}

\noindent{\sl Remarque :} L'expression $R\u{G}^\dag e_{\u{U}}^\dag R\u{G}^\dag$ d\'esigne le {\em $R$-module engendr\'e} par les distributions du type $\varphi*e_\Omega*\psi$ o\`u $\varphi,\psi\in R\u{G}^\dag$. Cette remarque s'applique \`a toutes les expressions de ce genre que l'on rencontrera dans la suite de ce paragraphe.

Remarquons aussi que l'idempotent associ\'e \`a un caract\`ere lisse $\theta:\u{G}^\dag \To{} R^\times$ normalis\'e par $\u{G}$ est essentiellement de niveau z\'ero \ssi\ $\theta_{|\u{U}^\dag}$ et $\theta_{|\o{\u{U}}^\dag}$ sont triviaux pour au moins une paire de  paraboliques oppos\'es $\u\GC$-admissibles minimaux. Le lemme suivant \'etudie une situation plus g\'en\'erale, que l'on rencontre souvent dans la th\'eorie des types.
\begin{prop} \label{constess}
Soit $\u{G}^*$ un sous-groupe ouvert normal de $\u{G}$ tel que $[\u{G}^\dag,\u{G}^\dag]\subseteq \u{G}^*\subseteq \u{G}^\dag$ et $\theta: \u{G}^*\To{} R^\times$ un caract\`ere lisse normalis\'e par $\u{G}$. On suppose que pour une paire $(\PC,\o\PC)$ de \paras oppos\'es $\u\GC$-admissibles minimaux, on a
\begin{enumerate}
        \item $\theta_{|\u{U}^*}$ et $\theta_{|\o{\u{U}}^*}$ sont triviaux (o\`u l'$*$ d\'esigne l'intersection avec $\u{G}^*$).
        \item L'accouplement $\application{}{\u{U}^\dag/\u{U}^*\times \o{\u{U}}^\dag/\o{\u{U}}^*}{R^\times}{(u,v)}{\theta([u,v])}$ est non-d\'eg\'en\'er\'e.
\end{enumerate}
Alors l'idempotent central $[\theta]$ de $R\u{G}^\dag$ associ\'e \`a $\theta$ est essentiellement de niveau z\'ero.
\end{prop}

Dans les exemples que l'auteur connait, le groupe $\u{G}^*$ est le groupe des points entiers d'un mod\`ele lisse connexe $\u\GC^*$ obtenu par dilatation dans $\u\GC$ d'un sous-groupe normal compris entre $[{^u\u\GC_k},{^u\u\GC_k}]$ et  $^u\u\GC_k$. 

Si $\MC$ est un \levi $\u\GC$-admissible de $\GC$, on a vu que $\u\MC$ est un mod\`ele lisse connexe de $\MC$. On peut donc appliquer \`a $\u\MC$ la notion d'idempotent essentiellement de niveau $0$ de $R\u{M}$ d\'efinie en \ref{defniv0}. %, laquelle est intrins\`eque par \ref{mod1}.

Le th\'eor\`eme suivant est une g\'en\'eralisation du th\'eor\`eme principal de \cite{HL}.

\begin{theo} \label{theomod}
Soit $\u\GC$ un mod\`ele lisse et connexe de $\GC$ et $(\PC=\MC\UC,\o\PC=\MC\o\UC)$ une paire de paraboliques $\u\GC$-admissibles oppos\'es de $\GC$. Pour tout idempotent central essentiellement de niveau z\'ero $\ve$ de $R\u{M}$, on a
$$ e_{\u{U}^\dag}e_{\o{\u{U}}} \ve \in R\u{G} e_{\u{U}} e_{\o{\u{U}}} \ve.$$ 
\end{theo}

Ce r\'esultat g\'en\'eral ne sera  toutefois pas suffisant pour les applications. On se donne
 maintenant  un autre mod\`ele lisse connexe $\u\GC'$ de $\GC$ ainsi qu'un morphisme $\u\GC'\To{\varphi} \u\GC$ de $\OC_K$-sch\'emas en groupes  dont la fibre g\'en\'erique est un isomorphisme.
\begin{lemme} \label{lemmod}
Le noyau de $\varphi_k: \u\GC'_k\To{} \u\GC_k$ est unipotent, lisse et connexe.
\end{lemme}
 
 En particulier, si $\u\SC'$ est un sous-tore  de $\u\GC'$, alors $\varphi_{|\u\SC'}$ est un monomorphisme  (et donc une immersion ferm\'ee par \cite[Exp VIII Cor. 5.7]{SGA3}) de $\u\SC'$ dans $\u\GC$. Il s'ensuit que si un \para de $\GC$ est $\u\GC'$-admissible, alors il est $\u\GC$-admissible. Nous d\'ecorerons d'un $'$ les objets d\'efinis comme ci-dessus relatifs \`a $\u\GC'$.  Par exemple, ${\u\GC'}^\dag$ d\'esignera la dilatation du radical unipotent de $\u\GC'$. 
Notons que d'apr\`es le lemme ci-dessus, on a ${^u\u\GC'_k}=\varphi_k^{-1}({^u\im(\varphi_k)})$.
On en d\'eduit aussi les \'equivalences suivantes :
$$ 
\left(\u{G}'\cap \u{G}^\dag\subseteq \u{G'}^\dag  \right) \equ \left(\varphi^{-1}({^u\u\GC_k})\supseteq{^u\u\GC'_k}\right) \equ  \left( \im(\varphi_k) \cap {^u\u\GC_k} \hbox{ connexe} \right)$$
$$ \left(\u{G}'\cap \u{G}^\dag = \u{G'}^\dag\right) \equ \left(\varphi^{-1}({^u\u\GC_k})={^u\u\GC'_k}\right) \equ  \left({^u(\im\varphi_k)} = \im\varphi_k  \cap {^u\u\GC_k}\right) $$

La situation un peu alambiqu\'ee \'etudi\'ee dans le r\'esultat suivant est, encore une fois, motiv\'ee par les exemples connus de la th\'eorie des types.

\begin{coro} \label{utile}
Gardons les notations ci-dessus avec l'hypoth\`ese $\u{G}'\cap \u{G}^\dag \subseteq \u{G'}^\dag$. Soit $\ve'$ un idempotent central essentiellement de niveau z\'ero de $R\u{M}'$, et supposons qu'il existe un idempotent $\wt\ve'$ de $R{\u{G}'}^\dag$, centralis\'e par $\u{G}'$ et
satisfaisant les propri\'et\'es suivantes :
\begin{enumerate}
        \item $e_{\u{\o{U}'}^\dag} \wt\ve' e_{{\u{U}'}^\dag}=e_{\u{\o{U}'}^\dag} \ve' e_{{\u{U}'}^\dag}$ et  $e_{{\u{U}'}^\dag} \wt\ve' e_{{\u{\o{U}}'}^\dag}=e_{{\u{U}'}^\dag} \ve' e_{{\u{\o{U}}'}^\dag}$. 
        \item l'ensemble d'entrelacement $\hbox{Int}_{\u{U}}(\wt\ve'):= \{u\in \u{U}, \wt\ve' u\wt\ve'\neq 0\}$ est contenu dans (donc \'egal \`a) $\u{U}'$.
        \end{enumerate}
 Alors $$e_{\u{U}^\dag} e_{\o{\u{U}}} \ve' \in {R \u{G}}e_{\u{U}}e_{\o{\u{U}}}\ve'.$$
\end{coro}

\alin{Les preuves}
Rappelons que si $\u\SC$ est un $\OC_K$-tore, alors tout morphisme  $ \u\SC_k \To{} \u\GC_k$ se rel\`eve de mani\`ere unique \`a conjugaison pr\`es en un morphisme $\u\SC \To{} \u\GC$, voir \cite[Exp IX, Thm 3.6]{SGA3}. En particulier les sous-tores d\'eploy\'es maximaux de la fibre sp\'eciale de $\u\GC_k$ se rel\`event en des sous-tores d\'eploy\'es de $\u\GC$, n\'ecessairement ferm\'es par \cite[Exp VIII, Cor. 5.6]{SGA3}. Leurs fibres g\'en\'eriques sont des tores d\'eploy\'es, mais g\'en\'eralement {\em pas maximaux} dans $\GC$. Fixons un tel sous-tore $\u\SC$ et notons $\SC$ sa fibre g\'en\'erique, $X^*(\SC)$ son r\'eseau de caract\`eres et $V:=X^*(\SC)\otimes \RM$. Posons aussi $\Phi:= 
\Phi(\u\SC_k,\hbox{Lie}(\u\GC_k))= \Phi(\SC,\hbox{Lie}(\GC))\subset V$,
 l'ensemble des poids non-nuls  de $\SC$ dans l'alg\`ebre de Lie de $\GC$.
C'est un ensemble sym\'etrique puisque $\GC$ est r\'eductif.
Reprenant la terminologie de Bruhat-Tits \cite[1.1.2]{BT2}, une demi-droite ouverte $\wt\alpha$ de $V$ 
contenant une racine de $\Phi$ sera appel\'ee {\em rayon radiciel de $\Phi$}. L'ensemble de ces rayons radiciels  sera not\'e $\wt\Phi$, et l'image d'un sous-ensemble $\Omega$ de $\Phi$ dans $\wt\Phi$ sera not\'ee $\wt\Omega$.

\`A tout rayon radiciel $\wt\alpha\in \wt\Phi$ est associ\'e un $K$-sous-groupe alg\'ebrique ferm\'e  connexe $\UC_{\wt\alpha}$ de $\GC$, dont les poids de $\SC$ dans l'alg\`ebre de Lie appartiennent \`a $\wt\alpha$, et qui est maximal pour ces propri\'et\'es {\em cf} \cite[1.1.3]{BT2}.
On  prolonge cette d\'efinition \`a $\wt\alpha=0$ en posant $\UC_0:=\ZC_\GC(\SC)$ le centralisateur de $\SC$ dans $\GC$.
Plus g\'en\'eralement, si $\Omega\subset \Phi\cup\{0\}$, on note $\UC_\Omega$ le sous-groupe ferm\'e de $\GC$ engendr\'e par les $\UC_{\wt\alpha}$ pour $\wt\alpha\in\wt\Omega$. 

Rappelons qu'un sous-ensemble $\Omega$ de $\Phi$ est dit
 {\em parabolique} s'il est l'intersection de $\Phi$ avec un demi-espace ferm\'e,
et {\em unipotent} s'il est l'intersection de $\Phi$ avec un demi-espace ouvert.
Si $v^*\in V^*$ est une forme lin\'eaire sur $V$, la composante de Levi du sous-ensemble parabolique $\Omega(v^*):=\{\alpha\in \Phi,\; v^*(\alpha)\geq 0\}$  est par d\'efinition $\Omega^0(v^*):=\{\alpha\in \Phi,\; v^*(\alpha)= 0\}=\Omega\cap (-\Omega)$, tandis que sa composante unipotente est $\Omega^+(v^*):=\{\alpha\in \Phi,\; v^*(\alpha)> 0\}$ et sa composante unipotente {\em oppos\'ee} est  $\Omega^-(v^*):=-\Omega^+=\Phi\setminus \Omega$. 
Le groupe $\PC_\Omega:=\UC_{\Omega\cup\{0\}}$ est un sous-groupe parabolique de $\GC$ dont le radical unipotent est $\UC_{\Omega^+}$ et la composante de Levi est $\MC_\Omega:=\UC_{\Omega^0\cup\{0\}}$. L'oppos\'e de $\PC_\Omega$ par rapport \`a $\MC_\Omega$ est $\PC_{-\Omega}$ dont le radical unipotent est $\UC_{\Omega^-}$.

\begin{fact} \label{std}
Un sous-groupe parabolique $\PC$ de $\GC$ est $\u\GC$-admissible \ssi\ il est conjugu\'e par un \'el\'ement de $\u\GC(\OC_K)$ \`a un groupe de la forme $\PC_\Omega$ pour un sous-ensemble parabolique $\Omega\subset \Phi$.
\end{fact}
\begin{proof}
V\'erifions tout d'abord que $\MC_\Omega$ est $\u\GC$-admissible (et donc $\PC_\Omega$ aussi).
Soit $\SC_\Omega
:=(\cap_{\alpha\in\Omega^0}\ker\alpha)^\circ$. En tant que sous-tore de $\SC$, il se prolonge  en un sous-tore d\'eploy\'e de $\u\SC$, et ses poids dans $\hbox{Lie}(\MC_\Omega)$, resp. dans $\hbox{Lie}(\UC_{\Omega^\pm})$, sont par construction tous nuls, resp. tous non-nuls. Puisque $\ZC_\GC(\SC_\Omega)$ est lisse et connexe, et puisque $\hbox{Lie}(\GC)=\hbox{Lie}(\UC_{\Omega^+})+\hbox{Lie}(\MC_\Omega)+\hbox{Lie}(\UC_{\Omega^-})$, on a $\ZC_\GC(\SC_\Omega)=\MC_\Omega$ et $\MC_\Omega$ est $\u\GC$-admissible.

R\'eciproquement, soit $\PC$ un parabolique $\u\GC$-admissible, et $\MC$ une composante de Levi de la forme $\ZC_\GC(\TC)$ pour un tore $\TC$  d\'eploy\'e  qui se prolonge en un sous-tore d\'eploy\'e $\u\TC$ de $\u\GC$. Comme la fibre sp\'eciale de $\u\TC$ est conjugu\'ee par un \'el\'ement de $\u\GC_k(k)$ \`a un sous-tore de $\u\SC_k$, il r\'esulte de \cite[Exp. XI, Thm 5.2 bis]{SGA3} que $\TC$ est lui-m\^eme conjugu\'e par un \'el\'ement de $\u\GC(\OC_K)$ \`a un sous-tore de $\SC$. On peut donc supposer $\MC\supset \SC$. On a alors $\PC=\PC_\Omega$ pour $\Omega:=\Phi(\SC,\hbox{Lie}(\PC))$.

\end{proof}

Les paraboliques de la forme $\PC_\Omega$ seront dits {\em semi-standards} (relativement au choix du tore d\'eploy\'e maximal $\u\SC$ de $\u\GC$).

\begin{fact} \label{BT} (Bruhat-Tits) Soit $(\PC,\o\PC)$ une paire de \paras de $\GC$ de composante de Levi commune $\MC$ suppos\'ee $\u\GC$-admissible. Alors
\begin{enumerate}
        \item le morphisme produit induit un isomorphisme $\u\UC\rtimes \u\MC \To{} \u\PC$.
        \item le morphisme produit induit une immersion ouverte $\u\CC:=\u\UC \times \u\MC \times \o{\u\UC}\To{\mu} \u\GC$. En particulier, $\u\UC$ et $\o{\u\UC}$ sont lisses et connexes. De plus l'image de $\mu$ est le compl\'ementaire d'un diviseur.
        \item le morphisme produit induit un isomorphisme $({^u\u\GC_k}\cap \u\UC_k)\times ({^u\u\GC_k}\cap \u\MC_k)\times({^u\u\GC_k}\cap \o{\u\UC_k})\simto {^u\u\GC_k}$, et de plus on a ${^u\u\GC_k}\cap \u\MC_k = {^u\u\MC_k}$.
  \item  Supposons $\PC$ et $\o\PC$ semi-standards et soit  $\QC$  un autre parabolique $\u\GC$-admissible semi-standard, de radical unipotent $\VC$. Alors l'application produit $(\u{V}\cap \u{U})\times(\u{V}\cap \u{M})\times (\u{V}\cap\u{\o{U}}) \To{} \u{V}$ est bijective.
\end{enumerate}

\end{fact}

\begin{proof}
Compte tenu de \ref{std}, les deux premiers points sont donn\'es par \cite[2.2.3 ii) et iii)]{BT2}. Pour le  troisi\`eme point, %le groupe  ${^u\u\GC_k}\cap \u\MC_k$ est un sous-groupe unipotent normal de $\u\MC_k$ donc est contenu dans ${^u\u\MC_k}$. l'autre inclusion et 
la d\'ecomposition de ${^u\u\GC_k}$ est donn\'ee par \cite[1.1.11]{BT2} et la relation entre les radicaux unipotents a d\'eja \'et\'e mentionn\'ee plus haut. %Pour le point iv), on se ram\`ene au cas $\PC=\PC_\Omega$ et $\QC=\PC_\Theta$,  on remarque que $\u\UC\cap \u\VC= \u\UC_{\Omega^+\cap \Theta^+}$, $\u{\o\UC}\cap \u\VC= \u\UC_{\Omega^-\cap \Theta^+}$, et $\u\MC\cap \u\VC= \u\UC_{\Omega^0\cap \Theta^+}$, et on applique 
Une fois choisis des ensembles $\Omega$ et $\Theta$ tels que $\PC=\PC_\Omega$ et $\QC=\PC_\Theta$, le point iv) est donn\'e par \cite[2.2.3 i)]{BT2}.
\end{proof}

Si dans le point iv) on prend $\PC$ minimal, de sorte que $\VC\cap \MC=\{1\}$, on obtient pour les idempotents associ\'es $e_{V^\dag}= e_{V^\dag \cap U} e_{V^\dag\cap\o{U}}$. Ceci combin\'e \`a \ref{std} montre l'implication $ii)\Rightarrow i) $ de \ref{defniv0}. Pour voir l'implication r\'eciproque, fixons une paire $(\PC,\o\PC)$ comme dans \ref{defniv0} ii).
Sous l'hypoth\`ese \ref{defniv0} i), on a $\ve \in R\u{G}^\dag e_{\u{U}^\dag} R\u{G}^\dag \cap
R\u{G}^\dag e_{\u{\o{U}}^\dag} R\u{G}^\dag$ donc $\ve\in R\u{G}^\dag e_{\u{U}^\dag} R\u{G}^\dag e_{\u{\o{U}}^\dag} R\u{G}^\dag$ puisque $\ve$ est idempotent. Or par le lemme ci-dessous appliqu\'e \`a $\u\GC'=\u\GC^\dag$, on a une d\'ecomposition d'Iwahori $\u{G}^\dag=\u{U}^\dag \u{M}^\dag \u{\o{U}}^\dag$, et donc $e_{\u{U}^\dag} R\u{G}^\dag e_{\u{\o{U}}^\dag} = R\u{M}^\dag e_{\u{U}^\dag} e_{\u{\o{U}}^\dag}$. Ceci cl\^ot la preuve de \ref{defniv0}.

\begin{lemme} \label{dilatIwa}
Soit $(\PC,\o\PC)$ une paire de paraboliques $\u\GC$-admissibles oppos\'es, et $\u\GC'\To{\varphi} \u\GC$ la dilatation dans $\u\GC$ d'un sous-groupe lisse connexe $\HC_k$ de $\u\GC_k$ de la forme $(\HC_k\cap \u\UC_k)(\HC_k\cap \u\MC_k)(\HC_k\cap \o{\u\UC}_k)$. Alors le triplet $(\u{U}',\u{M}',\u{\o{U}}')$ induit une d\'ecomposition d'Iwahori de  $\u{G}'=\u\GC'(\OC_K)$ au sens de \ref{defIwa}. 
\end{lemme}

%Notons que ce r\'esultat s'applique en particulier \`a $\GC'=\GC^\dag$. 

\begin{proof}
Il est sous-entendu dans l'\'enonc\'e que les notations $\u{U}'$, etc..., d\'esignent les groupes de points entiers des  adh\'erences sch\'ematiques $\u\UC'$, etc..., de $\UC$, etc... dans $\u\GC'$. Ces groupes sont bien ferm\'es dans $\u{G}'$ et satisfont les hypoth\`eses du pr\'eambule de \ref{defIwa}.

Montrons qu'ils satisfont le point i) de \ref{defIwa}.
On a d\'eja remarqu\'e que la restriction $\u\UC' \To{} \u\UC$ de $\varphi$ s'identifie \`a la dilatation de $\HC_k\cap \u\UC_k$ dans $\u\UC$. Reprenons alors les notations de \ref{BT} ii). L'ouvert r\'eciproque $\u\CC'$ de $\u\CC$ dans $\u\GC'$ s'identifie \`a la dilatation dans $\u\CC$ de $\HC_k$ (qui par hypoth\`ese est bien contenu dans $\u\CC_k$). Par commutation des dilatations aux produits ou par application de \cite[2.2.3]{BT2}, on a $\u\CC' = \u\UC'\times \u\MC'\times \u\UC'$ et l'immersion ouverte dans $\u\GC'$ est donn\'ee par le produit $\mu'$ dans $\u\GC'$.
Soit alors $d\in \OC_K[\u\GC]$ une \'equation du diviseur compl\'ementaire de $\u\CC$ dans $\u\GC$.
On a donc $\OC_K[\u\CC] =\OC_K[\u\GC][1/d]$ et  $\OC_K[\u\CC']=\OC_K[\u\GC'][1/d]$.
Comme le support ${\HC_k}$ de la dilatation effectu\'ee dans $\u\GC$ est disjoint du diviseur $d=0$, l'\'el\'ement $d$ est inversible dans le localis\'e $\OC_K[\u\GC']_{(\varpi)}$, et on a donc $\OC_K[\u\CC']_{(\varpi)}\simeq\OC_K[\u\GC']_{(\varpi)}$. En particulier $\mu'$ induit une bijection $\u\CC'(\OC_K) \To{} \u\GC'(\OC_K)$. D'o\`u le point i) de \ref{defIwa}.

Pour obtenir le point ii), on remarque que la discussion ci-dessus s'applique \`a la dilatation 
$\u\GC^1\To{} \u\GC$  de l'unit\'e de $\u\GC_k$ dans $\u\GC$, et par r\'ecurrence, aux dilatations successives $\u\GC^n\To{} \u\GC^{n-1}$  de l'unit\'e de la fibre sp\'eciale. Or les groupes $\u{G}^n=\u\GC^n(\OC_K)$ sont les groupes de congruences de $\u\GC$ ; ils sont  normaux et forment un syst\`eme de voisinages ouverts de l'unit\'e, {\em cf} \cite[2.8]{Yumodel}.
\end{proof}

%\begin{lemme}
Si $\PC$ est un \para $\u\GC$-admissible de $\GC$, le groupe $\wt\PC_k:={^u\u\GC_k}\u\PC_k$ est un sous-groupe parabolique de $\u\GC_k$. Plus pr\'ecis\'ement, supposons $\PC=\PC_\Omega$ comme dans \ref{std}, et introduisons le quotient r\'eductif $^q\u\GC_k:=\u\GC_k/{^u\u\GC_k}$ de $\u\GC_k$ et son syst\`eme de racines $\Phi_\dag:=\Phi(\u\SC_k,{^q\u\GC_k})\subseteq \Phi \subset V$. On v\'erifie alors que $\wt\PC_k/{^u\u\GC_k}$ est le parabolique de $^q\u\GC_k$ contenant $\u\SC_k$ et associ\'e au sous-ensemble parabolique $\Omega_\dag:=\Omega\cap \Phi_\dag$ de $\Phi_\dag$. L'application 
$$\{\hbox{Paraboliques $\u\GC$-admissibles de $\GC$}\} \To{} \{\hbox{Paraboliques de $\u\GC_k$}\}$$
ainsi obtenue est croissante pour la relation de contenance, surjective, mais g\'en\'eralement pas injective.
Nous appelons ``co-rang r\'esiduel" de $\PC$ dans $\u\GC$  le co-rang du parabolique associ\'e  ${\wt\PC_k}$ dans ${\u\GC_k}$. Il est toujours inf\'erieur ou \'egal au co-rang de $\PC$ dans $\GC$.

On d\'efinit maintenant $\u\GC_\PC \To{\nu_\PC} \u\GC$ la dilatation dans $\u\GC$ de $\wt\PC_k$. 
Le $\OC_K$-sch\'ema $\u\GC_\PC$ est un mod\`ele lisse connexe de $\GC$ qui a le m\^eme ensemble de paraboliques admissibles que $\u\GC$. Notons que $\PC$ a un corang r\'esiduel nul dans $\u\GC_\PC$ et que $\nu_\PC$ est un isomorphisme \ssi\ $\PC$ a un corang r\'esiduel nul dans $\u\GC$. De plus si $\QC \subset \PC$ est admissible, alors $(\u\GC_\PC)_\QC\simeq \u\GC_\QC$.
 % ({\em i.e.} $\wt{\PC}_{\PC,k} = \u\GC_{\PC,k}$).
Soit $\u\GC_\PC^\dag\To{} \u\GC_\PC$ la dilatation du radical unipotent de la fibre sp\'eciale de $\u\GC_\PC$. On prendra garde au fait que cette dilatation {\em ne se factorise pas} par $\u\GC^\dag$ mais au contraire on a une fl\`eche canonique $\u\GC^\dag\To{} \u\GC_\PC^\dag$ (penser aux groupes parahoriques : les relations de contenance des pro-$p$-radicaux sont oppos\'ees \`a celles des groupes eux-m\^emes). 
En fait, $\u\GC_\PC^\dag \To{} \u\GC$ s'identifie \`a la dilatation de ${^u{\wt\PC_k}}$ dans $\u\GC$. On a donc, par d\'efinition,  $\u{G}_\PC = \u{G}^\dag \u{P}$ et $\u{G}_\PC^\dag = \u{G}^\dag \u{U}$, et par \ref{dilatIwa}, le triplet
$(\u{\o{U}}^\dag,\u{M},\u{U})$, resp. $(\u{\o{U}}^\dag,\u{M}^\dag,\u{U})$, induit une d\'ecomposition d'Iwahori de $\u{G}_\PC$, resp. $\u{G}_\PC^\dag$.

Par ailleurs, si $\QC$ est un autre parabolique $\u\GC$-admissible, la projection $\u{G} \To{} {^q\u\GC_k}(k)$ induit une bijection des doubles classes :
 $$ \u{Q} \ba \u{G} / \u{G}_\PC \simeq \u{G}_\QC\ba \u{G} / \u{P} \simto \wt\QC_k\ba \u\GC_k / \wt\PC_k.$$
Plus pr\'ecis\'ement, supposons  $\QC$ et $\PC$ semi-standards
 %{\em i.e.} de la forme $\PC=\PC_\Omega$ et $\QC=\PC_\Theta$ comme dans \ref{std}, 
pour un choix de tore d\'eploy\'e maximal $\u\SC$ dans $\u\GC$. % et notons $\LC$ et $\MC$ leurs composantes de Levi semi-standard respectives. 
D'apr\`es \cite[Exp. XI, Cor. 5.3.bis]{SGA3}, le $\OC_K$-foncteur en groupes ``normalisateur de $\u\SC$ dans $\u\GC$" est repr\'esentable par un sch\'ema lisse sur $\OC_K$. 
On sait \cite[1.1.13]{BT2} que le groupe de Weyl $W(\Phi_\dag)$  s'identifie \`a $W(\u\SC_k,{\u\GC_k}):=\NC_{\u\GC}(\u\SC)(k)/\ZC_{\u\GC}(\u\SC)(k)$. Par lissit\'e, on peut donc relever chaque $w\in W(\u\SC_k,\u\GC_k)$ en une $\OC_K$-section $n_w$ de $\NC_{\u\GC}(\u\SC)$.
On obtient ainsi une d\'ecomposition
\ini\begin{equation} \label{decompmod}
 \u{G} = \bigsqcup_{w\in [W(\u\SC_k,\wt\QC_k)\ba W(\u\SC_k,\u\GC_k)/W(\u\SC_k,\wt\PC_k)]} \u{P_\Theta} n_w \u{P}_\Omega \u{G}^\dag 
 \end{equation}
o\`u les crochets d\'esignent un ensemble de repr\'esentants des doubles classes dans $W(\u\SC_k,\u\GC_k)$.

\alin{Preuve de \ref{theomod} : r\'eduction au co-rang r\'esiduel $1$}
Commen{\c c}ons par remarquer que l'\'enonc\'e de \ref{theomod} est vide lorsque $\PC$ est de corang r\'esiduel nul dans $\GC$, puisqu'on a alors $\u{U}=\u{U}^\dag$. Nous supposerons donc ce corang r\'esiduel strictement positif et nous allons
nous ramener au cas o\`u il vaut $1$. 
Pour cela,  il sera pratique de supposer $\PC$ semi-standard, {\em i.e.} de la forme $\PC_\Omega$ pour un choix de tore maximal d\'eploy\'e $\u\SC$ dans $\u\GC$ et un sous-ensemble parabolique $\Omega$ de $\Phi(\SC,\GC)$. 
Rappelons alors qu'on peut trouver une suite 
$\Omega_0:=-\Omega ,\Omega_1,\cdots, \Omega_r:=\Omega$ de sous-ensembles paraboliques de $\Phi$, de composante de Levi commune $\Omega^0$ et tels que
\begin{enumerate}
        \item pour tout $i$, la r\'eunion $\Theta_i:= \Omega_{i-1}\cup \Omega_i$ est un sous-ensemble parabolique tel que $\hbox{Vect}_\RM( \Omega^0)$ est de codimension $1$ dans $\hbox{Vect}_\RM( \Theta_i^0)$.
        %dans lequel $\Omega_i$ et $\Omega_{i-1}$ sont de co-rang $1$ (\em i.e.} .        
        \item la suite $\Omega^+\cap \Omega_i^+$ est strictement croissante.
\end{enumerate}
(Pour m\'emoire : les paraboliques dont la composante de Levi contient $\Omega^0$ sont de la forme $\Omega(v^*)$ pour $v^*\in \hbox{Vect}_\RM( \Omega^0)^\perp$. Les classes d'\'equivalence pour la relation $\Omega(v^*_1)=\Omega(v^*_2)$ sur $\hbox{Vect}_\RM(\Omega^0)^\perp$ sont les facettes d'une d\'ecomposition en c\^ones convexes. Les c\^ones ouverts (chambres) correspondent aux paraboliques de Levi $\Omega^0$. On obtient des $\Omega_i$ comme dans l'\'enonc\'e en choisissant une galerie tendue ({\em i.e.} de longueur minimale, {\em cf} \cite[1.1]{BT1}) entre les c\^ones correspondant \`a $\Omega$ et $-\Omega$.)

Posons alors $\u\GC_i:=\u\GC_{\PC_{\Theta_i}}$. Par construction, $\PC$ est de co-rang r\'esiduel $\leq 1$ dans chaque $\u\GC_i$. 

\begin{lemme}
Supposons que pour chaque $i$, l'\'enonc\'e de \ref{theomod} est v\'erifi\'e lorsqu'on remplace $\u\GC$ par  $\u\GC_{i}$. Alors l'\'enonc\'e de \ref{theomod} est v\'erifi\'e pour $\u\GC$.
\end{lemme}

\begin{proof}
Restant coh\'erents avec le syst\`eme de notations utilis\'e jusqu'ici, nous notons $\u{\UC}_i$, resp. $\u\UC_i^\dag$, l'adh\'erence sch\'ematique de $\UC$ dans $\u\GC_i$, resp. dans $\u\GC_i^\dag$.
Insistons sur le fait que le $\dag$ de $\u\GC_i^\dag$ se rapporte \`a $\u\GC_{i}$ et non \`a $\u\GC$, {\em cf} plus haut. En particulier $\u{U}_i^\dag$ est g\'en\'eralement distinct de $ \u{U}^\dag\cap \u{G}_i$.

%Par exemple on a $\u{G}_{\Theta_i}^\dag = \u{G}^\dag \u{U}_{\Theta_i^+}$.
Nous allons prouver
\ini\begin{equation} \label{interm+}
\u{\o{U}}_1=\u{\o{U}},\;\; \u{\o{U}}_r^\dag=\u{\o{U}}^\dag \; \hbox{ et } \;
 \forall i=1,\cdots, r, \;\; \u{\o{U}}_i^\dag=\u{\o{U}}_{i+1} \end{equation}
et sym\'etriquement
\ini\begin{equation} \label{interm-}
\u{{U}}_1^\dag=\u{{U}}^\dag,\;\; \u{{U}}_r= \u{{U}} \; \hbox{ et } \;
 \forall i=1,\cdots, r, \;\;  \u{{U}}_i=\u{{U}}_{i+1}^\dag 
\end{equation}
De \ref{interm+} nous retiendrons en particulier $e_{\u{\o{U}}}=e_{\u{\o{U}}_i}e_{\u{\o{U}}}$ pour tout $i=1,\cdots,r$, de sorte que l'hypoth\`ese du lemme nous donne
$$\forall i=1,\cdots, r, \;\; e_{\u{{U}}_i^\dag}e_{\u{\o{U}}} \in R\u{G}_i e_{\u{{U}}_i}e_{\u{\o{U}}}
.$$ On termine alors la preuve du lemme en \'ecrivant gr\^ace \`a \ref{interm-}
$$ e_{\u{U}^{\dag}}e_{\u{\o{U}}} =  e_{\u{U}_1^\dag}e_{\u{\o{U}}} \in R\u{G} e_{\u{U}_1}e_{\u{\o{U}}} =  R\u{G} e_{\u{U}_2^\dag}e_{\u{\o{U}}} \subset R\u{G}
e_{\u{U}_2}e_{\u{\o{U}}} \cdots \subset R\u{G} e_{\u{U}_r}e_{\u{\o{U}}}= R\u{G} e_{\u{U}}e_{\u{\o{U}}}
$$ o\`u les $\cdots$ d\'esignent une r\'ecurrence \'evidente.
Reste donc \`a prouver \ref{interm+} et \ref{interm-}. Pour des raisons de sym\'etrie,on se contentera de prouver \ref{interm-}.
Vue la d\'efinition de $\GC_r$, l'\'egalit\'e $\u{U}_r=\u{U}$ \'equivaut \`a la relation $\Omega^+\subseteq \Theta_r$, laquelle est vraie puisque $\Theta_r\supset \Omega$. Vue la d\'efinition de $\GC_1^\dag$, l'\'egalit\'e $\u{U}_1^\dag=\u{U}^\dag$ \'equivaut \`a la relation $\Omega^+\subseteq -\Theta_1$, laquelle est assur\'ee par $\Theta_1\supset \Omega_0=-\Omega$. 
Enfin pour $i=1,\cdots,r$, vues les d\'efinitions de $\GC_i$ et $\GC_{i+1}^\dag$, l'\'egalit\'e $\u{U}_i=\u{U}_{i+1}^\dag$ \'equivaut \`a l'identit\'e $ \Omega^+\cap \Theta_{i}= \Omega^+\cap \Theta_{i+1}^+  $. Or,
%\begin{eqnarray*}
$$\Omega^+\cap \Theta_{i+1}^+ = \Omega^+\cap (\Omega_i^+\cap \Omega_{i+1}^+) = \Omega^+\cap \Omega_i^+ $$
%\end{eqnarray*}
par la propri\'et\'e ii) de la suite $(\Omega_i)_i$, tandis que
$$\Omega^+\cap \Theta_{i} = \Omega^+\cap (\Omega_i\cup\Omega_{i-1})=\Omega^+\cap (\Omega_i^+\cup\Omega_{i-1}^+)= \Omega^+\cap \Omega_i^+, $$
par cette m\^eme propri\'et\'e ii) (la deuxi\`eme \'egalit\'e vient de $\Omega^+\cap \Omega^0=\emptyset$).
 
\end{proof}

\alin{Preuve de \ref{theomod} : r\'ecurrence}
Fixons un mod\`ele lisse connexe $\u\GC$ de $\GC$.
Nous allons d\'emontrer le pr\'edicat \`a trois variables suivant : {\em 
Pour tout couple $\PC\subset \QC$ de paraboliques admissibles et tout idempotent essentiellement de niveau z\'ero $\ve$ de $R\u{M}$, l'\'enonc\'e de \ref{theomod} est vrai pour $\PC$, $\ve$ et pour $\u\GC_\QC$ \`a la place de $\u\GC$}, par r\'ecurrence sur $\PC$, c'est \`a dire en le supposant vrai pour tout $\PC'\subsetneq \PC$. Il suffira alors de prendre $\QC=\GC$ pour obtenir \ref{theomod}. 

D'apr\`es \ref{std}, on peut supposer $\PC$ et $\o\PC$ semi-standards relativement \`a un sous-tore d\'eploy\'e maximal $\u\SC$ de $\u\GC$.
D'apr\`es le lemme pr\'ec\'edent, on peut aussi supposer $\PC$ de corang r\'esiduel $1$ dans $\u\GC_\QC$. Par ailleurs, $\u\GC_Q$ ayant les m\^eme propri\'et\'es formelles que $\u\GC$, nous pouvons supposer  $\u\GC_\QC=\u\GC$ pour prouver le pas de r\'ecurrence, et $\QC$ n'apparaitra plus dans ce qui suit.

Enfin, en l'absence d'ambigu\"it\'e, nous all\`egerons les notations en cessant de souligner les groupes de points entiers ; $\u{G}$ devient donc $G$, $\u{U}^\dag$ devient $U^\dag$, etc...

Comme dans \cite{HL}, le principe de la preuve repose sur l'\'egalit\'e  suivante dans $RG$ :
$$ [U : U^\dag] e_{U^\dag} e_{\o{U}} e_{U} e_{\o{U}} =  
 e_{U^\dag} e_{\o{U}} e_{U^\dag} e_{\o{U}} + \sum_{u\in U/U^\dag \setminus \{1\}} e_{U^\dag} e_{\o{U}} u e_{U^\dag} e_{\o{U}}. $$

D'apr\`es  \ref{propIwa} appliqu\'e \`a $G_{\o\PC} = U^\dag M^\dag \o{U}$, il existe un \'el\'ement central inversible $z_\PC$ dans $RM$ tel que $e_{U^\dag} e_{\o{U}} e_{U^\dag} e_{\o{U}} = z_\PC e_{U^\dag} e_{\o{U}}$. Ainsi pour prouver l'\'enonc\'e de \ref{theomod}, il suffira de prouver
que pour tout $u\in U\setminus U^\dag$, on a
\ini\begin{equation} \label{xu}
 X(u):= \varepsilon e_{U^\dag} e_{\o{U}} u e_{U^\dag} e_{\o{U}} \varepsilon \in RG e_{U}e_{\o{U}}\ve .
 \end{equation}
D'apr\`es \ref{decompmod} appliqu\'ee \`a $u$, il existe un \'el\'ement $n_w\in
\NC_{\u\GC}(\u\SC)(\OC_K)$ d'image $w\in W(\SC_k,{^q\u\GC_k})\setminus W(\SC_k,{^q\u\MC_k})$, des \'el\'ements  $\o{u}_i\in \o{U}$ et $ m_i\in M$ pour $i=1,2$  et un \'el\'ement $g^\dag\in G^\dag$ tels que $u= m_1 \o{u}_1 n_w \o{u}_2 m_2 g^\dag$.
%$$ u e_{U^\dag} = m_1 u_1 n_w RG^\dag. RU_{\Omega^+}. RM_\Omega.$$
 Puisque $M$ normalise $\o{U}$ et $U^\dag$,
il s'ensuit que  
$$ X(u) \in  RM. \ve e_{U^\dag} e_{\o{U}} n_w .RG^\dag. e_{\o{U}} \ve RM .$$
%$$ e_{U} e_{\o{U}} e_{U} e_{\Omega} \varepsilon \in \sum_{w} \varepsilon RM_\Omega e_{\o{U}} $$
 Utilisons maintenant l'hypoth\`ese faite sur $\ve$ d'\^etre ``essentiellement de niveau $0$" ; elle nous permet d'\'ecrire que $\ve \in RM^\dag e_{M\cap {^w}U^\dag}RM^\dag$, o\`u on a pos\'e
 ${^w}?:= w ? w^{-1}$.
 %pos\'e ${^w}\Omega^-:=w\Omega^- w^{-1}$ et o\`u par convention $e_\emptyset^\dag=1$.
  On en d\'eduit
$$ X(u) \in  RM. \ve e_{U^\dag} e_{M\cap {^w{U^\dag}}} e_{\o{U}} n_w .RG^\dag. e_{\o{U}} \ve RM. $$
On d\'ecompose maintenant le premier $e_{\o{U}}$ en un produit $e_{\o{U}\cap {^w}U}e_{\o{U}\cap {^wM}} e_{\o{U}\cap {^w}\o{U}}$ gr\^ace au point iv) de \ref{BT}. En remarquant que $e_{\o{U}\cap {^w}\o{U}} n_w. RG^\dag. e_{\o{U}} \subseteq n_w.RG^\dag.e_{\o{U}}$, on obtient
$$ X(u) \in  RM. \ve e_{U^\dag} e_{M\cap {^w}U^\dag} e_{\o{U}\cap {^w}U} e_{\o{U} \cap {^w}M} n_w .RG^\dag. e_{\o{U}} \ve RM. $$
Observons que le produit 
$e_{U\cap {^w}U^\dag} e_{M\cap {^w}U^\dag} e_{\o{U}\cap {^w}U}$ est l'idempotent associ\'e au pro-$p$-groupe ${^w}U^\dag ({^w}U\cap \o{U})$ et s'\'ecrit encore 
%Or, $e_{U^\dag} e^\dag_{\Omega^0\cap {^w}\Omega^-} e_{\Omega^+\cap {^w}\Omega^-} = 
$e_{^wU^\dag} e_{\o{U}\cap {^w}U} =  e_{\o{U}\cap {^w}U} e_{^wU^\dag} $. Ainsi, en d\'ecomposant $e_{U^\dag}=e_{U^\dag\cap {^w}\o{U}} e_{U^\dag\cap {^w}M} e_{U^\dag\cap {^w}U}$, puis en commutant \`a $n_w$, on obtient en posant $?^w:=w^{-1}?w= {^{w^{-1}}}?$
\begin{eqnarray*} X(u) & \in &  \left(RM. \ve e_{U^\dag\cap{^w}\o{U}}\,   n_w\right)e_{M\cap U^{\dag w}}(e_{U\cap\o{U}^w} \,e_{U^\dag}) e_{M\cap\o{U}^w} .RG^\dag. e_{\o{U}} \ve RM \\
& \subset &  RG. e_{M\cap U^{\dag w}} e_{U\cap\o{U}^w} e_{M\cap\o{U}^w} e_{U^\dag} .RG^\dag. e_{\o{U}} \ve RM \\
& = &  RG. e_{M\cap U^{\dag w}} e_{U\cap\o{U}^w} e_{M\cap\o{U}^w} e_{U^\dag} e_{\o{U}} \ve RM \\
& = &  RG. e_{M\cap U^{\dag w}} e_{U\cap\o{U}^w}  e_{U^\dag}  e_{\o{U}} e_{M\cap\o{U}^w} . \ve RM \\
& = &  RG.  e_{U\cap\o{U}^w}  (e_{U^\dag} e_{M\cap U^{\dag w}})  (e_{\o{U}} e_{M\cap\o{U}^w}) . \ve RM  \;\;\;\;\;(*)
\end{eqnarray*}
Dans la deuxi\`eme ligne, on fait commuter $ e_{M\cap\o{U}^w}$ et $ e_{U^\dag} $ car $M$ normalise $U^\dag$.
\`A la troisi\`eme ligne, on utilise la d\'ecomposition \`a la Iwahori $G^\dag = U^\dag M^\dag \o{U}^\dag$ pour \'ecrire  $e_{U^\dag} .RG^\dag. e_{\o{U}} =  e_{U^\dag}  e_{\o{U}}. RM^\dag$.
La quatri\`eme ligne vient encore du fait que $M$ normalise $U^\dag$ et $\o{U}$, et la derni\`ere vient du fait que $M^\dag$ normalise le groupe $U^\dag  (U\cap \o{U}^w)$ puisqu'il normalise $U$, $U^\dag$ et agit trivialement sur le quotient.

Deux cas se pr\'esentent maintenant : supposons tout d'abord que 
$\MC\cap {\PC}^w$ soit un parabolique propre de $\MC$ et posons $\PC':=(\MC\cap \PC^w)\UC$. C'est un \para de $\GC$ qui contient strictement $\PC$. Son radical unipotent est $(\MC\cap \UC^w)\UC$, sa composante de Levi semi-standard est $\LC:=\MC\cap \MC^w$, et son oppos\'e semi-standard est $\o\PC':=(\MC\cap \o\PC^w)\o\UC$. Par ailleurs, si $\MC=\ZC_\GC(\TC)$ pour un sous-tore de $\SC$, on a $\LC=\ZC_\GC(\TC.\TC^w)$ et $\PC'$ est donc $\u\GC$-admissible.
 La  d\'ecomposition d'Iwahori $M^\dag = (M^\dag\cap U^w) L^\dag (M^\dag\cap\o{U}^w)$ 
nous fournit une unique distribution $\ve_L$ dans $RL^\dag$ telle que $e_{M^\dag\cap U^w} \ve e_{M^\dag\cap\o{U}^w}  = e_{M^\dag\cap U^w} \ve_L e_{M^\dag\cap \o{U}^w}$.
Par unicit\'e et par la proposition \ref{propIwa}, $\ve_L$ est un idempotent central de $RL^\dag$. Pour v\'erifier que $\ve_L$ est essentiellement de niveau $0$, fixons un parabolique $\u\GC$-admissible minimal $\RC$ de $\GC$ et contenant $\PC'$. Par la d\'efinition \ref{defniv0}, on a
$ \ve \in RM^\dag e_{M^\dag\cap \o{V}} e_{M^\dag\cap{V}} RM^\dag $. On en d\'eduit
\begin{eqnarray*}
 e_{M^\dag\cap U^w} \ve e_{M^\dag\cap \o{U}^w}  & \in &  
e_{M^\dag\cap U^w}  RM^\dag e_{M^\dag\cap \o{V}} e_{M^\dag\cap{V}} RM^\dag  e_{M^\dag\cap \o{U}^w}  \\
& = &
e_{M^\dag\cap U^w}  RL^\dag e_{M^\dag\cap \o{U}^w}  e_{L^\dag\cap \o{V}} e_{L^\dag\cap {V}}
e_{M^\dag\cap U^w}  RL^\dag  e_{M^\dag\cap \o{U}^w}  \\
 & = & 
e_{M^\dag\cap U^w}  RL^\dag  e_{L^\dag\cap \o{V}} e_{L^\dag\cap {V}}
  RL^\dag  e_{M^\dag\cap \o{U}^w} 
\end{eqnarray*}
On a utilis\'e la d\'ecomposition $M^\dag=(M^\dag\cap U^w)L^\dag(M^\dag\cap \o{U}^w)$ pour changer $M$ en $L$ dans la seconde ligne, et on a appliqu\'e \ref{propIwa} \`a cette m\^eme d\'ecomposition pour passer \`a la troisi\`eme ligne.
%On a utilis\'e l'\'egalit\'e $e_{M^\dag\cap U^w} e_{M^\dag\cap U^w} e_{M^\dag\cap U^w} e_{M^\dag\cap U^w} =e_{M^\dag\cap U^w} e_{M^\dag\cap U^w} $
Par unicit\'e, il s'ensuit que  $\ve_L$ est essentiellement de niveau $0$ dans $RL$.
On peut donc  appliquer l'hypoth\`ese de r\'ecurrence \`a $\PC'$ et $\ve_L$. 
On obtient, \`a partir de la ligne $(*)$
\begin{eqnarray*}
 X(u) & \in  & RG. (e_{U} e_{M\cap U^w})  (e_{\o{U}} e_{M\cap\u{U}^w}) . \ve RM \\
    & \subset & RG  e_{U}e_{\o{U}} \ve
 \end{eqnarray*}
toujours en utilisant le fait que $M$ normalise $U$ et $\o{U}$.

Supposons au contraire que $\MC\cap \PC^w=\MC$, ou de mani\`ere \'equivalente, que $\MC=w^{-1}\MC w$. C'est ici qu'intervient l'hypoth\`ese de corang r\'esiduel $1$. En effet,
$w^{-1}$ normalise aussi $^q\u\MC_k$  qui sous cette hypoth\`ese est un Levi maximal (et propre) de $^q\u\GC_k$. Il s'ensuit que le groupe de Weyl $W(\u\SC_k,\u\MC_k)$ est normal et d'indice $2$ dans le groupe de Weyl $W(\u\SC_k,\u\GC_k)$, et puisque $w$ n'est pas dans  $W(\u\SC_k,\u\MC_k)$, son action par conjugaison \'echange le radical unipotent de $\wt\PC_k$ et celui de $\wt{\o\PC}_k$. Par cons\'equent, la conjugaison par $n_w$ induit un isomorphisme $\GC_\PC \simto \GC_{\o\PC}$, induisant \`a son tour un isomorphisme  $\GC_\PC^\dag \simto \GC_{\o\PC}^\dag$ qui sur les point entiers signifie simplement que $w^{-1} G^\dag \o{U} w = G^\dag U$. On en tire imm\'ediatement que $U = U^\dag(U\cap \o{U}^w)$, et  la ligne $(*)$ ci-dessus nous donne
$$X(u) \in  RG. e_{U}e_{\o{U}} \ve.$$
  
On a donc termin\'e la preuve de \ref{xu} et par l\`a celle du th\'eor\`eme \ref{theomod}.

\alin{Preuve du corollaire \ref{utile}}
Comme dans la preuve pr\'ec\'edente, nous all\'egeons les notations en ne soulignant pas les groupes de points entiers.
D'apr\`es la proposition \ref{propIwa} appliqu\'ee \`a la d\'ecomposition d'Iwahori ${G}_{\o\PC}^\dag={U}^\dag {M}^\dag {\o{U}}$, on a 
${e_{U^\dag}}{e_{\o{U}}}\in RG {e_{\o{U}}}{e_{U^\dag}}{e_{\o{U}}}$. Ainsi, gr\^ace \`a l'hypoth\`ese i), %${U'}_{{\Omega'}^-}^\dag \subseteq U_{\Omega^-}^\dag$,
$$ {e_{U^\dag}}{e_{\o{U}}}\ve' \in RG{e_{\o{U}}}{e_{U^\dag}}{e_{\o{U}}}\ve' = RG {e_{\o{U}}} \wt\ve' {e_{U^\dag}} \wt\ve' {e_{\o{U}}} .$$
D'apr\`es notre hypoth\`ese sur l'entrelacement de $\wt\ve'$ et notre hypoth\`ese $U^\dag\cap U'\subseteq {U'}^\dag$, on a donc
$${e_{U^\dag}}{e_{\o{U}}}\ve' \in RG \wt\ve' e_{{U'}^\dag} \wt\ve' {e_{\o{U}}}
= RG \wt\ve' e_{{U'}^\dag} e_{{\o{U}}'} \ve'{e_{\o{U}}}.$$ D'apr\`es le th\'eor\`eme \ref{theomod} appliqu\'e \`a $\u\GC'$ et $\ve'$, on obtient
$${e_{U^\dag}}{e_{\o{U}}}\ve' \in RG \wt\ve' \, RG' \, e_{{U'}} e_{{\o{U}}'}  \ve' {e_{\o{U}}} = RG \wt\ve' \, RG' \,e_{{U'}} \wt\ve' {e_{\o{U}}} = RG \wt\ve' e_{{U'}} \wt\ve' {e_{\o{U}}} .$$
Or, par notre hypoth\`ese sur l'entrelacement, on a
%$\wt\ve' \, RG' \, e_{{U'}} u \, \wt\ve' =0$ pour tout $u\in U\setminus U'$, de sorte que pour toute $f\in RG' $, on a 
%$$ \wt\ve'. f.e_{{U'}} \wt\ve' = [U : U']^{-1} \wt\ve'. f .e_{U} \wt\ve' .$$
$\wt\ve'  e_{{U'}} u \, \wt\ve' =0$ pour tout $u\in U\setminus U'$, et donc 
$ \wt\ve' e_{{U'}} \wt\ve' = [U : U']^{-1} \wt\ve' e_{U} \wt\ve' .$
On en d\'eduit que 
$$ {e_{U^\dag}}{e_{\o{U}}}\ve' \in RG e_U \wt\ve' {e_{\o{U}}} = RG e_{U} {e_{\o{U}}} \ve'.$$

\alin{Preuve de \ref{constess}}
Comme dans les preuves pr\'ec\'edentes, nous all\'egeons les notations en ne soulignant pas les groupes de points entiers.
Pour un caract\`ere $\chi: U^\dag/U^* \To{} R^\times$, on note $[\chi]$ l'idempotent de $RU^\dag$ associ\'e. L'hypoth\`ese ii) implique que $R$ est suffisamment gros pour que $\sum_\chi [\chi]= e_{U^*}$, la somme \'etant sur tous les caract\`eres $\chi$. On a donc l'\'egalit\'e dans $RG^\dag$
$$ [\theta] = \sum_\chi [\chi][\theta] .$$
Fixons maintenant un tel $\chi$ ainsi qu'un \'el\'ement $u\in \o{U}^\dag$ et calculons l'expression $[\chi]u[\chi][\theta]$. On a
\begin{eqnarray*}
|U^\dag / U^*|.  [\chi]u[\chi][\theta] & = & \sum_{v\in U^\dag/U^*} \chi^{-1}(v)e_{U^*} v u [\chi][\theta] \\
   & = & \sum_{v} \chi^{-1}(v) e_{U^*} v u v^{-1} \chi(v) [\chi][\theta] \\
   &  = & e_{U^*} u  \sum_v [u^{-1},v] [\chi][\theta] = e_{U^*} u \sum_v \theta([u^{-1},v]) [\chi][\theta]  
\end{eqnarray*}
D'apr\`es le ii) de notre hypoth\`ese, la somme $\sum_v \theta([u^{-1},v])$ est non-nulle seulement si $u\in \o{U}^*$. On en d\'eduit que
\begin{eqnarray*}
|\o{U}^\dag / \o{U}^*| [\chi]e_{\o{U}^\dag}[\chi][\theta] & = & \sum_{u\in \o{U}^\dag / \o{U}^*} [\chi]e_{\o{U}^*} u [\chi][\theta] \\
& = & [\chi]e_{\o{U}^*}[\chi][\theta] = [\chi][\theta]
\end{eqnarray*}
et donc que
\begin{eqnarray*}
 [\theta] & = & \sum_\chi [\chi][\theta] = |\o{U}^\dag / \o{U}^*| \sum_\chi [\chi] e_{\o{U}^\dag} [\chi][\theta] \\
  & \in & RG^\dag e_{\o{U}^\dag} RG^\dag
\end{eqnarray*}
De m\^eme on prouve que $[\theta]\in RG^\dag e_{{U}^\dag} RG^\dag$ et par produit et d\'ecomposition d'Iwahori, on en d\'eduit que $[\theta]$ v\'erifie le point ii) de la d\'efinition \ref{defniv0}.

\alin{Preuve de \ref{lemmod}} La premi\`ere \'etape consiste \`a d\'evisser au cas o\`u $\varphi$ est une dilatation. Notons pour cela que, en tant que quotient de $\u\GC'_k$, le groupe image $\im(\varphi_k)$ est r\'eduit, donc  lisse. Soit $\u\GC^1\To{\nu^1} \u\GC$ la dilatation de $\im(\varphi_k)$. On sait que $\u\GC^1$ est lisse et, par la propri\'et\'e universelle des dilatations sur les $\OC_K$-sch\'emas plats, que $\varphi$ se factorise en $\nu^1\circ \varphi^1$ o\`u $\varphi^1 :\u\GC'\To{} \u\GC^1$ induit un isomorphisme des fibres g\'en\'eriques. It\'erons le proc\'ed\'e : on obtient une suite de dilatations $\u\GC^n\To{\nu^n} \u\GC^{n-1}$ et une suite de factorisations $\varphi^{n-1}=\nu^n\circ \varphi^n$. Notons que si $\nu^n$ n'est pas un isomorphisme, alors l'indice de $\u{G}^n$ dans $\u{G}^{n-1}$ est un entier $>1$. Comme l'indice de $\u{G}'$ dans $\u{G}$ est fini, il existe  un entier $n$ \`a partir duquel la suite devient stationnaire. On a alors $\u{G}'=\u{G^n}$, et par unicit\'e des mod\`eles lisses \cite[1.7]{BT2}, $\varphi^n$ est un isomorphisme. On a donc pr\'esent\'e $\varphi$ comme une composition de dilatations de $\u\GC$ \`a centres lisses. En particulier $\ker(\varphi_k)$ a une suite de composition dont les sous-quotients sont les $\ker(\nu^n)$.

On est donc ramen\'e au cas o\`u $\varphi_k$ est une dilatation (la dilatation de centre $\im(\varphi_k)$). Dans ce cas le noyau est m\^eme vectoriel  ; en effet, revenant \`a la d\'efinition originale d'une dilatation comme ouvert de l'\'eclatement du centre,
on peut identifier $\u\GC'_k$ comme l'extension vectorielle de $\im(\varphi_k)$ associ\'ee au fibr\'e conormal de $\im(\varphi_k)$ dans $\u\GC_k$.

\section{Paraboliques minimaux, niveau z\'ero}
Nous reprenons maintenant les notations du paragraphe \ref{notaimmeuble}. En particulier si $x\in B(\GC,K)$, on note $\GC_x$ le mod\`ele lisse de $\GC$ associ\'e \`a $x$ par Bruhat-Tits, et $\GC_x^\circ$ sa composante neutre. Commen{\c c}ons par la remarque suivante :
\begin{rema} \label{remapara}
Soit $x\in B(\GC,K)$ et $\MC$ un Levi de $\GC$. Les propri\'et\'es suivantes sont \'equivalentes :
\begin{enumerate}
        \item $x\in B(\MC,K)$
        \item L'adh\'erence sch\'ematique dans $\GC_x$ du tore d\'eploy\'e maximal du centre de $\MC$ est un tore.
        \item $\MC$ est $\GC_x^\circ$-admissible au sens de \ref{defmod}.
\end{enumerate}
\end{rema}
\begin{proof}
Supposons $x\in B(\MC,K)$, et soit $\SC$ un tore d\'eploy\'e maximal de $\MC$ dont l'appartement associ\'e contient $x$. Par construction, le groupe de Bruhat-Tits $\GC_x$ contient le prolongement canonique de $\SC$ \`a $\OC_K$. Comme $\SC$ contient la partie d\'eploy\'ee du centre de $\MC$, on en d\'eduit $i)\Rightarrow ii)$. 
L'implication $ii)\Rightarrow iii)$ est tautologique, vue la d\'efinition d'admissibilit\'e.
Prouvons donc $iii)\Rightarrow i)$.  Supposons que $\MC$ est $\GC_x^\circ$-admissible et choisissons un tore d\'eploy\'e $\TC$ de $\MC$  se prolongeant en un sous-tore $\TC_x$ de $\GC_x$ et dont le centralisateur dans $\GC$ est $\MC$. 
Choisissons alors  un tore d\'eploy\'e maximal $\SC_x$ de $\GC_x$ contenant $\TC_x$. (Pour ce faire, on choisit d'abord $S_{x,k}$ dans $\GC_{x,k}$  contenant $\TC_{x,k}$, on rel\`eve $\SC_{x,k}$ en un tore $\SC'_x$ par \cite[Exp IX. Thm 3.6]{SGA3}, puis $\TC_{x,k}$ en $\TC'_x\subset \SC'_x$, puis on utilise l'unicit\'e \`a conjugaison pr\`es des rel\`evements de $\TC_{x,k}$ pour conjuguer $\SC'_x$ en un tore contenant $\TC_x$).
 Alors la fibre g\'en\'erique $\SC$ de $\SC_x$ est un tore d\'eploy\'e maximal de $\GC$ contenu dans le centralisateur de $\TC$, c'est \`a dire $\MC$. Par d\'efinition $x$ appartient \`a l'appartement associ\'e \`a $\SC$ et par cons\'equent \`a $B(\MC,K)$. 
\end{proof}

\begin{prop} \label{pfcommutmin}
 \label{pfprinc} Soit $\PC=\MC\UC$ un \para  de $\GC$  et $x$ un point de l'immeuble $B(\MC,K)$. %Soit $\ve$ un idempotent central de $RM_x$ ``essentiellement de niveau $0$" au sens du lemme \ref{lemDRM2}. Alors
 \begin{enumerate}
        \item (niveau z\'ero)  $\euxp\eubx e_{M_x^+} \in (R G_x)\eux\eubx e_{M_x^+}$.
         \item Si $\PC$ est minimal, alors  $\euxp\eubx \in (R G_x)\eux\eubx$. 
        \end{enumerate}
 \end{prop}
\begin{proof}
Remarquons pour commencer que les groupes $\UC_x$ et $\o{\UC}_x$ obtenus par adh\'erence sch\'ema\-tique de $\UC$ dans $\GC_x$ sont {\em connexes} donc contenus dans $\GC_x^\circ$.
Pour le voir, on peut se ramener aux adh\'erences de groupes radiciels $\UC_{\alpha,x}$, lesquelles sont explicit\'ees par Bruhat-Tits (construction en \cite[4.3 et 5.2.2]{BT2} et propri\'et\'e d'immersion ferm\'ee de \cite[3.8.1 (S2)]{BT2}) qui montrent que leurs sch\'emas
sous-jacents sont des vectoriels sur $\OC_K$. Comme $G_x^+=(\GC_x^\circ)^\dag(\OC_K)$, on voit que tous les idempotents de l'\'enonc\'e sont  en r\'ealit\'e dans $G_x^\circ$.

Dans le point i), l'idempotent $e_{M_x^+}$ est clairement ``essentiellement de niveau z\'ero" au sens de  \ref{defniv0}, donc on peut appliquer le th\'eor\`eme \ref{theomod}. Dans le point ii), on remarque que l'idempotent $1_{M_x^+}$ (\'el\'ement unit\'e) satisfait aussi  les hypoth\`eses de \ref{defniv0} puisque $\MC$ n'a pas de parabolique $\MC_x^\circ$-admissible propre.
\end{proof}

En combinant ce corollaire avec \ref{commut}, on obtient bien la propri\'et\'e de commutation \ref{commutmin} annonc\'ee dans l'introduction pour les paraboliques minimaux, et par \ref{cororescent}, la propri\'et\'e de seconde adjonction dans ce cas.

En ce qui concerne le niveau z\'ero, voici le r\'esultat obtenu :

\begin{prop} \label{niveau0}
   Soit $\Mo{R}{G}_0$ la sous-cat{\'e}gorie pleine des objets engendr\'es par la r\'eunion de leurs $G_x^+$-invariants, pour $x\in B(\GC,K)$ (objets de ``niveau z\'ero"), 
 et $\Mo{R}{G}_{>0}$ celle des objets dont tous les  $G_x^+$-invariants sont nuls, pour $x\in B(\GC,K)$ (objets de ``niveau positif").%$W^{G_x^+}=0$, $x\in A$.
  \begin{enumerate}
  \item On a une d\'ecomposition $\Mo{R}{G}=\Mo{R}{G}_0 \oplus \Mo{R}{G}_{>0}$.
  \item Les  foncteurs paraboliques envoient objets de niveau z\'ero, resp. positif, sur objets de niveau z\'ero, resp. positif. 
    \item Pour tout \para $\PC$ de $\GC$, la restriction du foncteur $\Rp{G,P}{M}$
{\`a} la cat{\'e}gorie $\Mo{R}{G}_0$ est adjointe {\`a} droite de la restriction du foncteur
$\delta_P\Ip{M,\o{P}}{G}$ {\`a} la cat{\'e}gorie $\Mo{R}{M}_0$.
\item La cat\'egorie $\Mo{R}{G}_0$ est noeth\'erienne.
\end{enumerate}
\end{prop}
\begin{proof} La d\'ecomposition du point i) est expliqu\'ee dans l'appendice. 
Pour le point ii), soit $\PC=\MC\UC$ un sous-groupe parabolique. 
D'apr\`es le corollaire pr\'ec\'edent et \ref{commut}, on a $\ip{M,P}{G}{\cind{M_x^+}{M}{R}}\simeq \cind{G_x^+}{G_x}{R}$ pour tout $x\in B(\MC,K)$. Par exactitude de $\Ip{M,P}{G}$, on en d\'eduit que celui-ci envoie $\Mo{R}{M}_0$ dans $\Mo{R}{G}_0$. Par fid\'elit\'e, et puisque toute $G$-orbite dans $B(\GC,K)$ rencontre $B(\MC,K)$, on en d\'eduit qu'il envoie aussi $\Mo{R}{M}_{>0}$ dans $\Mo{R}{G}_{>0}$.
Par r\'eciprocit\'e de Frobenius, on en d\'eduit les propri\'et\'es analogues pour $\Rp{G,P}{M}$.

Pour le point iii), on recopie la preuve de \ref{cororescent} en utilisant le fait que les $e_{M_x^+}$ pour $x\in B(\MC,K)$ forment une famille g\'en\'eratrice de $\Mo{R}{M}_0$ d'idempotents $P$-bons.
Enfin, les arguments de  la partie \ref{noether} montrent que iii) implique iv).
\end{proof}

\section{GL(N)} \label{lineaire}

Pour coller aux notations de Bushnell, Kutzko et Stevens, nous noterons $F$ le corps local que nous notions $K$ jusqu'ici.

\alin{Dictionnaire Bruhat-Tits/Bushnell-Kutzko}
Ce dictionnaire est tr\`es bien expliqu\'e dans \cite{BrLem}
auquel on renvoie  le lecteur pour les d\'etails.
Soit $V$ un $F$-espace vectoriel.
Une fonction r\'eseau sur $V$ relativement \`a $F$ est une fonction 
${\Lambda:\;}{\RM}\To{} \{\OC_F-\hbox{r\'eseaux de }V\}$ qui est
d\'ecroissante, continue \`a gauche, et telle que $\Lambda(r+v_F)=\PG_F
\Lambda(r)$ pour tout $r\in\RM$, o\`u $v_F\in \RM_+$ d\'esigne la valuation d'une uniformisante de $F$. % que nous supposerons dor\'enavant \'egale \`a $1$.
L'ensemble $\FC\RC_F(V)$ de ces fonctions est muni
d'une action de $G$ et d'une action de $\RM$ par translations que nous noterons $\Lambda\mapsto \Lambda[t] :r \mapsto \Lambda(r-t)$. 
Soit $\GC$ le $F$-sch\'ema en groupes des automorphismes $F$-lin\'eaires de $V$, dont le groupe des points rationnels est $G=\aut{F}{V}$. D'apr\`es \cite[Prop. 2.4]{BrLem},
il y a une application naturelle $ B(\GC,F)\To{}\FC\RC_F(V)$. Celle-ci est bijective, 
$G$-\'equivariante, et $\RM$-\'equivariante si on identifie $X_*(\ZC(\GC))\otimes \RM\simto \RM$ en envoyant un endomorphisme scalaire de $V$ sur l'oppos\'e de la valuation de ce scalaire.

\def\aG{\mathfrak a}
\def\uG{\mathfrak u}
Notons $A:=\endo{F}{V}$. La fonction r\'eseau  $\Lambda\in \FC\RC_F(V)$
d\'etermine une fonction r\'eseau  $\aG(\Lambda)\in \FC\RC_F(A)$ d\'efinie par
$\aG_r(\Lambda)=\{x\in A,\,\forall u\in\RM,\, x\Lambda(u) \subseteq \Lambda(r+u)\}$ pour $r\in \RM$. On
note aussi $\aG_{r+}(\Lambda):= \bigcup_{s>r} \aG_s(\Lambda)$. Alors
$\aG_0(\Lambda)$ est un ordre h\'er\'editaire, $\aG_{0+}(\Lambda)$ est son
radical de Jacobson et tous les $\aG_r(\Lambda)$ en sont des id\'eaux
fractionnaires.
Posons aussi $\uG_0(\Lambda):=\aG_0(\Lambda)^\times$, c'est un
sous-groupe compact ouvert de $GL(V)$ dont la famille
$\uG_r(\Lambda):=1+\aG_r(\Lambda)$ pour $r>0$ est une filtration par
des pro-$p$-sous-groupes ouverts normaux.
Si $x\in B(\GC,F)$ est le point correspondant \`a 
$\Lambda$, on a $\uG_0(\Lambda)=G_x$ et $\uG_{0+}(\Lambda)=G_x^+$, et d'apr\`es \cite[Appendice A]{BrLem}, la filtration $G_{x,r}:=\uG_r(\Lambda)$, $r\geq 0$ est celle de Moy et Prasad. 

Si $\MC\subset \GC$ est le sous-groupe de Levi correspondant \`a une
d\'ecomposition $V=\bigoplus_{i\in I}V_i$, le ``sous-immeuble''
$B(\MC,F)$ de $B(\GC,F)$ s'identifie au sous-ensemble des fonctions
r\'eseaux d\'ecompos\'ees par $M$ au sens o\`u $\forall r\in \RM,
\Lambda(r)=\bigoplus_{i\in I} \Lambda(r)\cap V_i$. Soit $x$ le point de $B(\MC,F)$ associ\'e \`a une telle fonction r\'eseau, alors l'ensemble  $x+a_M$ est l'ensemble des fonctions r\'eseau de la forme $\bigoplus \Lambda_i[t_i]$ o\`u les $t_i$ sont des r\'eels.

Nous dirons que $\Lambda$ est rationnelle si ses sauts sont dans
$\QM v_F\subset\RM$. Il existe alors un plus petit entier positif
$e=e(\Lambda)$ tel que la fonction $\wt\Lambda: r\mapsto \Lambda(er/v_F)$
soit une suite de  r\'eseaux  au sens de
Bushnell-Kutzko \cite[2.1]{BK2} La p\'eriode de $\wt\Lambda$ est justement
$e(\Lambda)$. R\'eciproquement, une suite de r\'eseaux d\'etermine une
unique fonction r\'eseau ; il suffit de rendre la p\'eriode \'egale
\`a $v_F$.

\alin{Strates et caract\`eres semi-simples}
Soit $V$ un $F$-vectoriel et $A=\endo{F}{V}$. Une strate dans $V$ est un quadruplet $[\Lambda,n,r,\gamma]$ o\`u $\Lambda$ est une fonction r\'eseau, $\gamma\in \aG_{-n}(\Lambda)$ et $r<n$. On d\'efinit l'\'equivalence de telles strates comme dans le cas des suites de r\'eseaux consid\'er\'e par Bushnell et Kutzko. D'ailleurs, 
lorsque $\Lambda$ est rationnelle et $n\in \frac{v_F}{e(\Lambda)}\NM$, seul cas que l'on consid\`erera ici, le quadruplet $[\wt\Lambda,\frac{e(\Lambda)}{v_F}n,\frac{e(\Lambda)}{v_F}r,\gamma]$ est une honn\^ete strate au sens de \cite[3.1]{BK}. On dira alors que $[\Lambda,n,r,\gamma]$ est fondamentale, simple, semi-simple si $[\wt\Lambda,\frac{e(\Lambda)}{v_F}n,\frac{e(\Lambda)}{v_F}r,\gamma]$ l'est.

\`A toute strate  $[\Lambda,n,0,\beta]$ simple, resp. semi-simple,  Bushnell et Kutzko \cite[5]{BK2}, resp Stevens \cite[3]{Stevens}, associent  deux sous-ordres $\jG(\Lambda,\beta)\supseteq \hG(\Lambda,\beta)$ de $\aG_0(\Lambda)$, et un ensemble de caract\`eres dits {\em simples}, resp. {\em semi-simples}, du sous-groupe $H^+(\Lambda,\beta):=\hG(\Lambda,\beta)\cap \uG_{0+}(\Lambda)$ de $J(\Lambda,\beta):=\jG(\Lambda,\beta)\cap \uG_0(\Lambda)$.

Comme le groupe $H^+(\Lambda,\beta)$ est pro-$p$, les caract\`eres (semi-)simples sont \`a valeurs dans l'anneau $\ZM_{p-cycl}$ des entiers de l'extension $p^{\infty}$-cyclotomique de $\QM$. Rappelons aussi que leur d\'efinition d\'epend du choix d'un caract\`ere $\psi:F/\PG_F\To{} \ZM_{p-cycl}^\times$.

\begin{prop} \label{p1}
Soit $[\Lambda,n,0,\beta]$ une strate semi-simple et $\PC=\MC\UC$ un \para de $\GC=\GC\LC(V)$ tel que $M$ contienne le tore $F[\beta]^\times$ et $B(\MC,F)$ contienne le point $x$ de $B(\GC,F)$ associ\'e \`a $\Lambda$. %Alors, le groupe $H^+(\Lambda,\beta)$ admet une $(\o{P}.P)$ d\'ecomposition d'Iwahori et pour tout caract\`ere $\theta\in \CC(\Lambda,0,\beta)$, les restrictions $\theta_{|H^+\cap \o{U}}$ et $\theta_{H^+\cap U}$ sont triviales. 
Soit $\theta\in\CC(\Lambda,0,\beta)$,   $\theta_M$ sa restriction  \`a $H^+(\Lambda,\beta)\cap M$ et $\ve_{\theta_M}$ l'idempotent de $RM_x$ associ\'e, o\`u $R=\ZM_{p-cycl}[\frac{1}{p}]$.
On a 
$$e_{U_x^+}\eubx \ve_{\theta_M} \in RG_x \eux \eubx \ve_{\theta_M} .$$

\end{prop}

Partons maintenant d'un \para $\PC=\MC\UC$ dans $\GC$ et %et d'un point $x\in B(M,F)$. 
notons $V=\bigoplus_{i\in I}V_i$  la d\'ecomposition de $V$ associ\'ee au \levi $\MC$. Donnons-nous pour chaque $i$ une strate semi-simple $[\Lambda_i,n_i,0,\beta_i]$ dans $\endo{F}{V_i}$ et un caract\`ere semi-simple $\theta_i\in\CC(\Lambda_i,0,\beta_i)$.
La collection des $\Lambda_i$ correspond \`a un point de $B(\MC,F)$ et la collection des caract\`eres simples nous fournit un idempotent $\ve \in RM_x$ que nous qualifierons de {\em semi-simple}. 

\begin{prop} %(Stevens, Bushnell-Kutzko) 
\label{p2}
Soit $\Lambda:=\bigoplus_{i\in I} \Lambda_i$. Il existe une strate semisimple $[\Lambda,n,0,\beta]$ avec $F[\beta]^\times \subset M$ et un caract\`ere semi-simple $\theta\in\CC(\Lambda,0,\beta)$ tel que $\ve_{\theta_M}=\ve$.
\end{prop}

En appliquant ce r\'esultat aux collections translat\'ees $\Lambda_i[t_i]$, avec $t_i\in \QM$, on en d\'eduit que l'idempotent $\ve$ est $P$-bon au sens de \ref{pbon}.

\begin{prop} (Stevens, Bushnell-Kutzko) \label{p3}
La famille des idempotents semi-simples engendre la cat\'egorie $\Mo{R}{M}$, o\`u $R=\ZM_{p-cycl}[{1\over p}]$.
\end{prop}

Avant de prouver ces trois propositions, expliquons comment descendre ces r\'esultats \`a $\zp$. Remarquons que le groupe $\Gamma_p:=\gal(\QM_{p-cycl}|\QM)$ agit sur les caract\`eres $\psi: F/\PG_F \To{} \ZM_{p-cycl}^\times$. Soit $[\Lambda,n,0,\beta]$ une strate semi-simple, on peut donc avec des notations \'evidentes consid\'erer les ensembles 
$$ \o\CC(\Lambda,0,\beta):= \bigcup_{\gamma\in\Gamma_p} \CC_{\psi^\gamma}(\Lambda,0,\beta).$$
Alors la somme
$$\ve_{\Lambda,\beta}:= \sum_{\theta\in \o\CC(\Lambda,0,\beta)} \ve_\theta $$
est un idempotent de $\ZM[\frac{1}{p}]G_x$ (o\`u $x$ correspond \`a $\Lambda$).
La proposition \ref{p3} implique que la famille des idempotents de la forme $\times_{i\in I} \ve_{\Lambda_i,\beta_i} \in \zp M_x$ engendre $\Mo{\zp}{M}$, et les propositions \ref{p1} et \ref{p2} assurent que ces idempotents sont $P$-bons au sens de \ref{pbon}.

\alin{Preuve de la proposition \ref{p1}}
Nous voulons appliquer les \'enonc\'es ``g\'en\'eraux" \ref{constess} et \ref{utile}.
Notons pour cela $A(\beta):=\endo{F[\beta]}{V}\subset A=\endo{F}{V}$, et $\aG_r(\Lambda,\beta):=A(\beta)\cap \aG_r(\Lambda)$ pour $r\in\RM$. On sait alors que $\aG_0(\Lambda,\beta)$ est un $\OC_F$-ordre h\'er\'editaire dans la $F$-alg\`ebre semi-simple $A(\beta)$ et que $\aG_{0+}(\Lambda,\beta)$ est son radical de Jacobson.
De m\^eme notons $\jG_r(\Lambda,\beta):=\jG(\Lambda,\beta)\cap \aG_r(\Lambda)$. On sait que 
$\jG_{0+}(\Lambda,\beta)$ est le radical de Jacobson de $\jG(\Lambda,\beta)$ et on a par d\'efinition l'\'egalit\'e $\jG(\Lambda,\beta)=\jG_{0+}(\Lambda,\beta)+\aG_0(\Lambda,\beta)$.

Lorsqu'on a une $\OC_F$-alg\`ebre finie et plate $B$, le foncteur sur les $\OC_F$-alg\`ebres qui \`a $R$ associe $(B\otimes R)^\times$ est repr\'esentable par un sch\'ema en groupes affine lisse sur $\OC_F$ (un ouvert d'un espace affine sur $\OC_F$). Nous noterons $\GC_x$, resp. $\GC_{\beta,x}$, le sch\'ema en groupes associ\'e \`a $\aG_0(\Lambda)$, resp. $\aG_0(\Lambda,\beta)$ et $\JC$ celui associ\'e \`a $\jG(\Lambda,\beta)$.
Les relations de contenance de ces ordres 
%et leurs propri\'et\'es rappel\'ees ci-dessus 
induisent d'une part
un morphisme  $\varphi: \JC \To{} \GC_x$ qui sur la fibre g\'en\'erique induit l'identit\'e de $\GC\LC(V)$ et 
%dont la fibre sp\'eciale v\'erifie $\varphi_{k_F}^{-1}({^u\GC_{x,k_F}})={^u\JC_{k_F}}$, et
d'autre part une immersion ferm\'ee $\psi: \GC_{\beta,x}\To{} \JC$.
Puisque $\aG_{0+}(\Lambda)\cap \jG(\Lambda,\beta)$ est le radical de Jacobson de $\jG(\Lambda,\beta)$, on a $\varphi_k^{-1}({^u\GC_{x,k}})={^u\JC_k}$ ou, de mani\`ere \'equivalente, $\JC^\dag(\OC_F)=J^+(\Lambda,\beta):=G_x^+\cap J(\Lambda,\beta)$.
Puisque $\jG(\Lambda,\beta)=\jG_{0+}(\Lambda,\beta)+\aG_0(\Lambda,\beta)$, la
 fibre sp\'eciale de $\psi$ induit un isomorphisme des quotients r\'eductifs ${^q\JC_{k_F}}\simto {^q\GC_{\beta,x,k_F}}$. Notons que, par la d\'efinition d'une strate semi-simple, la fonction-r\'eseau $\Lambda$ d\'efinit un point $x_\beta$ de l'immeuble $B(\GC_\beta,x)$ associ\'e au groupe r\'eductif $\GC_\beta$ des inversibles de la $F$-alg\`ebre $A(\beta)$ (le centralisateur de $\beta$).
Le groupe $\GC_{\beta,x}$ s'identifie au  groupe parahorique de $\GC_\beta$ associ\'e \`a $x_\beta$. En particulier, comme $\GC_\beta$ est un groupe lin\'eaire, $\GC_{\beta,x}$ est connexe (ainsi que $\GC_x$), et par cons\'equent $\JC$ l'est aussi.

%Notons maintenant que $F[\beta]^\times$ est le (groupe des points $K$-rationnels du) centre de $\GC_\beta$.
Maintenant l'hypoth\`ese $M\supset F[\beta]^\times$ implique que le centre $\ZC(\MC)$ de $\MC$ est inclus dans $\GC_\beta$, et par cons\'equent l'intersection $\MC_\beta:=\MC\cap \GC_\beta=\ZC_{\GC_\beta}(\ZC(\MC))$ est un sous-groupe de Levi de $\GC_\beta$. %et $\PC_\beta:=\PC\cap \GC_\beta$ un sous-groupe parabolique de $\GC_\beta$ de composante de Levi $\MC_\beta$. 
L'hypoth\`ese $x\in B(\MC,F)$ implique que $x_\beta\in B(\MC_\beta,F)$ et par la remarque \ref{remapara}, que l'adh\'erence sch\'ematique de $\ZC(\MC)$ dans $\GC_{\beta,x}$ est un tore. 
Il en est donc de m\^eme de l'adh\'erence sch\'ematique de $\ZC(\MC)$ dans $\JC$,
 autrement dit  $\MC$ est $\JC$-admissible au sens de \ref{defmod}.

R\'eciproquement, partons d'un \levi $\LC$ dans $\GC$ qui est $\JC$-admissible 
 et choisissons un tore $\TC$ d\'eploy\'e dont le centralisateur est $\LC$ et qui se prolonge en un tore $\TC_x$ de $\JC$. Alors d'une part $\TC_x$ s'envoie sur un sous-tore de $\GC_x$ (car le noyau de $\JC_k \To{} \GC_{x,k}$ est unipotent) donc $\LC$ est $\GC_x$-admissible et, par \ref{remapara}, on a $x\in B(\LC,F)$. D'autre part, comme l'immersion ferm\'ee $\GC_{\beta,x,k_F}\To{} \JC_{k_F}$ induit un isomorphisme des quotients r\'eductifs, le tore $\TC_{x,k_F}$ est contenu dans $\GC_{\beta,x,k_F}$. Par \cite[Exp IX, Thm 3.6bis]{SGA3}, il s'ensuit que $\TC_x$ est conjugu\'e par un \'el\'ement de $J(\Lambda,\beta)$ \`a un tore de $\GC_{\beta,x}$, et par cons\'equent que $\LC$ contient un conjugu\'e sous $J(\Lambda,\beta)$ du tore $F[\beta]^\times$.
On d\'eduit aussi du paragraphe pr\'ec\'edent que le centre  de $\LC$ tout entier se prolonge en un sous-tore ferm\'e de $\JC$.

Rappelons maintenant les r\'esultats suivants de la th\'eorie des types :

\begin{fact} \label{fait2} (Stevens) Avec les hypoth\`eses et notations de la proposition \ref{p1}, \begin{enumerate}
        \item $J(\Lambda,\beta)$ normalise $H^+(\Lambda,\beta)$ et $\theta$.         \item $H^+(\Lambda,\beta)$ a la d\'ecomposition d'Iwahori par rapport \`a $P,\o{P}$ et les restrictions de tout caract\`ere semi-simple dans $\CC(\Lambda,0,\beta)$ \`a $H^+\cap U$, $H^+\cap \o{U}$ sont triviales.
        \item $[J^+(\Lambda,\beta),J^+(\Lambda,\beta)]\subseteq H^+(\Lambda,\beta)\subseteq J^+(\Lambda,\beta)$ et l'application $(u,\o{u})\mapsto \theta([u,\o{u}])$ induit un accouplement non-d\'eg\'en\'er\'e 
        $$ (J^+\cap U)/(H^+\cap U)\times (J^+\cap \o{U})/(H^+\cap \o{U}) \To{} R^\times. $$
         \item \'Ecrivons $V=\bigoplus V_i$ la d\'ecomposition associ\'ee \`a $M$ et $\beta=\oplus_i \beta_i$ la d\'ecomposition de $\beta$ correspondante. Alors $H^+(\Lambda,\beta)\cap M =\prod_i H^+(\Lambda_i,\beta_i)$ et la restriction de tout caract\`ere semi-simple dans $\CC(\Lambda,0,\beta)$ \`a $H^+\cap M$ est un produit de caract\`eres semi-simples dans $\CC(\Lambda_i,0,\beta_i)$.
           \item  L'ensemble d'entrelacement $\hbox{Int}_{G_x}(\theta)$ de $\theta$ dans $G_x$ est $J(\Lambda,\beta)$.
      \end{enumerate}
\end{fact}
\begin{proof} (r\'ef\'erences et commentaires) le premier point est prouv\'e en \cite[Coro 3.12 (iii)]{Stevens} et \cite[Lemma 3.15.(iii)]{Stevens}, et le dernier point d\'ecoule de \cite[Thm 3.22]{Stevens} (qui calcule l'entrelacement dans tout $G$).
Le point iv) est une cons\'equence \`a peu pr\`es directe de la  d\'efinition \cite[3.13]{Stevens} d'un caract\`ere semi-simple et de la proposition 3.4 de \cite{Stevens}. La d\'ecomposition d'Iwahori du groupe $H^+$ et des caract\`eres semisimples dans le point ii) d\'ecoule de leur d\'efinition inductive, suivant  le m\^eme argument  que \cite[7.1.19]{BK} (qui est le cas simple),  \cite[Prop 5.2.ii)]{BK2}  ou \cite[Lemma 3.15 i)]{Stevens}.
Enfin le point iii) se d\'eduit de \cite[Prop. 3.24]{Stevens} et du point ii) suivant la m\^eme observation que \cite[7.2.3 (i)]{BK}.

\end{proof}

Fixons un parabolique $\JC$-admissible minimal $\QC=\LC\VC$ de $\GC$ ; par la discussion pr\'ec\'edente on a $x\in B(\LC,F)$ et on peut supposer, quitte \`a conjuguer, que $F[\beta]^\times \subset \LC$. On peut donc appliquer  les points i), ii) et iii) ci-dessus  \`a $Q$ \`a la place de $P$, et obtenir gr\^ace \`a \ref{constess}  que tout idempotent de $RJ(\Lambda,\beta)$ associ\'e \`a un caract\`ere semi-simple  est ``essentiellement de niveau z\'ero", au sens de \ref{defniv0}. D'apr\`es le point iv) ci-dessus, l'idempotent $\ve_{\theta_M}$ de $R(J(\Lambda,\beta)\cap M)$ est donc essentiellement de niveau $0$ pour le mod\`ele lisse de $\MC$ obtenu comme adh\'erence sch\'ematique de $\MC$ dans $\JC$.
Mais d'apr\`es le point ii) \`a nouveau et le point v), les hypoth\`eses du corollaire \ref{utile} sont satisfaites pour le morphisme $\varphi : \JC\To{} \GC_x$  et pour $\wt\ve'=\ve_\theta$. Il ne reste donc plus qu'\`a appliquer ce corollaire.

\alin{Preuve de la proposition \ref{p2}}
Ce r\'esultat ne figure explicitement ni dans \cite{Stevens}, ni dans \cite{BK2}, mais d\'ecoule pourtant des techniques d\'evelopp\'ees dans ces deux articles. La discussion qui suit entend en convaincre le lecteur d\'eja familier de ces techniques.

Introduisons la famille de sous-groupes $H^r(\Lambda,\beta):=\hG(\Lambda,\beta)\cap \uG_r(\Lambda)$ %et $J^r(\Lambda,\beta):=\jG(\Lambda,\beta)\cap \uG_r(\Lambda)$  
pour $r\in \RM_+$. Cette famille est d\'ecroissante et ses sauts sont dans le mono\"ide discret $\frac{1}{e(\Lambda)}\NM$. Pour $r\in \RM$, nous conviendrons de noter $r+$, resp. $r-$ le plus petit saut  strictement plus grand que $r$, resp. le plus grand saut strictement plus petit que $r$.
 Lorsque $r$ est inf\'erieur \`a l'entier $k_0(\beta,\Lambda)=k_0(\beta)$ d\'efini en \cite[(3.6)]{Stevens}, Stevens d\'efinit \cite[3.13]{Stevens} un ensemble  de caract\`eres complexes $\CC(\Lambda,r,\beta)$ du groupe $H^{r+}(\Lambda,\beta)$, pour $r\in\RM_+$. D'apr\`es \cite{Stevens}, Rk 3.14.ii) et Lemma 3.15.i), les applications de restriction $\CC(\Lambda,r,\beta)\To{} \CC(\Lambda,r',\beta)$ pour $0\leq r \leq r' <k_0(\beta)$ sont surjectives. On peut donc prolonger la notation \`a  tout $r\in\RM_+$ en d\'efinissant  $\CC(\Lambda,r,\beta)$ comme l'ensemble des restrictions \`a $H^{r+}(\Lambda,\beta)$ des caract\`eres dans $\CC(\Lambda,0,\beta)$. En particulier, pour $r\geq \frac{n}{2}$, on a $H^{r+}(\Lambda,\beta)=\uG_{r+}(\Lambda)$ et $\CC(\Lambda,r,\beta)=\{{\psi_\beta}_{|H^{r+}}\}$, 
o\`u $\psi_\beta:\; \uG_{\frac{n}{2}+}(\Lambda) \To{}\ZM_{p-cycl}[\frac{1}{p}]^\times$ est le caract\`ere $x\mapsto \psi(\hbox{Tr}(\beta(x-1))$ associ\'e \`a $\beta$ et $\psi$. 
De plus, 
par \cite[3.14.(ii)]{Stevens}, on a  $\CC(\Lambda,r,\beta)=\CC(\Lambda,r,\gamma)$ pour toute strate  $[\Lambda,n,r,\gamma]$  semi-simple, \'equivalente \`a $[\Lambda,n,r,\beta]$ et telle que $F[\gamma]^\times\subset M(\beta)$ o\`u $M(\beta)$ est le Levi de $G$ d\'ecoup\'e par l'alg\`ebre semi-simple $F[\beta]\subset A$. 

Revenons \`a  l'\'enonc\'e de \ref{p2}. Posons $n:=\max(n_{\beta_i})_{i\in I}$ et notons $M$ le \levi associ\'e \`a la d\'ecomposition $\Lambda=\bigoplus_{i\in I} \Lambda_i$.
La strate semi-simple $\prod_i [\Lambda_i,n_{\beta_i},0,\beta_i]$ pour $M$ d\'ecoupe un Levi $M(\beta) \subset M$ dont nous noterons $V=\bigoplus_{j\in J_\beta} V_j$ la d\'ecomposition associ\'ee. %, {\em cf} \cite[Rk 3.3.i)]{Stevens}. 
On a donc une application surjective $J_\beta\To{} I$ qui \`a $j$ associe l'unique $i(j)$ tel que $V_j\subset \Lambda_{i(j)}\otimes F$,  et des d\'ecompositions $\Lambda_i=\bigoplus_{j\mapsto i} \Lambda_i\cap V_j$.

Nous allons prouver par r\'ecurrence descendante sur $n\geq t\geq 0$ (nombre fini de sauts !) l'assertion suivante :

{\em Il existe une strate semi-simple $[\Lambda,n,t,\gamma^{t}]$ avec $F[\gamma^{t}]^\times \subset M(\beta) \subset M(\gamma^t)$, et un caract\`ere semi-simple $\theta^t\in\CC(\Lambda,0,\gamma^t)$ tels que 
\begin{itemize}
        \item $H^{t+}(\Lambda,\gamma^t)\cap M = \prod_i H^{t+}(\Lambda_i,\beta_i)$,
        \item $\theta^t_{|H^{t+}(\Lambda,\gamma^t)\cap M}= \prod_i {\theta_i}_{|H^{t+}(\Lambda_i,\beta_i)}$
\end{itemize}
}
Le premier saut est $t=n$, pour lequel il suffit de prendre la strate nulle $[\Lambda,n,n,\gamma^n=0]$ et le caract\`ere trivial de $H^+(\Lambda,\gamma^n)=\uG_{0+}(\Lambda)$.

Supposons donc l'\'enonc\'e connu pour $t$ et d\'eduisons-le pour $t-$. \'Ecrivons $\gamma^t=\bigoplus_{i\in I} \gamma^t_i$ la d\'ecomposition de $\gamma^t$ comme \'el\'ement de $M$. Comme en \ref{fait2} iv),  pour tout $r\in\RM_+$ on a $H^{r+}(\Lambda,\gamma^t)\cap M = \prod_i H^{r+}(\Lambda_i,\gamma_i^t)$ et si $\theta\in \CC(\Lambda,r,\gamma^t)$, alors $\theta_{|H^{r+}(\Lambda,\gamma^t)\cap M}$ est un produit sur $i$ de caract\`eres semi-simples dans $\CC(\Lambda_i,r,\gamma_i^t)$.
%Comme dans \cite[8.4]{BK2}, on d\'eduit donc de \cite[3.5.9]{BK}
On d\'eduit alors du lemme \ref{lbk} ci-dessous
 et de l'hypoth\`ese de r\'ecurrence que  $H^{t}(\Lambda,\gamma^t)\cap M = \prod_i H^{t}(\Lambda_i,\beta_i)$. Il existe donc un \'el\'ement $b= \bigoplus_i b_i\in \bigoplus_i \aG_{-t}(\Lambda_i)\subset \aG_{-t}(\Lambda)$ tel que 
$$ \theta^t_{|H^t(\Lambda,\gamma^t)\cap M} .{\psi_b}_{|H^t(\Lambda,\gamma^t)\cap M} = \prod_{i\in I} {\theta_i}_{|H^t(\Lambda_i,\beta_i)}, $$
o\`u $\psi_b$ est le caract\`ere de $\uG_t(\Lambda)/\uG_{t+}(\Lambda)$ associ\'e \`a $\psi$ et $b$. En fait, par d\'ecomposition d'Iwahori \ref{fait2} ii)4 des caract\`eres semi-simples, on peut supposer que chaque $b_i$ se d\'ecompose en $b_i=\bigoplus_{j\mapsto i} b_j$ avec $b_j\in \aG_{-t}(\Lambda_i\cap V_j)$ et $j\in J_\beta$, de sorte que pour tout $j\in J_\beta$, on a
$$\theta^t_{|H^t(\Lambda_j,\gamma^t_j)} .{\psi_{b_j}}_{|H^t(\Lambda_j,\gamma^t_j)} = {\theta_{i(j)}}_{|H^t(\Lambda_j,\beta_j)}$$
o\`u les deux $\theta$ ainsi restreints sont des caract\`eres {\em simples}.

Soit alors $V=\bigoplus_{k\in K_t} V_k$ la d\'ecomposition d\'etermin\'ee par l'alg\`ebre semi-simple $F[\gamma^t]$ (correspondant au Levi $M(\gamma^t)$). Puisque $M(\gamma^t)\supset M(\beta)$, on a une application surjective $J_\beta \To{} K_t$ qui \`a $j$ associe l'unique $k(j)$ tel que $V_j\subset V_{k(j)}$. On a aussi la d\'ecomposition en produit de corps $F[\gamma^t]\simeq \prod_{k\in K_t} E_k$. Choisissons alors pour chaque $k\in K_t$ une corestriction mod\'er\'ee $s_k :\endo{F}{V_k} \To{} \endo{E_{k}}{V_k}$. Puisque $F[\gamma^t]^\times\subset M(\beta)$, celle-ci induit par restriction une corestriction mod\'er\'ee $s_j :\endo{F}{V_j} \To{} \endo{E_{k}}{V_j}$ pour tout $j$ tel que $k=k(j)$.
D'apr\`es  \cite[4.6]{BK2}, la $E_{k(j)}$-strate $[\Lambda_j,t,t-,s_j(b_j)]$ est \'equivalente \`a une strate simple, \'eventuellement nulle.
On en d\'eduit que pour  tout $k\in K_t$, la $E_k$-strate $[\Lambda_k,t,t-,s_k(b_k)]$, o\`u  $\Lambda_k=\bigoplus_{j\mapsto k} \Lambda_j$ et $b_k:= \bigoplus_{j\mapsto k} b_j$ est \'equivalente \`a une strate semi-simple. Par \cite[3.5]{Stevens}, il s'ensuit que la strate $[\Lambda,n,t-,\gamma^t+b]$ est \'equivalente \`a une strate semi-simple, disons $[\Lambda,n,t-,\gamma^{t-}]$. En outre, comme dans la preuve de \cite[3.4]{Stevens}, on peut choisit $\gamma^{t-}$ tel que  $F[\gamma^{t-}]^\times\subset M(\beta)$. D'apr\`es \cite[Rk 3.14.(i)]{Stevens}, on a $H^t(\Lambda,\gamma^{t-})=H^t(\Lambda,\gamma^t)$ et   $\CC(\Lambda,t-,\gamma^{t-}) = \psi_b. \CC(\Lambda,t-,\gamma^t)$. Comme l'application de restriction des caract\`eres induit une surjection $\CC(\Lambda,0,\gamma^{t-})\To{} \CC(\Lambda,t-,\gamma^{t-})$ on peut choisir un $\theta^{t-}\in \CC(\Lambda,0,\gamma^{t-})$ tel que
$$ \theta^{t-}_{|H^t(\Lambda,\gamma^{t-})} = (\psi_b\theta^t)_{|H^t(\Lambda,\gamma^t)}.$$
Celui-ci remplit le cahier des charges. 

Nous avons utilis\'e dans cette preuve le lemme suivant qui est une g\'en\'eralisation au cas semi-simple de \cite[3.5.9]{BK}.

\begin{lemme} \label{lbk}
Soient $[\Lambda,n,0,\beta_i]$, $i=1,2$, deux strates semi-simples et $r> 0$ telles que 
%\begin{enumerate}
        %\item 
        $H^{r+}(\Lambda,\beta_1)=H^{r+}(\Lambda,\beta_2)$,
        et $\CC(\Lambda,r,\beta_1)\cap \CC(\Lambda,r,\beta_2)\neq \emptyset$.
%\end{enumerate}
Alors $H^r(\Lambda,\beta_1)=H^r(\Lambda,\beta_2)$.
\end{lemme}

\begin{proof}
Choisissons deux strates $[\Lambda,n,2r,\gamma_i]$ semi-simples avec $\gamma_i\in M(\beta_i)$, et respectivement \'equivalentes \`a $[\Lambda,n,2r,\beta_i]$. On sait alors que pour tout $t\geq r$, on a $\hG^t(\Lambda,\beta_i) = \hG^t(\Lambda,\gamma_i)$ et $\CC(\Lambda,2t,\beta_i)=\CC(\Lambda,2t,\gamma_i)$. En particulier, puisque pour tous $t\geq t'$ l'application de restriction $\CC(\Lambda,t',\beta_i)\To{}\CC(\Lambda,t,\beta_i)$ est surjective, on a  $\CC(\Lambda,2r,\gamma_1)\cap \CC(\Lambda,2r,\gamma_2)\neq \emptyset$.
Soit $\theta$ un \'el\'ement de cette intersection, le th\'eor\`eme 3.22 de \cite{Stevens} calcule l'ensemble d'entrelacement de $\theta$ dans $G$ et nous fournit l'\'egalit\'e
$$ \Gamma_{2r}(\Lambda,\gamma_1) G_{\gamma_1} \Gamma_{2r}(\Lambda,\gamma_1) =\Gamma_{2r}(\Lambda,\gamma_2) G_{\gamma_2} \Gamma_{2r}(\Lambda,\gamma_2) $$
%, donc en particulier aussi celui dans $u_0(\Lambda)$. On en d\'eduit  
avec les notations de {\em loc. cit}.
% que $\Gamma_{2r}(\Lambda,\gamma_1)(u_0(\Lambda)\cap A_{\gamma_1}) =\Gamma_{2r}(\Lambda,\gamma_2)(u_0(\Lambda)\cap A_{\gamma_2})$.
En prenant l'intersection avec $\aG_r(\Lambda)$ et en prenant la cl\^oture additive, on obtient l'ind\'ependance de $i$ de l'ensemble suivant :
\ini\begin{equation}\label{ens}
 (\aG_r(\Lambda)\cap A_{\gamma_i}) + (\aG_r(\Lambda)\cap 
 A_{\gamma_i})(\nG_{-2r}(\Lambda,\gamma_{i})\cap \aG_{r_i-2r}(\Lambda)) + (\aG_r(\Lambda)\cap A_{\gamma_i}) \jG^{\frac{r_i}{2}}(\Lambda,\gamma_i) 
\end{equation}
o\`u on a pos\'e $r_i:=k_0(\gamma_i,\Lambda)$ ({\em cf} \cite[(3.6)]{Stevens}).
Soit $r'_i:= 2r-\frac{r_i}{2}<r$. On a
\begin{eqnarray*}
(\aG_{r'_i}(\Lambda)\cap A_{\gamma_i})(\nG_{-2r}(\Lambda,\gamma_{i})\cap \aG_{r_i-2r}(\Lambda)) 
& \subset & (\nG_{-\frac{r_i}{2}}(\Lambda,\gamma_{i})\cap \aG_{\frac{r_i}{2}}(\Lambda)) \\
& \subset & \jG^{\frac{r_i}{2}}(\Lambda,\gamma_i) \\
\end{eqnarray*}
par \cite[3.10.i)]{Stevens}, donc 
\begin{eqnarray*}
(\aG_r(\Lambda)\cap A_{\gamma_i})(\nG_{-2r}(\Lambda,\gamma_{i})\cap \aG_{r_i-2r}(\Lambda)) & \subset & (\aG_{r-r'_i}(\Lambda)\cap A_{\gamma_i})\jG^{\frac{r_i}{2}}(\Lambda,\gamma_i) \\
& \subset & \hG^{\frac{r_i}{2}+}(\Lambda,\gamma_i) \subset \hG^{r+}(\Lambda,\gamma_i) \\
\end{eqnarray*}  
par \cite[3.11.i)]{Stevens}. De plus, par \cite[3.11.ii)]{Stevens}, on a $(\aG_r(\Lambda)\cap A_{\gamma_i}) \jG^{\frac{r_i}{2}}(\Lambda,\gamma_i)\subset \hG^{r+}(\Lambda,\gamma_i)$.
Ajoutons alors \`a l'ensemble \ref{ens} le groupe  $\hG^{r+}(\Lambda,\gamma_i)$ qui par hypoth\`ese est aussi ind\'ependant de $i$. On obtient que l'ensemble
$  (\aG_r(\Lambda)\cap A_{\gamma_i}) + \hG^{r+}(\Lambda,\gamma_i)$ est ind\'ependant de $i$.
Mais celui-ci n'est autre que $\hG^{r}(\Lambda,\gamma_i)$ puisque $r_i>2r$.
\end{proof}

Pour la preuve de la proposition \ref{p3}, on renvoie \`a celle de la proposition \ref{p3c} au paragraphe \ref{pfp3c}.

\section{Groupes classiques} \label{classique}

Nous adoptons les notations de Stevens dans \cite{Stevens}. Cette fois le corps de base, que nous notions $K$ pr\'ec\'edemment sera not\'e $F_0$ et sera suppos\'e \^etre le corps des points fixes d'une involution $x\mapsto\o{x}$ sur un corps $F$ de caract\'eristique r\'esiduelle diff\'erente de $2$. On n'\'ecarte pas le cas o\`u $F_0=F$ et l'involution est l'identit\'e. Soit $V$ un $F$-espace vectoriel muni d'une forme bilin\'eaire $h$ $\ve$-hermitienne (pour $\ve=\pm 1$) non-d\'eg\'en\'er\'ee. La $F$-alg\`ebre $A$ est munie de  l'anti-involution ``adjoint pour $h$" qui prolonge l'involution donn\'ee sur $F$ et que nous noterons encore $x\mapsto \o{x}$. Notons $\wt\GC$ le $F$-groupe $\GC\LC(V)$ ;
alors cette anti-involution munit le $F_0$-sch\'ema en groupes $\hbox{Res}_{F|F_0}(\wt\GC)$ d'une involution $\sigma$ dont le sous-sch\'ema des points fixes $\GC$ s'identifie au
$F_0$-sch\'ema en groupes unitaire, orthogonal ou symplectique  associ\'e \`a $(V,h)$.

\alin{Immeuble et fonctions r\'eseaux autoduales} \label{introclassique}
La r\'ef\'erence ici est \cite{BrSt}.
Pour un $\OC_F$-r\'eseau $L$ de $V$, on note $L^\#:=\{v\in V, \; h(v,L)\subset \PG_F\}$, et pour une fonction r\'eseau $\Lambda\in \FC\RC_F(V)$, on note $\Lambda^\#$ la fontion r\'eseau $r\mapsto \Lambda((-r)+)^\#$. La bijection naturelle $B(\hbox{Res}_{F|F_0}(\wt\GC),F_0)=B(\wt\GC,F)\simto \FC\RC_F(V)$ est compatible avec les involutions $\sigma$ sur $B(\hbox{Res}_{F|F_0}(\wt\GC),F_0)$ et $\#$ sur $\FC\RC_F(V)$. Par l'hypoth\`ese de caract\'eristique r\'esiduelle $\neq 2$, elle induit en prenant les invariants une application bijective et $G$-\'equivariante $B(\GC,F_0) \To{} \FC\RC_F(V)^\#$. % o\`u  

Si  $\Lambda\in \FC\RC_F(V)^\#$, on a $\o{\aG_r(\Lambda)}=\aG_r(\Lambda)$, resp.  $\sigma(\uG_r(\Lambda))=\uG_r(\Lambda)$, pour tout $r\in\RM$, resp. $r\in\RM_+$. 
Si $x\in B(\GC,F_0)$ correspond \`a $\Lambda$, alors $G_x= G\cap \uG_0(\Lambda) = \uG_0(\Lambda)^\sigma$ et $G_x^+ = G\cap \uG_{0+}(\Lambda)=\uG_{0+}(\Lambda)^\sigma$, et les sp\'ecialistes s'accordent \`a penser que la filtration $(\uG_r(\Lambda)^\sigma)_{r\in\RM_+}$ co\"incide avec celle de Moy et Prasad $(G_{x,r})_{r\in\RM}$ (ce qui ne nous importe gu\`ere ici). 
Un sous-groupe de Levi $\MC$ de $\GC$ correspond \`a une d\'ecomposition $h$-orthogonale $V=\bigoplus_{i\in I} (V_i\oplus V_{-i}) \bigoplus V_0$ o\`u $h_{|V_0\times V_0}$ est non-d\'eg\'en\'er\'ee et pour chaque $i\in I$, $V_i$ est totalement isotrope et $h_{|V_i\times V_{-i}}$ est un accouplement parfait. On a alors $M:=\MC(F_0)\simeq \prod_{i\in I} \aut{F}{V_i} \times \aut{F}{V_0}^\sigma$. Le sous-immeuble $B(\MC,F_0)$ correspond aux $\Lambda$ qui sont d\'ecompos\'ees sous la forme $\Lambda=\bigoplus_{i\in I} (\Lambda_i\oplus \Lambda_{-i})\bigoplus \Lambda_0$ o\`u $\Lambda_0\in \FC\RC_F(V_0)^\#$ et les $\Lambda_{\pm i} \in \FC\RC_F(V_{\pm i})$ sont telles que $\Lambda_{-i}=\Lambda_i^\#$.
Si $x$ est le point correspondant \`a $\Lambda$, alors le sous-espace affine $x+ a_M$ de $B(\MC,F_0)$ est l'ensemble des $\Lambda'$ de la forme $\bigoplus_{i\in I} (\Lambda_i[t_i]\oplus \Lambda_i[t_i]^\#) \bigoplus \Lambda_0$ o\`u $t_i\in \RM$.

\alin{Strates semi-simples autoduales}
Une $F$-strate $[\Lambda,n,r,\gamma]$ dans $V$ est dite autoduale si $\Lambda=\Lambda^\#$ et $\o\gamma=-\gamma$. 
Soit $[\Lambda,n,r,\beta]$ une strate semi-simple et autoduale. La sous-alg\`ebre $F[\beta]\subset A$ est stable par $x\mapsto \o{x}$ et l'ensemble de ses idempotents centraux primitifs aussi. On peut donc arranger la d\'ecomposition de $V$ associ\'ee \`a $\beta$ sous la forme $V=\bigoplus_{i\in I_\beta} (V_i\oplus V_{-i}) \bigoplus_{j\in J_\beta} V_j$ o\`u les espaces index\'es par $I_\beta\cup J_\beta$ sont deux \`a deux orthogonaux et pour $i\in I_\beta$, $V_i$ et $V_{-i}$ sont isotropes maximaux dans $V_i\oplus V_{-i}$.
Suivant Stevens, la strate semisimple autoduale $[\Lambda,n,r,\beta]$ est dite {\em gauche} (skew en anglais), si $I_\beta=\emptyset$, c'est-\`a-dire si $\beta$ est elliptique. 

Dans \cite[3.6]{Stevens}, Stevens d\'efinit les caract\`eres semi-simples pour $G$ associ\'es \`a une strate semi-simple gauche $[\Lambda,n,0,\beta]$. L'hypoth\`ese gauche n'est pas n\'ecessaire pour cette d\'efinition, il suffit de supposer la strate autoduale ; le point est que pour tout $r>0$ , on peut trouver une  strate semi-simple {\em autoduale} $[\Lambda,n,r,\gamma]$  \'equivalente \`a $[\Lambda,n,r,\beta]$, avec de plus  $\gamma$ dans le Levi de $GL(V)$ d\'ecoup\'e par $\beta$ (combiner \cite[3.4]{Stevens} et  \cite[(1.10)]{StPLMS}).
Il s'ensuit que les ordres $\hG(\Lambda,\beta)$ et $\jG(\Lambda,\beta)$ sont stables par l'involution $x\mapsto \o{x}$, les groupes $H^{r+}(\Lambda,\beta)$ et $J^{r+}(\Lambda,\beta)$ sont stables par $\sigma$, ainsi que les ensembles de caract\`eres semi-simples $\CC(\Lambda,r,\beta)$. Un caract\`ere semi-simple de $H^{r+}(\Lambda,\beta)^\sigma$ est alors, par d\'efinition, la restriction d'un caract\`ere semi-simple $\sigma$-invariant de $H^{r+}(\Lambda,\beta)^\sigma$.

\begin{prop} \label{p1c}
Soit $[\Lambda,n,0,\beta]$ une strate semi-simple autoduale et $\PC=\MC\UC$ un \para de $\GC$ tel que $M$ contienne le ``tore" $(F[\beta]^\times)^\sigma$ et $B(\MC,F_0)$ contienne le point $x$ de $B(\GC,F_0)$ associ\'e \`a $\Lambda$. %Alors, le groupe $H^+(\Lambda,\beta)$ admet une $(\o{P}.P)$ d\'ecomposition d'Iwahori et pour tout caract\`ere $\theta\in \CC(\Lambda,0,\beta)$, les restrictions $\theta_{|H^+\cap \o{U}}$ et $\theta_{H^+\cap U}$ sont triviales. 
Soit $\theta\in\CC(\Lambda,0,\beta)^\sigma$,   $\theta_M$ sa restriction  \`a $H^+(\Lambda,\beta)^\sigma\cap M$ et $\ve_{\theta_M}$ l'idempotent de $RM_x$ associ\'e, o\`u $R=\ZM_{p-cycl}[\frac{1}{p}]$.
On a 
$$e_{U_x^+}\eubx \ve_{\theta_M} \in RG_x \eux \eubx \ve_{\theta_M} .$$

\end{prop}

Partons maintenant d'un \para $\PC=\MC\UC$ dans $\GC$ et %et d'un point $x\in B(M,F)$. 
notons $V=\bigoplus_{i\in I} (V_i\oplus V_{-i}) \bigoplus V_0$  la d\'ecomposition orthogonale de $V$ associ\'ee au \levi $\MC$. Donnons-nous pour chaque $i\neq 0$ une strate semi-simple $[\Lambda_i,n_i,0,\beta_i]$ dans $\endo{F}{V_i}$ et un caract\`ere semi-simple $\theta_i\in\CC(\Lambda_i,0,\beta_i)$, ainsi qu'une strate semi-simple autoduale $[\Lambda_0,n_0,0,\beta_0]$ et un caract\`ere semi-simple $\theta_0\in\CC(\Lambda_0,0,\beta_0)^\sigma$.
La collection des $\Lambda_i$ correspond \`a un point de $B(\MC,F_0)$ et la collection des caract\`eres semi-simples nous fournit un idempotent $\ve \in RM_x$ que nous qualifierons de {\em semi-simple}. 

\begin{prop}  \label{p2c}
Soit $\Lambda:=\bigoplus_{i\in I} (\Lambda_i\oplus \Lambda_{-i}^\#)\bigoplus \Lambda_0$. Il existe une strate semisimple autoduale $[\Lambda,n,0,\beta]$ avec $(F[\beta]^\times)^\sigma \subset M$ et un caract\`ere semi-simple $\theta\in\CC(\Lambda,0,\beta)^\sigma$ tel que $\ve_{\theta_M}=\ve$.
\end{prop}

En appliquant ce r\'esultat aux collections translat\'ees $\Lambda_i[t_i]$, avec $t_i\in \QM$ pour $i\neq 0$, on en d\'eduit que l'idempotent $\ve$ est $P$-bon au sens de \ref{pbon}.

\begin{prop} (Stevens) \label{p3c}
Les idempotents semi-simples forment une famille g\'en\'eratrice de la cat\'egorie $\Mo{R}{M}$.
\end{prop}

\alin{Preuve de la proposition \ref{p1c}} Comme dans la preuve de la proposition \ref{p1}, 
on veut appliquer \ref{constess} et \ref{utile}. Rappelons que dans ladite preuve, nous avons introduit  et utilis\'e
%Pour cela, reprenons momentan\'ement les notations de cette preuve et notamment 
des morphismes de  $\OC_F$-sch\'emas en groupes lisses 
$$ \wt\GC_{\beta,x}\To{\wt\psi_\beta} \wt\JC\To{\wt\varphi} \wt\GC_{x}$$
(on rajoute ici des $\wt{}$ pour \^etre coh\'erent avec les notations du paragraphe \ref{introclassique}).
 Comme les ordres auxquels ils sont associ\'es, ces sch\'emas en groupes sont munis d'une action semi-lin\'eaire de $\sigma$. Appliquons-leur le foncteur $\hbox{Res}_{\OC_F|\OC_{F_0}}(-)^\sigma$. On sait que la restriction des scalaires pr\'eserve la lissit\'e, et d'apr\`es \cite[3.4]{Edix} et l'hypoth\`ese de caract\'eristique r\'esiduelle $\neq 2$, le passage aux $\sigma$-invariants aussi. %{\em Changeant alors de notations}, on 
 On obtient donc des morphismes de $\OC_{F_0}$-sch\'emas en groupes lisses
$$ \GC_{\beta,x} \To{\psi_\beta}  \JC \To{\varphi} \GC_{x}, $$
le premier \'etant une immersion ferm\'ee et le second induisant un isomorphisme des fibres g\'en\'eriques.
%dont les fibres sp\'eciales et g\'en\'eriques satisfont les m\^emes propri\'et\'es que la suite initiale, et
 Les points entiers sont donn\'es par $\GC_x(\OC_{F_0})=G_x$ et $\JC(\OC_{F_0})=J(\Lambda,\beta)^\sigma$. En particulier $\GC_x$ est le mod\`ele lisse de $\GC$ associ\'e par  Bruhat-Tits  \`a $x\in B(\GC,F)$.
Comme dans le cas lin\'eaire, par d\'efinition d'une strate semi-simple auto-duale, la fonction-r\'eseau $\Lambda$ d\'efinit un point, disons $x_\beta$ de l'immeuble  du centralisateur $\GC_\beta$, et le groupe $G_{\beta,x}$ s'identifie au fixateur de $x_\beta$ (voir aussi \cite[part. 6]{BrSt}). Il s'ensuit que $\GC_{\beta,x}$ est le mod\`ele lisse de $\GC$ associ\'e par Bruhat-Tits \`a $x_\beta$.

Vient maintenant une difficult\'e technique par rapport au cas lin\'eraire : si la restriction des scalaires pr\'eserve la connexit\'e, il n'en va pas de m\^eme du passage aux $\sigma$-invariants. En fait, dans les cas unitaires et symplectiques, o\`u $\GC$ est connexe et simplement connexe, on sait que $\GC_x$ et $\GC_{\beta,x}$ sont connexes et donc, comme on va le voir ci-dessous,   $\JC$ l'est aussi. On peut alors suivre mot pour mot la m\^eme preuve que \ref{p1}, gr\^ace \`a \ref{fait2c} ci-dessous. Dans le cas orthogonal impair, $\GC$ est connexe mais pas simplement connexe et les $\GC_x$ et $\GC_{\beta,x}$ peuvent ne pas \^etre connexes. Dans le cas orthogonal pair, $\GC$ lui-m\^eme n'est d\'eja pas connexe !  Cependant on observe comme dans la preuve de \ref{pfprinc} que les idempotents $e_{U_x}$, $\eubx$ de l'\'enonc\'e de \ref{p1c} vivent dans $RG_x^\circ$, et que par ailleurs $\ve_{\theta_M}\in RJ^+(\Lambda,\beta)=R\JC^\dag(\OC_{F_0})\subset R\JC^\circ(\OC_{F_0})$. On peut donc essayer de raisonner sur les composantes neutres de ces groupes, ce qui permettra d'appliquer les r\'esultats de la partie \ref{mod}.

Commen{\c c}ons par v\'erifier que le Levi $\MC$ de l'\'enonc\'e de \ref{p1c} est $\JC^\circ$-admissible.
Soit $\wt\MC$ le \levi de $\wt\GC$ d\'ecoup\'e par la d\'ecomposition associ\'ee \`a $\MC$.
On a donc $\MC=\hbox{Res}_{F|F_0}(\wt\MC)^\sigma$. D'apr\`es les hypoth\`eses faites sur $x$ et $\MC$ et la preuve de la proposition \ref{p1}, $\wt\MC$ est $\wt\JC$-admissible et m\^eme, plus pr\'ecis\'ement,  $\ZC(\wt\MC)$ se prolonge en un tore d\'eploy\'e de $\wt\JC$. Il s'ensuit que la partie d\'eploy\'ee du centre $\hbox{Res}_{F|F_0}(\ZC(\wt\MC))$ de $\hbox{Res}_{F|F_0}(\wt\MC)$ se prolonge en un $\OC_{F_0}$-tore d\'eploy\'e de $\hbox{Res}_{\OC_{F}|\OC_{F_0}}(\wt\JC)$, puis, passant aux points fixes, on en d\'eduit que la partie d\'eploy\'ee connexe de $\ZC(\MC)$ se prolonge en un tore de $\JC$, et donc que $\MC$ est $\JC^\circ$-admissible. %R\'eciproquement, si un Levi $\NC$ est $\JC^\circ$-admissible

Pour suivre la strat\'egie de la preuve de \ref{p1}, l'analogue de
 \ref{fait2} est :
\begin{fact} \label{fait2c} (Stevens) Avec les hypoth\`eses et notations de la proposition \ref{p1}, \begin{enumerate}
        \item $J(\Lambda,\beta)^\sigma$ normalise $H^+(\Lambda,\beta)^\sigma$ et $\theta$.         \item $H^+(\Lambda,\beta)^\sigma$ a la d\'ecomposition d'Iwahori par rapport \`a $P,\o{P}$ et les restrictions de tout caract\`ere semi-simple dans $\CC(\Lambda,0,\beta)^\sigma$ \`a $H^{+}\cap U$, $H^{+}\cap \o{U}$ sont triviales.
        \item $[J^+(\Lambda,\beta)^\sigma,J^+(\Lambda,\beta)^\sigma]\subseteq H^+(\Lambda,\beta)^\sigma\subseteq J^+(\Lambda,\beta)^\sigma$ et l'application $(u,\o{u})\mapsto \theta([u,\o{u}])$ induit un accouplement non-d\'eg\'en\'er\'e 
        $$ (J^{+}\cap U)/(H^{+}\cap U)\times (J^{+}\cap \o{U})/(H^{+}\cap \o{U}) \To{} R^\times. $$
         \item \'Ecrivons $V=\bigoplus_i (V_i\oplus V_i) \bigoplus V_0$ la d\'ecomposition orthogonale de $V$ associ\'ee \`a $M$, $\beta=\oplus_i (\beta_i\oplus \beta_{-i}) \oplus \beta_0$ la d\'ecomposition de $\beta$ correspondante, et identifions $M$ \`a $\prod_i GL(V_i) \times GL(V_0)^\sigma$. Alors $H^+(\Lambda,\beta)\cap M \simeq \prod_i H^+(\Lambda_i,\beta_i) \times H^+(\Lambda_0,\beta_0)^\sigma$ et la restriction de tout caract\`ere semi-simple dans $\CC(\Lambda,0,\beta)$ \`a $H^+\cap M$ s'identifie \`a  un produit de caract\`eres semi-simples dans $\CC(\Lambda_i,0,\beta_i)$ par un caract\`ere semi-simple autodual dans $\CC(\Lambda_0,0,\beta_0)^\sigma$.

                    \item  L'ensemble d'entrelacement $\hbox{Int}_{G_x}(\theta)$ de $\theta$ dans $G_x$ est $J(\Lambda,\beta)^\sigma$.
      \end{enumerate}
\end{fact}
\begin{proof} (r\'ef\'erences et commentaires) les deux premiers points d\'ecoulent imm\'ediatement des points correspondants de \ref{fait2}.  Le point iii) est prouv\'e dans \cite[3.28]{Stevens}. Le point iv) se prouve comme le point correspondant de \ref{fait2}, en combinant \cite[(1.10)]{StPLMS} et \cite[Prop 3.4]{Stevens}.
Enfin, le dernier point d\'ecoule de \cite[Thm 3.27]{Stevens} qui calcule l'entrelacement dans tout $G$.
\end{proof}

 On remarque que les points ii) iii) et iv) de \ref{fait2c} concernent des objets relatifs \`a $\JC^\circ$ et on peut toujours restreindre le point i) \`a $\JC^\circ(\OC_{F_0})$.
On en d\'eduit en particulier comme dans le cas lin\'eaire que $\ve_{\theta_M}$ est un idempotent essentiellement de niveau z\'ero pour le mod\`ele lisse de $\MC$ obtenu par adh\'erence sch\'ematique dans $\JC^\circ$. 
 
Le point v) nous donne  $\hbox{Int}_{G_x^\circ}(\theta)=J(\Lambda,\beta)^\sigma\cap G_x^\circ$, ce qui, compte tenu de ce que $\UC_x$ est connexe (voir preuve de \ref{pfprinc}) implique $\hbox{Int}_{U_x}(\ve_\theta)= \JC(\OC_{F_0})\cap U_x$.
 Ainsi, pour pouvoir appliquer \ref{utile} au morphisme $\JC^\circ \To{} \GC_x^\circ$, avec $\ve':=\ve_{\theta_M}$ et $\wt\ve':=\ve_\theta$, et terminer la preuve de \ref{p1c} comme dans le cas lin\'eaire, il reste deux choses \`a prouver :
\begin{enumerate}
\item $U_x^+\cap \JC(\OC_{F_0}) \subseteq U_x\cap \JC^\dag(\OC_{F_0})$
\item  $\JC(\OC_{F_0})\cap U_x = \JC^\circ(\OC_{F_0})\cap U_x$.
\end{enumerate}

Pour cela, l'ingr\'edient essentiel est : {\em soit $\NC_k$ un groupe unipotent sur un corps parfait de caract\'eristique $\neq 2$ et $\sigma$ une involution. 
Alors $H^1(\la\sigma\ra,\NC_k)=0$, et si de plus $\NC_k$ est lisse et connexe, alors $\NC_k^\sigma$ l'est aussi.}
Utilisant une s\'erie centrale caract\'eristique, l'assertion sur le $H^1$ se d\'evisse imm\'ediatement au cas ab\'elien  o\`u elle est \'evidente. Pour l'assertion de connexit\'e, on peut commencer par \'etendre les scalaires \`a une cl\^oture alg\'ebrique. Utilisant ensuite une s\'erie centrale caract\'eristique dont les sous-quotients sont des produits de $\GM_a$ \cite[Exp. 4.1.1 iii) et 4.1.5]{SGA3} et la nullit\'e du $H^1$ qu'on vient de v\'erifier, on est ramen\'e  au cas $\NC_k\simeq \GM_a^r$. On peut alors diagonaliser la matrice de $\sigma$ et r\'eduire encore au cas $r=1$ o\`u c'est \'evident.

Prouvons alors ii) : il suffit de voir que  l'adh\'erence sch\'ematique de $\UC$ dans $\JC$ est connexe. Or, celle-ci est la partie $\sigma$-invariante de l'adh\'erence sch\'ematique de $\hbox{Res}_{F|F_0}(\wt\UC)$ dans $\hbox{Res}_{\OC_F|\OC_{F_0}}(\wt\JC)$, laquelle est connexe par \ref{BT} ii), puisque $\hbox{Res}_{\OC_F|\OC_{F_0}}(\wt\JC)$ est connexe. L'assertion de connexit\'e de l'ingr\'edient ci-dessus permet donc de conclure. 
Pour obtenir i), nous prouverons l'assertion plus forte $\varphi^{-1}({^u\GC_{x,k}})= {^u\JC_k}$, ce qui \'equivaut \`a $\JC^\dag(\OC_{F_0})=\JC(\OC_{F_0})\cap G_x^+$.  Posons $\breve\JC:=\hbox{Res}_{\OC_F|\OC_{F_0}}(\wt\JC)$ et $k:=k_{F_0}$ pour all\'eger les notations. La nullit\'e du $H^1$ ci-dessus montre qu'en appliquant le foncteur des $\sigma$-invariants, la suite exacte $^u\breve\JC_k \injo \breve\JC_k \twoheadrightarrow {^q\breve\JC_k}$ reste exacte. L'assertion de connexit\'e, et le fait que le groupe des invariants d'un groupe r\'eductif par une involution est r\'eductif assurent alors que 
${^u\breve\JC_k}^\sigma={^u\JC_k}$. De m\^eme avec des notations similaires, on a ${^u\breve\GC_{x,k}}^\sigma={^u\GC_{x,k}}$. Mais on sait que $\breve\varphi^{-1}({^u\breve\GC_{x,k}} ) = {^u\breve\JC_k}$ par le cas lin\'eaire (par compatibilit\'e entre les radicaux de Jacobson des ordres auxquels sont associ\'es ces groupes), et il ne reste plus qu'\`a prendre les $\sigma$-invariants.

On en d\'eduit aussi au passage que $\psi_\beta$ induit un isomorphisme $ \pi_0(\GC_{\beta,x,k})\simto\pi_0(\JC_k)$, {\em i.e.} que le d\'efaut de connexit\'e de $\JC$ est le m\^eme que celui de $\GC_{\beta,x}$.

\alin{Preuve de la proposition \ref{p2c}} On peut faire le m\^eme raisonnement inductif que pour la preuve de \ref{p2}, en demandant que les strates interm\'ediaires $[\Lambda,n,t,\gamma^t]$ v\'erifient $\gamma^t =-\o{\gamma^t}$ et que les caract\`eres $\theta^t\in\CC(\Lambda,0,\gamma^t)$ soient invariants par $\sigma$. Pour assurer la propri\'et\'e requise des $\gamma^{t-}$ dans la construction inductive, on utilise \cite[(1.10)]{StPLMS}. L'invariance des caract\`eres $\theta^{t-}$ sous $\sigma$ est alors automatique par d\'ecomposition d'Iwahori.

\alin{Preuve de la proposition \ref{p3c}} \label{pfp3c}
Le groupe $M$ est un produit de groupes lin\'eaires et d'un groupe classique, il suffit  donc de traiter chacun de ces  groupes s\'epar\'ement. En fait nous ne traiterons que le cas classique car c'est le cadre dans lequel les outils n\'ecessaires ont \'et\'e  d\'evelopp\'es par Stevens (notamment le lemme 5.4. de \cite{Stevens}). Nous laisserons le lecteur se convaincre que la m\^eme preuve fonctionne dans le cas lin\'eaire, en admettant que les outils correspondants sont encore valables (ils sont en fait plus faciles \`a obtenir et souvent, un analogue ``simple" se trouve dans \cite[Ch. 8.1]{BK}).

Supposons donc $M=G$,  $R=\ZM_{p-cycl}[\frac{1}{p}]$ et $V\in \Mo{R}{G}$. On veut trouver un idempotent semi-simple $\ve$ tel que $\ve V\neq 0$. Il suffit bien-s\^ur de le faire pour $V$ irr\'eductible. Quitte \`a \'etendre les scalaires on peut supposer que $V$ est d\'efinie sur un corps alg\'ebriquement clos de caract\'eristique diff\'erente de $p$.

Donnons-nous  une strate $[\Lambda,n,r,\beta]$ semi-simple autoduale avec $r\leq n$  et un caract\`ere $\theta\in \CC(\Lambda,r,\beta)^\sigma$ tels que $\ve_\theta V\neq 0$. Remarquons que pour $r$ assez grand, il existe de telles donn\'ees. Choisissons
un prolongement $\wt\theta\in\CC(\Lambda,r-,\beta)$ de $\theta$. Nous noterons $A_-:=\{x\in A=\endo{F}{V}, x=-\o{x}\}$ et $\aG_\bullet(\Lambda)_-:=A_-\cap \aG_\bullet(\Lambda)$.
Comme dans le d\'ebut de la preuve du th\'eor\`eme 5.1 de \cite{Stevens}, il existe un \'el\'ement $c\in \aG_{-r}(\Lambda)_-$ tel que le caract\`ere $\vartheta:=\wt\theta {\psi_c}_{|H^r(\Lambda,\beta)^\sigma}$ apparaisse dans $V_{|H^r(\Lambda,\beta)^\sigma}$. 
D\'ecomposons  $F[\beta]=\prod_{i\in I_\beta} (E_i\times E_{-i}) \times \prod_{j\in J_\beta} E_j$ en un produit de corps o\`u l'involution $x\mapsto \o{x}$ identifie $E_i$ et $E_{-i}$ et stabilise chaque $E_j$. On a aussi la d\'ecomposition  orthogonale $V=\bigoplus_{i\in I_\beta} (V_i\oplus V_{-i}) \bigoplus_{j\in J_\beta} V_j$ selon les idempotents primitifs de ce produit. 
Par le m\^eme argument que le lemme 5.2 et le paragraphe qui le suit dans \cite{Stevens}, on peut supposer que $c$ se d\'ecompose en $c=\bigoplus_{i\in I_\beta} (c_i\oplus \o{c_i}) \bigoplus_{j\in J_\beta} b_j$
avec $c_i\in \aG_{-r}(\Lambda_i)$ pour $i\in I_\beta$ et $c_j\in \aG_{-r}(\Lambda_j)_-$ pour $j\in J_\beta$.

Choisissons des corestrictions mod\'er\'ees $s_k : \endo{F}{V_k} \To{} \endo{E_k}{V_k}$ pour $k\in  I_\beta\sqcup J_\beta$. On d\'efinit comme dans \cite[(5.3)]{Stevens} la $F[\beta]$-strate d\'eriv\'ee (autoduale) $[\Lambda,r,r-,s(c)]$ comme la somme
$$ \bigoplus_{i\in I_\beta} ([\Lambda_i,r,r-,s_i(c_i)]\oplus [\Lambda_{-i},r,r-,s_i(c_{-i})]) \bigoplus_{j\in J_\beta} [\Lambda_j,r,r-,s_j(c_j)]. $$
L'\'el\'ement $s(c)$ est donc dans l'alg\`ebre $A_\beta\subset A$ centralisatrice de $\beta$ et v\'erifie $\o{s(c)}=-s(c)$.

\begin{lemme} \cite[Thm 4.4]{StPLMS} \label{ll1}
Il existe une $F[\beta]$-strate autoduale {\em semi-simple} $[\Lambda',r,r-,\alpha']$ dans $V$ telle que
$$ s(c) + \left(\aG_{-(r_-)}(\Lambda)\cap A_\beta\right) \subset \alpha' + \left(\aG_{-(r_-)}(\Lambda')\cap A_\beta\right). $$
\end{lemme}
\begin{proof} Cet \'enonc\'e est une ``somme" d'\'enonc\'es analogues pour chaque $i\in I_\beta$, $j\in J_\beta$. Nous traitons seulement le cas $j\in J_\beta$, les autres cas se traitant de la m\^eme mani\`ere mais sans les complications ``autoduales".
 D'apr\`es \cite[Prop 4.2]{StPLMS}, il existe une $E_j$-fontion r\'eseau autoduale dans $V_j$ telle que $\aG_{-(r_-)}(\Lambda_j)\subset \aG_{-(r_-)}(\Lambda'_j)$ et  $s_j(c_j)\in \aG_{-r}(\Lambda_j)_-$
soit ``de r\'eduction semisimple" (voir {\em loc. cit} pour le sens de cette expression).
Si cette r\'eduction est nulle, on a $s_j(c_j)\in \aG_{-(r_-)}(\Lambda_j')$ et on peut prendre $\alpha'_j=0$.
 Sinon, le r\'eel $r$ est n\'ecessairement de la forme $n/e(\Lambda'_j)$ et on peut suivre la proc\'edure de la preuve de \cite[Thm 4.4]{StPLMS}.
%\'Ecrivons le polyn\^ome caract\'eristique $\varphi(X)$ de la strate $[\Lambda',r,r-,s(c)]$ sous la forme $$\varphi(X)= \prod_{k} \left(\varphi_k(X)\o\varphi_k(\eta X)\right)^{a_k} \prod_l \varphi_l(X)^{a_l} \times X^a $$
%o\`u $\eta$ est un signe ({\em cf} ...), $\varphi_l(X)=\o\varphi_l(\eta X)$ pour tout $l$, et tous ces polyn\^omes sont deux \`a deux premiers entre eux.
\end{proof}

L'outil fondamental est le lemme suivant  qui n'est \'enonc\'e que pour les strates {\em gauches} dans \cite{Stevens}, mais dont la preuve s'\'etend aux strates {\em autoduales} :
\begin{lemme} \label{ll2} \cite[Lemma 5.4]{Stevens}  Soit $[\Lambda',r',r'_-,\alpha']$ une $F[\beta]$-strate autoduale dans $V$ telle que
$$ s(c) + \left(\aG_{-(r_-)}(\Lambda)\cap A_\beta\right) \subset \alpha' + \left(\aG_{-(r'_-)}(\Lambda')\cap A_\beta\right). $$
Alors il existe $\wt\theta'\in \CC(\Lambda',r'-,\beta)$ et $c'\in \aG_{-r'}(\Lambda')$ de corestriction  $s(c')=\alpha'$  tels que la repr\'esentation $V$ contienne le caract\`ere $\vartheta':= \theta'{\psi_{c'}}_{|H^{r'}(\Lambda',\beta)^\sigma}$. Si de plus $\alpha'=0$, alors on peut choisir $c'=0$.
\end{lemme}
\begin{proof}
Nous nous contenterons de remarquer que les arguments de type ``th\'eorie des repr\'esenta\-tions" de la preuve de Stevens (lemme 5.9 de \cite{Stevens}) concernent des repr\'esentations de pro-$p$-groupes et restent valables dans  notre situation o\`u $\CM$ est remplac\'e par un corps alg\'ebriquement clos de caract\'eristique diff\'erente de $p$. On pourrait aussi facilement remplacer ces arguments en prouvant 
$$ \ve_{\wt\theta\psi_c} \in RG \ve_{\wt\theta'\psi'_{c'}}RG$$
de mani\`ere analogue \`a la preuve de la proposition \ref{constess}.
\end{proof}

Appliquons ce lemme \`a la strate $[\Lambda',r,r-,\alpha']$ du lemme \ref{ll1}. D'apr\`es \cite[Lemma 3.5]{Stevens}, la strate $[\Lambda',n,r-,\beta+c']$ est \'equivalente \`a une strate semi-simple, disons $[\Lambda',n,r-,\beta']$ et par \cite[Rk 3.14.ii)]{Stevens}, on a $\vartheta'\in \CC(\Lambda',r-,\beta)$.

On a donc une proc\'edure pour baisser {\em strictement} le ``niveau" d'un caract\`ere semi-simple auto\-dual intervenant dans $V$. Cependant, il n'est pas encore clair que cette proc\'edure produise apr\`es plusieurs it\'erations un caract\`ere semi-simple autodual de ``niveau" $0$ ; on peut en effet supposer les $r$ rationnels, mais on ne contr\^ole pas les d\'enominateurs. Pour les contr\^oler, il faut appliquer le lemme \ref{ll2} dans la situation o\`u
\begin{enumerate}
        \item $\wt\theta$ intervient dans $V$ et $c=0$.
        \item $\Lambda'$ est une somme sur $i,j$ de fonctions-r\'eseaux {\em optimales} (au sens de Moy-Prasad, {\em cf} \cite[(4.3)]{StPLMS} dans le contexte pr\'esent), $r'\leq r$ est tel que $\aG_{-r'_-}(\Lambda')\supset \aG_{-r_-}(\Lambda)$ et $\alpha'=0$.
\end{enumerate}
Le lemme \ref{ll2} nous dit alors que $V$ contient un caract\`ere de $\CC(\Lambda',r'-,\beta)$, mais cette fois $\Lambda'$ est somme de strates optimales et sa p\'eriode est born\'ee par un entier d\'ependant seulement de $\dim_F(V)$. Ainsi les $r\in\QM$ qui sont des sauts pour de telles strates ont leur d\'enominateurs born\'es, et la proc\'edure de raffinement ci-dessus produit  bien un caract\`ere dans un certain $\CC(\Lambda,0,\beta)$ intervenant dans $V$.

\section{Groupes mod\'er\'es} \label{modere}

Dans cette section, le corps de base redevient $K$ et le groupe r\'eductif connexe $\GC$ est suppos\'e mod\'er\'ement ramifi\'e. Nous allons appliquer les r\'esultats g\'en\'eraux de la partie \ref{mod} aux caract\`eres g\'en\'eriques introduits par Yu dans \cite{Yutame}.

\alin{Groupes de Yu}
Suivant \cite[sec. 2]{Yutame}, un sous-groupe ferm\'e $\GC^0$ de $\GC$ est appel\'e {\em \levi tordu mod\'er\'e} si apr\`es extension des scalaires de $K$ \`a une extension mod\'er\'ement ramifi\'ee, il devient un sous-groupe de Levi. On sait alors --toujours de mani\`ere non-canonique, mais peu importe-- identifier $B(\GC^0,K)$ \`a un sous-ensemble de $B(\GC,K)$, et ce de telle sorte que pour un point $x\in B(\GC^0,K)$, on ait $G^0_x=G^0\cap G_x$ et $G^{0+}_x=G^0\cap G_x^+$.

\'Etant donn\'ee
une suite de Levi tordus mod\'er\'ement ramifi\'es  $\vec\GC:=\{\GC^0\subset \GC^1 \subset \cdots \subset  \GC^d:=\GC\}$, et une suite $\vec{r}=\{0\leq r_0\leq \cdots \leq r_d\}$, Yu d\'efinit dans  \cite[sec 2]{Yutame} un groupe $\vec{G}_{x,\vec{r}}$ qu'il r\'ealise dans  dans \cite[sec 10]{Yumodel} comme groupe des points entiers d'un mod\`ele lisse $\vec\GC_{x,\vec{r}}$ de $\GC$ sur $\OC_K$.

D'apr\`es \cite[Prop. 10.4]{Yumodel}, la fibre sp\'eciale de $\vec\GC_{x,\vec{r}}$ est unipotente si $r_0>0$. Dans le cas contraire $r_0=0$, il y a une immersion ferm\'ee
 $\GC^0_x \injo \vec\GC_{x,\vec{r}}$ qui induit sur les fibres sp\'eciales un isomorphisme des quotients r\'eductifs. Comme dans la preuve de la proposition \ref{p1}, on en d\'eduit qu'un \levi $\MC$ de $\GC$ est $\vec\GC_{x,\vec{r}}^\circ$-admissible \ssi\  i) $x\in B(\MC,K)$, ii) $\MC$ contient le centre connexe $\ZC(\GC^0)^\circ$ \'eventuellement conjugu\'e par un \'el\'ement de $G_{x,\vec{r}}$. On en d\'eduit aussi que $\vec\GC_{x,\vec{r}}^\dag(\OC_K)=G_{x,\vec{r}}^+:=G_{x,\vec{r}}\cap G_x^+$.

\alin{Caract\`eres g\'en\'eriques}
Gardons les notations pr\'ec\'edentes et donnons-nous aussi une suite $\vec\phi:=\{\phi_0,\cdots,\phi_d\}$ o\`u chaque $\phi_i$ est un caract\`ere de $G^i:=\GC^i(F)$ qu'on supposera $G^{i+1}$-g\'en\'erique (pour $i<d$), au sens de \cite[sec. 5]{Yutame}.
On suppose que la suite $\vec{r}$ des niveaux $r_i$ des $\phi_i$ v\'erifie $0 <r_0< \cdots <r_{d-1} \leq r_d$,  on pose $s_i:=r_i/2$ et on note $\vec{s}:=(0,s_0,\cdots, s_{d-1})$ et $\vec{s}+:=(0+,s_0+,\cdots, s_{d-1}+)$.

 Selon \cite[Prop 4.1]{Yutame}, la donn\'ee de $\vec{\phi}$ d\'efinit 
 un caract\`ere $\theta=\prod_i\hat{\phi_i}$  de $\vec{G}_{x,\vec{s}+}$, normalis\'e par $\vec{G}_{x,\vec{s}}$, et tel que
\begin{enumerate}
\item $\hbox{Int}_{G_x}(\theta)=\vec{G}_{x,\vec{s}}$, {\em cf} \cite[Prop 4.1]{Yutame}
\item la forme bilin\'eaire $(x,y)\mapsto \theta([x,y])$ sur $\vec{G}_{x,\vec{s}}^+/\vec{G}_{x,\vec{s}+}$ est non-d\'eg\'en\'er\'ee ; c'est la somme de $i=0$ \`a $d-1$ des assertions de non-d\'eg\'en\'erescence de \cite[11.1]{Yutame} appliqu\'ees \`a $(G^i,G^{i+1})_{r_i,s_i}$ et $\hat{\phi_i}$.
\end{enumerate}
Soit alors $\MC$ un \levi $\vec\GC_{x,\vec{r}}^\circ$-admissible contenant $\ZC(\GC^0)^\circ$, et $(\PC,\o\PC)$ une paire de paraboliques oppos\'es de Levi commun $\MC$. Il r\'esulte imm\'ediatement des d\'efinitions qu'on a une d\'ecomposition d'Iwahori  $\vec{G}_{x,\vec{s}+} = \vec{U}_{x,\vec{s}+}\vec{M}_{x,\vec{s}+}\vec{\o{U}}_{x,\vec{s}+}$ pour laquelle les restrictions $\theta_{|\vec{U}_{x,\vec{s}+}}$ et $\theta_{|\vec{\o{U}}_{x,\vec{s}+}}$ sont triviales. De plus, le groupe $\vec{M}_{x,\vec{s}+}$ est le groupe de Yu associ\'e \`a la suite de Levis mod\'er\'es $\vec\MC:=\{\MC\cap\GC^0 \subset \cdots\subset \MC\cap\GC^d=\MC\}$ de $\MC$ et \`a la suite de ``r\'eels" $\vec{s}+$, et la restriction $\theta_M:=\theta_{|\vec{M}_{x,\vec{s}+}}$ est le caract\`ere associ\'e \`a la suite de caract\`eres g\'en\'eriques $\vec\phi_{|\vec\MC}$.

On a donc rassembl\'e tous les ingr\'edients pour prouver de la m\^eme mani\`ere que pour les groupes lin\'eaires et classiques  la proposition suivante qui est un analogue de \ref{p1} et \ref{p1c} :
\begin{prop} \label{p1m}
Gardons les notations ci-dessus, posons $R=\ZM_{p-cycl}[\frac{1}{p}]$, et notons  $\ve_{\theta_M}$ l'idempotent de $R\vec{M}_{x,\vec{s}}$ associ\'e \`a $\theta_M$.
Alors  
$$e_{U_x^+}\eubx \ve_{\theta_M} \in RG_x \eux \eubx \ve_{\theta_M} .$$
\end{prop}
\begin{proof}
Par ce qui pr\'ec\`ede et par \ref{constess}, $\ve_{\theta_M}$ est un idempotent essentiellement de niveau z\'ero pour le mod\`ele lisse $\vec\MC_{x,\vec{s}}^\circ$ de $\MC$. Par ce qui pr\'ec\`ede encore, on peut appliquer \ref{utile} avec $\GC'=\vec\GC_{x,\vec{s}}^\circ$ et $\GC=\GC_x^\circ$, et $\ve'=\ve_{\theta_M}$ et $\wt\ve'=\ve_\theta$.
On omet les d\'etails en renvoyant aux preuves de \ref{p1} et \ref{p1c}. 
\end{proof}

Remarquons maintenant que $\GC^0$ contient le centre connexe de $\MC$ et par cons\'equent, l'intersection $B(\GC^0,K)\cap B(\MC,K)$ est stable
par translations sous $a_M$. Appliquant la proposition pr\'ec\'edente aux points de $x+a_M$, on en d\'eduit que l'idempotent $\ve_{\theta_M}$ de $R\vec{M}_{x,\vec{r}}$ est $P$-bon au sens de \ref{pbon}. Il est fort probable qu'en utilisant des r\'esultats annonc\'es r\'ecemment par Ju-Lee Kim et Yu sur l'exhaustivit\'e de la construction de Yu pour les groupes mod\'er\'es (ceux dont tous les tores sont mod\'er\'ement ramifi\'es), on puisse prouver que la famille des idempotents du type $\ve_{\theta_M}$ comme ci-dessus est g\'en\'eratrice dans $\Mo{R}{M}$. Il faut pour cela attendre de lire les d\'etails de leur preuve.

\appendix

\section{D\'ecomposition ``par le niveau" de $\Mo{R}{G}$} \label{decomposition}

Le but de cette section est d'\'etendre \`a  la cat\'egorie
$\Mo{R}{G}$ la  d\'ecomposition ``par le niveau", implicite dans les travaux de Moy et
  Prasad lorsque $R=\CM$ et explicit\'ee dans le cas o\`u  $R$ est un corps par Vign\'eras dans
  \cite[II.5]{Vig}.

Les arguments reposent {\em in fine} sur les constructions de Moy et
Prasad dans \cite{MP2}, et notamment sur la comparaison entre deux
familles de filtrations concernant le groupe et l'alg\'ebre de
Lie. Pour cette comparaison,  des hypoth\`eses sont
n\'ecessaires, {\em cf} le commentaire qui suit  \cite[Cor
5.6.]{Yumodel}. Ces hypoth\`eses, peu contraignantes,  sont v\'erifi\'ees dans tous les cas
consid\'er\'es dans le pr\'esent article, et quoiqu'il en soit, Yu
explique dans \cite[5-6]{Yumodel} comment modifier la construction
originale de Moy-Prasad dans le cas g\'en\'eral.

\alin{D\'ecomposition de cat\'egories ab\'eliennes}
Nous rappelons ici un peu d'{\em abstract nonsense}.
Soit $\CC$ une cat\'egorie ab\'elienne (avec limites inductives exactes).
Pour une famille $(Q_n)_{n\in\NM}$ d'objets de $\CC$ on consid\`ere les propri\'et\'es suivantes :
\begin{itemize}
        \item (PROJ) Chaque $Q_n$ est projectif et de type fini (``compact").
        \item (DISJ) Si $n\neq m$, alors $\hom{Q_n}{Q_m}{\CC}=0$.
  \item (GEN) Pour tout objet $V$ de $\CC$, on a
    $\hom{\bigoplus_nQ_n}{V}{\CC}\neq 0$.
\end{itemize}
Par ailleurs, pour tout objet $V$ de
$\CC$, posons
$$V_n:= \sum_{\phi\in \hom{Q_n}{V}{G}} \im\phi \subseteq V,$$
un sous-objet de $V$. 
Les propri\'et\'es (PROJ) et (GEN) impliquent
que $V = \sum_n V_n$. La propri\'et\'e (DISJ), toujours
avec (PROJ),  assure que la somme est directe, {\em i.e.} 
$V=\bigoplus_n V_n$. On peut paraphraser cela en introduisant la 
sous-cat\'egorie pleine $\CC_n$ de $\CC$ form\'ee des  objets v\'erifiant $\hom{Q_m}{V}{\CC}=0$ pour tout $m\neq n$. On obtient en effet une d\'ecomposition de $\CC$ en une somme directe de sous-cat\'egories ``facteurs directs" 
$ \CC \simeq \bigoplus_n \CC_n$.

Plus g\'en\'eralement, pour $I\subset \NM$, notons $\CC_I$ la sous-cat\'egorie pleine de $\CC$
  form\'ee des  objets v\'erifiant $\hom{Q_m}{V}{\CC}=0$ pour $m\notin I$. 
Alors $\CC_I$ est une sous-cat\'egorie ``facteur direct" de $\CC$.

\alin{Types non raffin\'es de Moy-Prasad et d\'ecomposition de $\Mo{\zp}{G}$}
Soit $x\in B(\GC,K)$.
Moy et Prasad ont d\'efini (\cite{MP1}, \cite{MP2} et \cite[II.5]{Vig}), une certaine filtration d\'ecroissante de $G_x$ par des
pro-$p$-sous-groupes ouverts $G_{x,r}, r\in \RM_+$. Les sauts de cette filtration sont discrets et on a des relations de commutateurs $(G_{x,r},G_{x,s})\subset G_{x,r+s}$. Si l'on convient de noter $G_{x,r^+}:=\bigcup_{s>r} G_{x,s}$, alors $G_{x,0^+}=G_x^+$, et pour tout $r>0$, le groupe fini $G_{x,r}/G_{x,r^+}$ est naturellement un $\FM_p$-espace vectoriel.
Ils ont ensuite d\'efini certains caract\`eres complexes des gradu\'es $G_{x,r}/G_{x,r^+}$ appel\'es {\em types non raffin\'es minimaux de niveau $r$}, dont nous noterons l'ensemble $NR_{x,r}$. Ces caract\`eres sont donc \`a valeurs dans l'extension $\ZM[\frac{1}{p},\zeta_p]$ si $\zeta_p$ est une racine $p$-i\`eme de l'unit\'e.
Enfin, Moy et Prasad ont d\'efini un ensemble $PO$ de ``points optimaux" dans l'immeuble,  fini modulo action de $G$, et nous noterons $(r_n)_{n\in\NM}$ une \'enum\'eration des sauts des filtrations associ\'ees aux points de $PO$.

Posons maintenant $Q_0 : = \bigoplus_{x} \cind{G_{x,0^+}}{G}{\zp}$ o\`u $x$ d\'ecrit un
ensemble (fini) de repr\'esentants des $G$-orbites de sommets de $\IC$.
Pour $r\in \RM_+$, posons (comme dans la remarque de \cite[p. 136]{Vig}) 
$$ P(r) := \bigoplus_{x\in PO,\chi\in NR_{x,r}} \cind{G_{x,r}}{G}{\chi} $$
que l'on voit comme une
repr\'esentation de type fini  \`a coefficients dans $\ZM[\frac{1}{p}]$. % Puisqu'un type non raffin\'e minimal tordu par l'action d'un \'el\'ement de $\gal(\QM(\zeta_p)|\QM)$ est encore un type non raffin\'e minimal, cette repr\'esentation %est invariante sous $\gal(\QM(\zeta_p)|\QM)$ et 
%se descend  en une $\zp$-repr\'esentation que nous noterons encore $P(r)$.

\begin{lemme}
La famille $Q_n:=P(r_n)$, $n\in\NM$ d'objets de $\Mo{\zp}{G}$ v\'erifie les propri\'et\'es (PROJ), (GEN) et (DISJ).
\end{lemme}
\begin{proof}
D'apr\`es \cite[II.5]{Vig}, pour tout corps alg\'ebriquement clos $R$ de caract\'eristique $\neq p$, la famille de repr\'esentations $(Q_n\otimes_{\zp} R)_{n\in\NM}$ de $\Mo{R}{G}$ v\'erifie les propri\'et\'es (PROJ), (DISJ) \cite[II.5.8]{Vig} et (GEN) \cite[II.5.3]{Vig} de la section pr\'ec\'edente. Nous allons montrer que cela implique formellement qu'il en est de m\^eme de la famille $(Q_n)_{n\in\NM}$ dans $\Mo{\zp}{G}$.

(PROJ) :  En tant que somme d'induites de $\zp$-repr\'esentations de type fini  de pro-$p$-sous-groupes ouverts, $P(r)$ est  
projective et de type fini dans $\Mo{\zp}{G}$. 

(DISJ) :  puisque $Q_m$ est sans torsion,
$\hom{Q_n}{Q_m}{G} \injo \hom{Q_n\otimes \CM}{Q_m\otimes \CM}{G}$. Ce
dernier est nul par \cite[II.5.8]{Vig} appliqu\'e \`a $R=\CM$.

(GEN):  soit $V$ un objet de $\Mo{\zp}{G}$ tel qu'il existe $l\neq
p$ premier tel que $V_l:=\{v\in V, lv=0\} \neq 0$. On peut voir $V_l$
comme une $\FM_l$-repr\'esentation de $G$. 
 On sait alors par
\cite[II.5.3]{Vig} qu'il existe $n\in \NM$ et un morphisme non nul
$\phi :\; Q_n \To{} V_l\otimes \o\FM_l$. Par engendrement fini de $Q_n$,
ce morphisme se factorise par $V_l\otimes \FM_{l^k}$ pour un certain
$k\in \NM$. D'o\`u un morphisme non nul $Q_n \To{} (V_l)^k$
et par suite l'existence d'un morphisme  non nul $Q_n \To{} V_l$
que l'on peut composer avec l'injection $V_l\injo V$.

Si maintenant $V_l=0$ pour tout $l\neq p$, c'est \`a dire si $V$ n'a
pas de torsion, $V$ se plonge dans $V\otimes \QM$. Comme
pr\'ec\'edemment on d\'eduit de \cite[II.5.3]{Vig} l'existence d'un
morphisme non nul $Q_n \To{} V\otimes \QM$. Par engendrement fini de
$Q_n$, on peut multiplier par un ``d\'enominateur commun'' pour
obtenir un morphisme \`a image dans $V$. 

\end{proof} 

Bien-s\^ur on obtient des d\'ecompositions similaires en \'etendant
les scalaires \`a toute $\zp$-alg\`ebre $R$. On obtient aussi la
d\'ecomposition annonc\'ee dans la preuve de \ref{niveau0}. Enfin, on
d\'eduit de cette d\'ecomposition que si une repr\'esentation est
engendr\'ee par ses invariants sous un sous-groupe ouvert compact,
alors tous ses sous-objets ont la m\^eme propri\'et\'e. Ceci justifie
le corollaire \ref{corodecomposition}.

\end{document}